\pdfoutput=1
\RequirePackage{silence}
\documentclass[a4paper,svgnames]{amsart}
\usepackage[hmarginratio={1:1},vmarginratio={1:1},lmargin=65.0pt,tmargin=60.0pt]{geometry}

\allowdisplaybreaks

\usepackage[T1]{fontenc}
\usepackage{latexsym,exscale,mathtools,textcomp,bbm,aliascnt}
\usepackage{amsmath,amsthm,amsfonts,amssymb,amscd}
\usepackage{mathrsfs,stackrel}
\usepackage{tikz,tikz-cd}
\usetikzlibrary{matrix,arrows.meta,decorations.markings,shapes.multipart,decorations.pathreplacing,decorations.pathmorphing}
\usepackage{enumitem}
\setlist[enumerate]{itemsep=0.15cm,label=\emph{\upshape(\alph*)}}
\setlist[enumerate,2]{itemsep=0.15cm,label=\emph{\upshape(\roman*)}}
\usepackage{ytableau}
\usepackage{xparse}
\usepackage{cite}

\usepackage{array}
\newcolumntype{C}{>{$}c<{$}}
\usepackage{xltabular,booktabs}

\usepackage{tabulary}
\usepackage{booktabs}
\newcommand{\mystrut}{\rule[-0.2\baselineskip]{0pt}{1.1\baselineskip}}

\definecolor{mygray}{gray}{0.6}
\definecolor{mygraydark}{gray}{0.4}
\definecolor{mygraylight}{gray}{0.85}
\definecolor{spinach}{RGB}{46,139,87}
\definecolor{tomato}{RGB}{255,99,71}
\definecolor{orchid}{RGB}{143,40,194}
\definecolor{neon}{RGB}{77,77,255}
\definecolor{pumpkin}{RGB}{224,180,80}
\definecolor{citron}{RGB}{190,180,90}

\definecolor{lava}{RGB}{207,16,32}
\definecolor{cream}{RGB}{255,253,208}
\definecolor{verdigris}{RGB}{67,179,174}
\definecolor{Black}{RGB}{0,0,0}
\definecolor{mydarkblue}{RGB}{10,10,170}
\definecolor{darkspinach}{RGB}{20,70,20}
\definecolor{darktomato}{RGB}{155,40,30}
\definecolor{darkorchid}{RGB}{50,10,100}
\definecolor{darklava}{RGB}{150,8,16}

\let\emph\relax
\DeclareTextFontCommand{\emph}{\bfseries\em}

\newcommand{\leftsquigarrow}{\,\reflectbox{$\rightsquigarrow$}\,}
\newcommand{\placeholder}{{}_{-}}

\newcommand{\ie}{\text{i.e.}}
\newcommand{\eg}{\text{e.g.}}
\newcommand{\cf}{\text{cf.}}
\newcommand{\etc}{\text{etc.}}

\newcommand{\muta}{\text{mutatis mutandis}}

\renewcommand{\dots}{\text{...}}

\usepackage{todonotes}

\newcommand\floor[2]{\lfloor\tfrac{#1}{#2}\rfloor}

\newcommand\bsig{{\boldsymbol\sigma}}

\newcommand\blam{{\boldsymbol\lambda}}

\newcommand\brho{{\boldsymbol\rho}}
\newcommand\bmu{{\boldsymbol\mu}}
\newcommand\bnu{{\boldsymbol\nu}}
\newcommand\bom{{\boldsymbol\omega}}
\newcommand\charge{{\boldsymbol\kappa}}
\newcommand\chargetwo{{\boldsymbol\nu}}

\newcommand\ba{{\boldsymbol{a}}}

\newcommand\bfi{\boldsymbol{f}}
\newcommand\bu{{\boldsymbol{u}}}
\newcommand\bx{\boldsymbol{x}}
\newcommand\by{\boldsymbol{y}}

\newcommand\bS{{\boldsymbol{S}}}
\newcommand\bT{{\boldsymbol{T}}}

\newcommand\bs{\boldsymbol{s}}
\newcommand\bt{\boldsymbol{t}}

\newcommand\Pcal{\mathcal{P}}

\newcommand\Sub[1]{\dot{#1}}

\newcommand\Web{\mathbb{W}}

\NewDocumentCommand\Webaa{ O{\bx,\bi} O{\by,\bj} }{\Web_{#1}^{#2}}
\NewDocumentCommand\Webab{ O{(\charge,\bx),\bi} O{(\chargetwo,\by),\bj} }{\Web_{#1}^{#2}}
\NewDocumentCommand\Webabs{ O{\bx} O{\by} O{\bi} O{\bj} }{\Web_{\Sub{#1},\Sub{#3}}^{\Sub{#2},\Sub{#4}}}

\NewDocumentCommand\WA{ O{\beta} O{\brho} }{\mathscr{W}_{#1}^{#2}}
\NewDocumentCommand\WAc{ O{\beta} O{\brho} }{\mathscr{R}_{#1}^{#2}}
\NewDocumentCommand\WAs{ O{\beta} O{\brho} }{\mathscr{W}_{\Sub{#1}}^{{#2}}}
\NewDocumentCommand\WAsc{ O{\beta} O{\brho} }{\mathscr{R}^{\Sub{#1}}^{\Sub{#2}}}
\newcommand\WABasis{\mathcal{B}_{\beta}}
\NewDocumentCommand\SA{ O{\beta} O{\brho} }{\dot{\mathscr{W}}_{#1}^{#2}}

\NewDocumentCommand\TA{ O{\beta} O{\brho} }{\mathcal{W}_{#1}^{#2}}
\NewDocumentCommand\TAc{ O{\beta} O{\brho} }{\mathcal{R}_{#1}^{#2}}

\NewDocumentCommand\WAlam{ s O{\blam} } {\mathscr{W}_{n}^{\IfBooleanTF{#1}{\gdom}{\gedom}#2}}

\newcommand\idem[1][\blam]{\1_{#1}}
\newcommand\dotidem[1][\blam]{\1_{#1}^{y}}
\newcommand\sandwich[3]{\mathbb{S}_{#1}^{#2,#3}}
\newcommand\fsandwich[2]{\mathbb{S}_{#1}^{#2}}
\newcommand\daffine[1][\ba]{y^{#1}}
\newcommand\zeetwo[1][\bfi]{y^{#1}}

\newcommand\sand[1][\blam]{\mathbb{S}_{#1}}
\newcommand\sandbasis[1][\blam]{B_{#1}}

\newcommand\affine[1]{\underline{#1}}
\newcommand\Sym[1][n]{{\mathfrak{S}_{#1}}}
\newcommand\hell{\affine{\ell}}
\NewDocumentCommand\Parts{ D(){\ell} O{n} }{\mathrm{P}_{#1,#2}^{\brho}}
\newcommand\hParts[1][n]{\affine{\mathrm{P}}_{\hell,#1}^{\affine{\brho}}}
\newcommand\hPartsall[1][n]{\affine{\mathrm{P}}_{\hell,#1}^{\affine{\brho},all}}
\newcommand\Partsgood[1][n]{(\mathrm{P}_{\ell,#1}^{\brho})^{\neq 0}}
\newcommand\hPartsgood[1][n]{(\affine{\mathrm{P}}_{\hell,#1}^{\affine{\brho}})^{\neq 0}}

\newcommand\hcoord{\mathtt{x}_{\affine{\charge}}}
\newcommand\hcoorda{\mathtt{x}_{\affine{\charge}}^{A}}
\newcommand\hcoordf{\mathtt{max}_{\R}(\text{fin})}
\NewDocumentCommand\res{sO{}}{\mathop{\rm res}\nolimits_{#2\rho}}
\newcommand\Affch[1][\blam]{A^{y}(#1)}
\newcommand\Finch[1][\blam]{F^{y}(#1)}
\newcommand\Sandch[1][\blam]{S^{y}(#1)}
\newcommand\BX[1][{\WA[n](X)}]{\mathrm{D}_{#1}}
\newcommand\BXc[1][{\WAc[n](X)}]{\mathrm{D}_{#1}}
\newcommand\EX[1][{\TA[n](X)}]{\mathrm{E}_{#1}}
\newcommand\EXc[1][{\TAc[n](X)}]{\mathrm{E}_{#1}}
\newcommand\co{\mathsf{c}}
\newcommand\re{\mathsf{r}}

\def\Item(#1){\item\textbf{\upshape(#1\upshape)}}
\newcounter{relation}

\def\relationautorefname~#1\null{\upshape(W$_{#1}$\upshape)}
\let\ref\autoref
\let\eqref\autoref

\newcommand\Std{\mathop{\rm Std}\nolimits}
\newcommand\SStd[1][\charge]{\mathop{\rm SStd}\nolimits_{#1}}
\newcommand\hSStd[1][\underline{\charge}]{\mathop{\rm SStd}\nolimits_{#1}}

\newcommand\N{\mathbb{Z}_{\geq 0}}

\newcommand\Z{\mathbb{Z}}

\newcommand\typea[1][\Z]{A_{#1}}
\newcommand\typeb[1][\N]{B_{#1}}
\newcommand\typec[1][\N]{C_{#1}}
\newcommand\ac{AC}

\newcommand\aone[1][e]{A^{(1)}_{#1}}
\newcommand\cone[1][e]{C^{(1)}_{#1}}
\newcommand\atwo[1][e]{A^{(2)}_{2\cdot#1}}
\newcommand\dtwo[1][e]{D^{(2)}_{#1+1}}
\newcommand\bad{BAD}

\let\<=\langle
\let\>=\rangle

\NewDocumentCommand\mI{ sd() }{\mathtt{m}_{\IfBooleanT{#1}{\Sub}\Gamma}\IfNoValueF{#2}{(#2)}}
\def\pmod#1{\space(\text{mod }#1)}

\renewcommand\i{\hat i}
\newcommand\bi{\mathbf{i}}
\newcommand\bj{\mathbf{j}}

\newcommand\R{\mathbb{R}}

\def\map#1#2{\,{:}\,#1\!\longrightarrow\!#2}

\let\gedom=\trianglerighteq
\let\gdom=\vartriangleright

\usepackage{etoolbox}
\newcommand{\DeclareMyOperator}[1]{%
\expandafter\DeclareMathOperator\csname #1\endcsname{#1}
}
\forcsvlist{\DeclareMyOperator}{Shape,defect,gr,End}

\def\1{\mathbf{1}}

\def\Item(#1){\item\textbf{\upshape(#1\upshape)}}

\tikzset{
anchorbase/.style={baseline={([yshift=#1]current bounding box.center)}},
anchorbase/.default={-0.5ex},
dot colour/.initial=black,
dot colour/.default=black,
tinynodes/.style={font=\tiny,text height=0.25ex,text depth=0.05ex},
smallnodes/.style={font=\scriptsize,text height=0.75ex,text depth=0.15ex},
mor/.style={line width=0.75,color=black,fill=cream},
dots/.style={line width=1pt,line cap=round, gray, dash pattern=on 0pt off 2\pgflinewidth},
redstring/.style = {draw=red!50,fill=none,line width=0.35mm,preaction={draw=red,line width=2.5pt,-},nodes={color=red}},
affine/.style= {draw=citron!50,fill=none,
line width=0.35mm,preaction={draw=citron,line width=2.5pt,-},nodes={color=citron}},
solid/.style = {draw=blue,fill=none,dot colour=blue,line width=0.4mm,nodes={color=blue}},
ghost/.style = {draw=darkgray,fill=none,dot colour=darkgray,
densely dashed,line width=0.4mm,nodes={color=darkgray}},
dghost/.style = {draw=darkgray,double,fill=none,dot colour=darkgray,
densely dashed,line width=0.4mm,nodes={color=darkgray}},
crossline/.style={preaction={draw=white,line width=4.75pt,-},preaction={draw=black,line width=0.9pt,-}},
dot/.style = {
decoration={markings,
post length=0.25mm,
pre length=0.25mm,
mark=at position #1 with {\node[circle,radius=0.3cm,inner sep=-2.0pt,color=\pgfkeysvalueof{/tikz/dot colour},fill=\pgfkeysvalueof{/tikz/dot colour}]{};}
},
postaction={decorate}
},
dot/.default=0.5,
}

\tikzstyle directed=[postaction={decorate,decoration={markings,
mark=at position #1 with {\arrow[line width=0.25mm, black]{>}}}}]

\NewDocumentCommand\crossing{ O{1ex} mmO{} mmO{} } {
\tikz[centered=#1]{
\draw[#2, dot/.list={#4}](0,0)node[below]{$#3$}--++(1,1);
\draw[#5, dot/.list={#7}](1,0)node[below]{$#6$}--++(-1,1);
}%
}

\NewDocumentCommand\DottedIdempotentC{ D(){1.5} omm} {
\begin{tikzpicture}[anchorbase]
\def\residues{{#4}}
\foreach \res [count=\c,
evaluate=\res as \pos using {\res<=#3 ? #1*\res-\c*#1/15 : #1*(\res-2)-\c*#1/15}
] in {#4} {
\coordinate (\c) at (\pos,0);
\draw[solid](\c)node[below]{$\res$}--++(0,1);
\ifnum\res=#3\relax\else
\draw[ghost](\pos+#1,0)--++(0,1)node[above]{$\res$};
\fi
\ifnum\c=1\relax
\draw[redstring](\pos+#1/15,0)node[below]{$\res$}--++(0,1);
\fi
}
\foreach \pt [evaluate=\pt as \good using {\residues[\pt-1]==#3 ? 0 : 1}
] in {#2} {
\draw[solid,dot](\pt)--++(0,1);
\ifnum\good=1
\draw[ghost,dot]([shift={(#1,0)}]\pt)--++(0,1);
\fi
}
\end{tikzpicture}
}

\NewDocumentCommand\DottedIdempotentAA{ s D(){1.5} O{} mm} {
\begin{tikzpicture}[anchorbase]
\def\residues{{#5}}
\foreach \res [count=\c,
evaluate=\res as \pos using {\res<=#4 ? #2*\res-\c*#2/15 : #2*(\res-2)-\c*#2/15}
] in {#5} {
\coordinate (\c) at (\pos,0);
\draw[solid](\c)node[below]{$\res$}--++(0,1);
\ifnum\res=#4\relax\else
\draw[ghost](\pos+#2,0)--++(0,1)node[above]{$\res$};
\fi
\ifnum\c=1\relax
\IfBooleanTF{#1}{\draw[affine](\pos+#2/15,0)node[below]{$\res$}--++(0,1);}
{\draw[redstring](\pos+#2/15,0)node[below]{$\res$}--++(0,1);}
\fi
}
\foreach \pt [evaluate=\pt as \good using {\residues[\pt-1]==#4 ? 0 : 1}
] in {#3} {
\draw[solid,dot](\pt)--++(0,1);
\ifnum\good=1
\draw[ghost,dot]([shift={(#2,0)}]\pt)--++(0,1);
\fi
}
\end{tikzpicture}
}

\NewDocumentCommand\DottedIdempotentD{ s D(){3} O{} mm} {
\begin{tikzpicture}[anchorbase]
\def\residues{{#5}}
\foreach \res [count=\c, evaluate=\res as \pos using {\res==0 ? #2*\res-\c*#2/15 : #2*(\res-2)-\c*#2/15} ] in {#5} {
\coordinate (\c) at (\pos,0);
\draw[solid](\c)node[below]{$\res$}--++(0,1);
\ifnum\res=#4\relax\else
\ifnum\res=1
\draw[double,ghost](\pos+#2,0)--++(0,1)node[above]{$\res$};
\else\ifnum\res>0
\draw[ghost](\pos+#2,0)--++(0,1)node[above]{$\res$};
\fi
\fi
\fi
\ifnum\c=1\relax
\IfBooleanTF{#1}{\draw[affine](\pos+#2/15,0)node[below]{$\res$}--++(0,1);}
{\draw[redstring](\pos+#2/15,0)node[below]{$\res$}--++(0,1);}
\fi
}
\foreach \pt [evaluate=\pt as \res using {\residues[\pt-1]}] in {#3} {
\draw[solid,dot](\pt)--++(0,1);
\ifnum\res>1
\ifnum\res<#4 \draw[ghost,dot]([shift={(#2,0)}]\pt)--++(0,1);\fi
\fi
}
\end{tikzpicture}
}

\NewDocumentCommand\DoubleCrossing{ O{1ex} mmO{} mmO{} } {
\begin{tikzpicture}[centered=#1]
\draw[#2, dot/.list={#4}] (0,0)node[below]{$#3$} .. controls (1,0.5) .. (0,1);
\draw[#5, dot/.list={#7}] (1,0)node[below]{$#6$} .. controls (0,0.5) .. (1,1);
\end{tikzpicture}\bigskip
}

\NewDocumentCommand\dotstring{ O{1ex} mm}{%
\tikz[centered=#1]{\draw[#2,dot](0,0)node[below]{$#3$}--++(0,1);}}

\NewDocumentCommand\stringdot{ O{1ex} mmmm}{%
\tikz[centered=#1]{
\draw[#2](0,0)node[below]{$#3$}--++(0,1);
\draw[#4,dot](1,0)node[below]{$#5$}--++(0,1);
}%
}%

\DeclarePairedDelimiterX{\set}[1]{\{}{\}}{\setargs{#1}}
\NewDocumentCommand{\setargs}{>{\SplitArgument{1}{|}}m}{\setargsaux#1}
\NewDocumentCommand{\setargsaux}{mm}
{\IfNoValueTF{#2}{#1} {#1\,\delimsize|\,\mathopen{}#2}}

\def\NewTheorem#1{%
\newaliascnt{#1}{equation}%
\newtheorem{#1}[#1]{#1}%
\aliascntresetthe{#1}%
\expandafter\def\csname #1autorefname\endcsname{#1}%
}
\def\equationautorefname~#1\null{(#1)\null}

\numberwithin{equation}{subsection}

\NewTheorem{Proposition}
\NewTheorem{Theorem}
\NewTheorem{Corollary}
\AtEndEnvironment{Corollary}{\null\hfill$\square$}%
\NewTheorem{Lemma}
\theoremstyle{definition}
\NewTheorem{Definition}
\NewTheorem{Notation}
\NewTheorem{Example}
\AtEndEnvironment{Example}{\null\hfill$\Diamond$}%
\NewTheorem{Examples}
\AtEndEnvironment{Examples}{\vskip-10mm\null\hfill$\Diamond$}%

\theoremstyle{remark}
\NewTheorem{Remark}
\NewTheorem{Assumption}

\usepackage[hypertexnames=false]{hyperref}
\usepackage{bookmark}
\hypersetup{
pdftoolbar=true,
pdfmenubar=true,
pdffitwindow=false,
pdfstartview={FitH},
pdftitle={Cellularity for weighted KLRW algebras of types \texorpdfstring{$B$, $A^{(2)}$, $D^{(2)}$}{B,A2,D2}},
pdfauthor={Andrew Mathas and Daniel Tubbenhauer},
pdfsubject={},
pdfcreator={Andrew Mathas and Daniel Tubbenhauer},
pdfproducer={Andrew Mathas and Daniel Tubbenhauer},
pdfkeywords={},
pdfnewwindow=true,
colorlinks=true,
linkcolor=mydarkblue,
citecolor=teal,
filecolor=magenta,
urlcolor=orchid,
linkbordercolor=lava,
citebordercolor=teal,
urlbordercolor=orchid,
linktocpage=true
}
\newcommand{\nnfootnote}[1]{%
\begin{NoHyper}
\renewcommand\thefootnote{}\footnote{#1}%
\addtocounter{footnote}{-1}%
\end{NoHyper}
}

\def\makeautorefname#1#2{\csdef{#1autorefname}{#2}}
\makeautorefname{section}{Section}%
\makeautorefname{subsection}{Section}%
\makeautorefname{subsubsection}{Section}%

\begin{document}
\title[Cellularity for weighted KLRW algebras of types \texorpdfstring{$B$, $A^{(2)}$, $D^{(2)}$}{B,A2,D2}]{Cellularity for weighted KLRW algebras of types \texorpdfstring{$B$, $A^{(2)}$, $D^{(2)}$}{B,A2,D2}}
\author[A. Mathas and D. Tubbenhauer]{Andrew Mathas and Daniel Tubbenhauer}

\address{A.M.: The University of Sydney, School of Mathematics and Statistics F07, Office Carslaw 718, NSW 2006, Australia,
\href{http://www.maths.usyd.edu.au/u/mathas/}{www.maths.usyd.edu.au/u/mathas/}}
\email{andrew.mathas@sydney.edu.au}

\address{D.T.: The University of Sydney, School of Mathematics and Statistics F07, Office Carslaw 827, NSW 2006, Australia, \href{http://www.dtubbenhauer.com}{www.dtubbenhauer.com}}
\email{daniel.tubbenhauer@sydney.edu.au}

\begin{abstract}
This paper constructs homogeneous
affine sandwich cellular bases of weighted KLRW algebras in types $\typeb[\N]$, $\atwo[e]$, $\dtwo[e]$.
Our construction immediately gives homogeneous sandwich cellular bases for the finite dimensional quotients of these algebras.
Since weighted KLRW algebras generalize
KLR algebras, we also obtain bases and cellularity results for the (infinite and finite dimensional) KLR algebras.
\end{abstract}

\nnfootnote{\textit{Mathematics Subject Classification 2020.} Primary: 16G99, 20C08; Secondary: 20C30, 20G43.}
\nnfootnote{\textit{Keywords.} KLR algebras, diagram algebras, cellular bases, Hecke and Schur algebras.}

\addtocontents{toc}{\protect\setcounter{tocdepth}{1}}

\maketitle

\tableofcontents

\ytableausetup{centertableaux,mathmode,boxsize=0.64cm}


\section{Introduction}\label{S:Introduction}


\emph{KLR algebras}, or quiver Hecke algebras,
are graded infinite dimensional
algebras attached to a quiver. These algebras arise in
categorification and categorical representation theory,
see {\eg} \cite{KhLa-cat-quantum-sln-first},
\cite{KhLa-cat-quantum-sln-second}, \cite{Ro-2-kac-moody} or \cite{Ro-quiver-hecke}.
They admit finite dimensional quotients, called
\emph{cyclotomic KLR algebras}.
The KLR algebras and their cyclotomic quotients
play crucial roles in modern representation theory
and have attracted a lot of attention since their introduction.

A crucial problem is to determine whether
these algebras are (graded)
cellular in the sense of \cite{GrLe-cellular}, or some variation
of cellularity such as affine cellular \cite{KoXi-affine-cellular}.
It turns out that they tend to be graded affine cellular in the infinite dimensional
case, see {\eg} \cite{KlLo-klr-affine-cellular-finite-type}, and
graded cellular
in the finite dimensional case, see {\eg} \cite{HuMa-klr-basis}.
This paper gives
explicit homogeneous (affine) cellular bases for these algebras. Moreover, we will see that there is a close relationship between the cellular bases of the infinite dimensional KLR algebras and their finite dimensional quotients, even though this is very hard to see directly in the KLR setting.

\emph{Weighted KLRW algebras} are generalizations of KLR algebras that were introduced by Webster; see
\cite{We-weighted-klr}, \cite{We-rouquier-dia-algebra}, \cite{Bo-many-cellular-structures} or \cite{MaTu-klrw-algebras}.
Like the KLR algebras, the weighted KLRW algebras admit finite dimensional quotients. The KLR algebras can be constructed as idempotent subalgebras of the weighted KLRW algebras. Some care must be taken when comparing the weighted KLRW and KLR algebras because, for example, the functor induced by idempotent truncation annihilates some of the simple modules.

In the \emph{$\boldsymbol{\ac}$ types}, that is, types $\typea$, $\aone$, $\typec$, $\cone$, we showed \cite{MaTu-klrw-algebras} that one of the key properties of the weighted KLRW algebras is that they have homogeneous affine cellular bases
constructed in the style of low dimensional topology, and these bases automatically descend to
the finite dimensional quotients.
As a result, we obtained homogeneous (affine) cellular bases for the corresponding
KLR algebras and their finite dimensional quotients. Even better, all of these
bases are defined over any commutative integral domain $R$ (for example $R=\Z$).

It is natural to ask whether these constructions can be extended to other types. In this paper we answer this question affirmatively
for types $\typeb[\N]$, $\atwo[e]$, $\dtwo[e]$: the \emph{$\boldsymbol{\bad}$ types}. That is,
we construct homogeneous affine sandwich cellular bases
for weighted KLRW algebras that descend to the finite dimensional
quotients. As in the $\boldsymbol{\ac}$ types, we obtain bases for the corresponding KLR algebras.
The sandwich part of our bases, corresponding to the finite dimensional quotients, is given by copies of the dual numbers $R[X]/(X^{2})$,
which explains the powers of $2$ that appear in the dimension formulas
of the finite dimensional KLR algebras. For completeness,
we note that the weighted KLRW algebras
and their finite dimensional quotients of $\bad$ types are actually
(affine) cellular;
see \autoref{R:SandwichMain3} for a more precise statement. As a consequence, we obtain the usual cellularity results, such as a construction of the graded simple modules and that the decomposition matrix is unitriangular.

This paper is a sequel to \cite{MaTu-klrw-algebras}, so we assume
some familiarity to the definitions and results of that paper.
In sequels to this paper we hope to discuss weighted KLRW algebras of other types and their simple modules.

\begin{Remark}
The colors
used in this paper are not essential
and all strings are distinguishable by their thickness and if
they are solid or dashed.
\end{Remark}

\noindent\textbf{Acknowledgments.}
Want to thank Chris Bowman for a helpful zoom discussion that made us realize
several questions related to (affine) sandwich cellular algebras.
Special thanks to the referee: their comments and remarks were enormously helpful and improved the readability of the paper.

Both authors were supported, in part, by the Australian Research Council.
In these COVID19 infested times, we thank the first author's office for hosting us for the
one hour where the bulk of the mathematics in this paper was discovered.


\section{Sandwich cellular algebras}\label{S:BasesSandwich}
\stepcounter{subsection}

\emph{Cellular algebras} were introduced by Graham--Lehrer \cite{GrLe-cellular}.
Over the years many generalizations were
discovered. For example, incorporating a grading \cite{HuMa-klr-basis},
sandwiching a polynomial ring \cite{KoXi-affine-cellular} or sandwiching any
(potentially noncommutative) algebra \cite{GuWi-almost-cellular}, \cite{TuVa-handlebody}, \cite{Tu-sandwich-cellular}. (This final reference also discusses the history of  cellularity.) In order, these variations of cellularity are called \emph{graded cellular}, \emph{affine cellular} and \emph{sandwich cellular} algebras.

\begin{Remark}\label{R:BasesSandwichAffine}
The affine cellular algebras of \cite{KoXi-affine-cellular} are special case of sandwich cellular algebras were the sandwiched algebras are polynomial
rings or quotients of a polynomial rings. As we will see,
for weighted KLRW algebras of $\bad$ types we only need affine
cellularity. As we reserve the term affine cellular for the cases where the sandwich algebra is a full polynomial ring, we refer to our bases as sandwich cellular bases.
\end{Remark}

The following is a slight reformulation of \cite[Section 2]{TuVa-handlebody}.

\begin{Definition}\label{D:SandwichCellularAlgebra}
Let $R$ be a commutative ring with a unit.
Let $A$ be a locally unital graded $R$-algebra.
A \emph{graded sandwich cell datum} for $A$ is a tuple $(\Pcal,S,T,\sand[],\sandbasis[],C,\deg)$, where:
\begin{itemize}

\item $\Pcal=(\Pcal,\leq)$ is a poset (the \emph{middle set}),

\item $S=\bigcup_{\lambda\in\Pcal}S(\lambda)$ and $T=\bigcup_{\lambda\in\Pcal}T(\lambda)$ are collections of finite sets
(the \emph{bottom/top sets}),

\item $\sand[]=\bigoplus_{\lambda\in\Pcal}\sand[\lambda]$ is a direct sum of
graded algebras $\sand[\lambda]$
(the \emph{sandwiched algebras}) such that $\sandbasis[\lambda]$ is a homogeneous
basis of $\sand[\lambda]$ (we write $\deg$ for the degree function on $\sand[\lambda]$),

\item $C\map{\coprod_{\lambda\in\Pcal}S(\lambda)\times\sandbasis[\lambda]\times T(\lambda)}{A};(S,b,T)\mapsto C_{ST}^{b}$ is an injective map
(the \emph{basis}),

\item $\deg\map{\coprod_{\lambda\in\Pcal}S(\lambda)\cup T(\lambda)}{\Z}$ is a function
(the \emph{degree}),

\end{itemize}
such that:
\begin{enumerate}[label=\upshape(AC${}_{\arabic*}$\upshape)]

\item For $b\in\sandbasis[\lambda]$, $S\in S(\lambda),T\in T(\lambda)$ and $b\in\sandbasis[\lambda]$, $\lambda\in\Pcal$ the element $C_{ST}^{b}$
is homogeneous of degree $\deg(S)+\deg(b)+\deg(T)$.

\item The set $\set{C^{b}_{ST}|\lambda\in\Pcal,S\in S(\lambda),T\in T(\lambda),b\in\sandbasis[\lambda]}$
is an $R$-basis of $A$.

\item For all $x\in A$ there exist scalars $r^{c}_{SU}=r^{c}_{SU}(x)\in R$, that do not on
on $T$ and $b$, such that
\begin{gather*}
xC_{ST}^{b}\equiv
\sum_{U\in T(\lambda),c\in\sandbasis[\lambda]}r^{c}_{SU}C_{UT}^{c}\pmod{A^{>\lambda}},
\end{gather*}
where $A^{>\lambda}$ is the $R$-submodule of $A$ spanned by
$\set{C^{c}_{UV}|\mu\in\Pcal,\mu>\lambda,U\in S(\mu),V\in T(\mu),\in\sandbasis[\mu]}$.

\item Let $A(\lambda)=A^{\geq\lambda}/A^{>\lambda}$, where
$A^{\geq\lambda}$ is the $R$-submodule of $A$ spanned by $\set{C^{c}_{UV}|\mu\in\Pcal,
\mu\geq\lambda,U\in S(\mu),V\in T(\mu),c\in\sandbasis[\mu]}$. Then there exist free graded right and left
$\sand[\lambda]$-modules $\Delta(\lambda)$ and $\nabla(\lambda)$, respectively, such that $A(\lambda)\cong\Delta(\lambda)\otimes_{\sand[\lambda]}\nabla(\lambda)$.

\end{enumerate}
The algebra $A$ is a graded \emph{sandwich cellular algebra} if it
has a graded sandwich cell datum.

Assume that $S(\lambda)=T(\lambda)$ for all $\lambda\in\Pcal$, and that there is an antiinvolution $(\placeholder)^{\star}\colon A\to A$ such that:
\begin{enumerate}[label=\upshape(AC${}_{5}$\upshape)]

\item We have $(C^{b}_{ST})^{\star}\equiv C_{TS}^{b}\pmod{A^{>\lambda}}$.

\end{enumerate}
In this case, we write $(\Pcal,S=T,\sand[],\sandbasis[], C,\deg,(\placeholder)^{\star})$ and call this datum \emph{involutive}.
\end{Definition}

\begin{Remark}\label{R:SandwichCellularAlgebra}
The picture for elements in $\set{C^{b}_{ST}|\lambda\in\Pcal,S\in S(\lambda),T\in T(\lambda),b\in\sandbasis[\lambda]}$ is:
\begin{gather*}
C_{ST}^{s}
\leftrightsquigarrow
\begin{tikzpicture}[anchorbase,scale=1]
\draw[line width=0.75,color=black,fill=cream] (0,-0.5) to (0.25,0) to (0.75,0) to (1,-0.5) to (0,-0.5);
\node at (0.5,-0.25){$T$};
\draw[line width=0.75,color=black,fill=cream] (0,1) to (0.25,0.5) to (0.75,0.5) to (1,1) to (0,1);
\node at (0.5,0.75){$S$};
\draw[line width=0.75,color=black,fill=cream] (0.25,0) to (0.25,0.5) to (0.75,0.5) to (0.75,0) to (0.25,0);
\node at (0.5,0.25){$b$};
\end{tikzpicture}
,\quad
\begin{aligned}
\begin{tikzpicture}[anchorbase,scale=1]
\draw[line width=0.75,color=black,fill=cream] (0,1) to (0.25,0.5) to (0.75,0.5) to (1,1) to (0,1);
\node at (0.5,0.75){$S$};
\end{tikzpicture}
&\text{ top part $S\in S(\lambda)$},
\\
\begin{tikzpicture}[anchorbase,scale=1]
\draw[white,ultra thin] (0,0) to (1,0);
\draw[line width=0.75,color=black,fill=cream] (0.25,0) to (0.25,0.5) to (0.75,0.5) to (0.75,0) to (0.25,0);
\node at (0.5,0.25){$b$};
\end{tikzpicture}
&\text{ the middle $b\in\sandbasis[\lambda]$,}
\\
\begin{tikzpicture}[anchorbase,scale=1]
\draw[line width=0.75,color=black,fill=cream] (0,-0.5) to (0.25,0) to (0.75,0) to (1,-0.5) to (0,-0.5);
\node at (0.5,-0.25){$T$};
\end{tikzpicture}
&\text{ bottom part $T\in T(\lambda)$.}
\end{aligned}
\end{gather*}
This picture has to be interpreted with care because, in general, we do not necessarily have such a factorization.
\end{Remark}

We use the following terminology for special cases:
\begin{enumerate}

\item A graded \emph{affine sandwich cell datum} for $A$ is a graded sandwich cell datum such that, for all $\lambda\in\Pcal$ and for some $n(\lambda)\in\N$, we have $\sand[\lambda]\cong R[X_{1},\dots,X_{n(\lambda)}]$.

\item A graded \emph{cell datum} for $A$ is a graded sandwich cell datum such
that $\sand[\lambda]\cong R$ for all $\lambda\in\Pcal$.

\end{enumerate}

The image of $C$ in $A$ is an \emph{homogeneous sandwich cellular basis} for $A$. Similarly, we refer to affine sandwich cellular bases {\etc}

The following \emph{Clifford--Munn--Ponizovski\u{\i} theorem} parametrizes the simple modules of sandwich cellular algebras. To state this result
we need some notions.

For each $\lambda\in\Pcal$ there exists a cell module
$\Delta(\lambda)$ and a cellular pairing $\phi^{\lambda}$ on $\Delta(\lambda)$, see \cite[Section 2B]{TuVa-handlebody}.
The pairing $\phi^{\lambda}$ is a symmetric bilinear form.
Let $\Pcal^{\neq 0}\subset\Pcal$ be the subset
of those $\lambda$ for which $\phi^{\lambda}$ is nonzero.
The illustration to keep in mind is derived from \cite[Section 2B]{TuVa-handlebody}, the details can be copied from \cite[Section 2.2]{KoXi-affine-cellular}:
\begin{gather*}
\phi^{\lambda}
\left(
\begin{tikzpicture}[anchorbase,scale=1]
\draw[mor] (0,0.5) to (0.25,0) to (0.75,0) to (1,0.5) to (0,0.5);
\node at (0.5,0.25){$U$};
\draw[mor] (0.25,-0.5) to (0.25,0) to (0.75,0) to (0.75,-0.5) to (0.25,-0.5);
\node at (0.5,-0.25){$b$};
\draw[mor] (0,0.5) to (0,1) to (1,1) to (1,0.5) to (0,0.5);
\node at (0.5,0.75){$a$};
\end{tikzpicture}
,
\begin{tikzpicture}[anchorbase,scale=1]
\draw[mor] (0,1) to (0.25,1.5) to (0.75,1.5) to (1,1) to (0,1);
\node at (0.5,1.25){$D^{\prime}$};
\draw[mor] (0.25,1.5) to (0.25,2) to (0.75,2) to (0.75,1.5) to (0.25,1.5);
\node at (0.5,1.75){$b^{\prime}$};
\end{tikzpicture}
\right)
\quad\rightsquigarrow\quad
\begin{tikzpicture}[anchorbase,scale=1]
\draw[mor] (0,0.5) to (0.25,0) to (0.75,0) to (1,0.5) to (0,0.5);
\node at (0.5,0.25){$U$};
\draw[mor] (0.25,-0.5) to (0.25,0) to (0.75,0) to (0.75,-0.5) to (0.25,-0.5);
\node at (0.5,-0.25){$b$};
\draw[mor] (0,0.5) to (0,1) to (1,1) to (1,0.5) to (0,0.5);
\node at (0.5,0.75){$a$};
\draw[mor] (0,1) to (0.25,1.5) to (0.75,1.5) to (1,1) to (0,1);
\node at (0.5,1.25){$D^{\prime}$};
\draw[mor] (0.25,1.5) to (0.25,2) to (0.75,2) to (0.75,1.5) to (0.25,1.5);
\node at (0.5,1.75){$b^{\prime}$};
\end{tikzpicture}
\quad\leftsquigarrow\quad
\phi^{\lambda}
\left(
\begin{tikzpicture}[anchorbase,scale=1]
\draw[mor] (0,0.5) to (0.25,0) to (0.75,0) to (1,0.5) to (0,0.5);
\node at (0.5,0.25){$U$};
\draw[mor] (0.25,-0.5) to (0.25,0) to (0.75,0) to (0.75,-0.5) to (0.25,-0.5);
\node at (0.5,-0.25){$b$};
\end{tikzpicture}
,
\begin{tikzpicture}[anchorbase,scale=1]
\draw[mor] (0,0.5) to (0,1) to (1,1) to (1,0.5) to (0,0.5);
\node at (0.5,0.75){$a$};
\draw[mor] (0,1) to (0.25,1.5) to (0.75,1.5) to (1,1) to (0,1);
\node at (0.5,1.25){$D^{\prime}$};
\draw[mor] (0.25,1.5) to (0.25,2) to (0.75,2) to (0.75,1.5) to (0.25,1.5);
\node at (0.5,1.75){$b^{\prime}$};
\end{tikzpicture}
\right)
.
\end{gather*}

\begin{Theorem}\label{T:SandwichMain}
Let $R$ be a field, and let $A$ be a graded sandwich cellular algebra.
\begin{enumerate}

\item All (graded) simple $A$-modules are uniquely associated to a $\lambda\in\Pcal^{\neq 0}$, called their \emph{apex}.

\item Assume that $\sand[\lambda]$ is unital Artinian or commutative.
For a fixed apex $\lambda\in\Pcal^{\neq 0}$ there is
a 1:1-correspondence
\begin{gather*}
\set{\text{simple $A$-modules with apex $\lambda$}}/{\cong}
\;\xleftrightarrow{1:1}
\set{\text{simple $\sand[\lambda]$-modules}}/{\cong}
.
\end{gather*}

\item Assume that $\sand[\lambda]$ is unital graded Artinian or commutative.
For a fixed apex $\lambda\in\Pcal^{\neq 0}$ there is
a 1:1-correspondence
\begin{gather*}
\set{\text{graded simple $A$-modules with apex $\lambda$}}/{\cong}
\;\xleftrightarrow{1:1}
\set{\text{graded simple $\sand[\lambda]$-modules}}/{\cong}
.
\end{gather*}

\end{enumerate}
\end{Theorem}

Even for noncommutative $\sand[\lambda]$ the assumption of being unital (graded) Artinian
can be avoided under certain conditions, see \cite[Section 2C]{TuVa-handlebody} for details.

\begin{proof}
The proof is not very different from the general theory
as in \cite{GrLe-cellular}, \cite{KoXi-affine-cellular}, \cite{AnStTu-cellular-tilting}, \cite{EhTu-relcell},
\cite{GuWi-almost-cellular}, \cite{TuVa-handlebody} or\cite{Tu-sandwich-cellular}.
In particular, the result is just a (graded) reformulation of
\cite[Theorem 3]{GuWi-almost-cellular} and \cite[Section 2C]{TuVa-handlebody}. Details are omitted.
\end{proof}

\begin{Remark}\label{R:SandwichMain}
The bijection in \autoref{T:SandwichMain} can be made explicit:
the simple $A$-modules of apex $\lambda\in\Pcal^{\neq 0}$
can be constructed as the simple heads of (adjustments of) the cell modules.
\end{Remark}

\begin{Remark}\label{R:SandwichMain2}
The formulation of \autoref{T:SandwichMain}
is strongly inspired by
Green's theory of cells (Green's relations) \cite{Gr-structure-semigroups},
and the Clifford--Munn--Ponizovski\u{\i} theorem of semigroup theory,
see {\eg} \cite{GaMaSt-irreps-semigroups}
for a modern formulation.
\end{Remark}

\begin{Remark}\label{R:SandwichMain3}
If $A$ is involutive and all sandwiched algebras have an (affine sandwich) cellular datum
compatible with the sandwich structure on $A$, then $A$ also has a
cell datum that can be constructed by refining the sandwich cell datum
on $A$ in an appropriate sense; see \cite[Proposition 2.9]{TuVa-handlebody}. However, refining the datum often makes
the natural sandwich datum more cumbersome with very little gain.
In our case, the sandwiched algebras are (copies of)
polynomial rings and dual numbers, which are affine sandwich cellular
algebras, so all of our algebras can be refined to an (affine sandwich) cell datum.
Since these the representation theory of these algebras is well-understood, refining the sandwich cell datum
for the weighted KLRW algebras in this paper appears to be unnecessary.
\end{Remark}

\begin{Notation}\label{N:SandwichDefinition}
From now on all of our algebras are assumed to be involutive, and we will omit
the use of this word. We have separated the involutive condition
\upshape(AC${}_{5}$\upshape) in \autoref{D:SandwichCellularAlgebra} from the other axioms because
being involutive is not necessary for
\autoref{T:SandwichMain} to hold.
\end{Notation}


\section{Reminders about weighted KLRW algebras}\label{S:Recollection}


We recall the basic constructions and statements regarding weighted
KLRW algebras.
As we assume some familiarity with \cite{MaTu-klrw-algebras},
we often give examples instead of formal definitions.
The details for the material in \autoref{SS:RecollectionDefinition}--\autoref{SS:RecollectionDiagrams} can be found in
\cite[Section 2]{MaTu-klrw-algebras} and for \autoref{SS:RecollectionFaithful} in
\cite[Section 3]{MaTu-klrw-algebras}.

\begin{Notation}\label{N:RecollectionReadingDiagrams}
The following conventions used throughout the paper.
\begin{enumerate}

\item We work over a commutative integral domain $R$, the ground ring.

\item Graded algebra or module
will always mean a $\Z$-graded algebra or module.

\item We use the same diagrammatic conventions as in \cite{MaTu-klrw-algebras}:
\begin{gather*}
E\circ D
=
\begin{tikzpicture}[anchorbase,smallnodes,rounded corners]
\node[rectangle,draw,minimum width=0.5cm,minimum height=0.5cm,ultra thick] at(0,0){\raisebox{-0.05cm}{$D$}};
\node[rectangle,draw,minimum width=0.5cm,minimum height=0.5cm,ultra thick] at(0,0.5){\raisebox{-0.05cm}{$E$}};
\end{tikzpicture}
.
\end{gather*}
In particular, left actions and left modules are given by acting from the top.
Modules will always be left modules.

\end{enumerate}
\end{Notation}


\subsection{Weighted KLRW algebras in a nutshell}\label{SS:RecollectionDefinition}


The \emph{weighted KLRW algebras} are diagram algebras consisting of \emph{weighted KLRW diagrams}
(diagrams for short). These diagrams have three types
of strings: \emph{solid}, \emph{ghost} and \emph{red strings}.
We will also draw \emph{affine red strings}, but these are only a visual aid and are not actually part of the diagrams.
All of these strings are labeled by \emph{residues}, and we will
illustrate these strings as
\begin{gather}\label{Eq:RecollectionStrings}
\text{solid}:
\begin{tikzpicture}[anchorbase,smallnodes,rounded corners]
\draw[solid] (0,0)node[below]{$i$} to (0,0.5)node[above,yshift=-1pt]{$\phantom{i}$};
\end{tikzpicture}
,\quad
\text{ghost}:
\begin{tikzpicture}[anchorbase,smallnodes,rounded corners]
\draw[ghost] (0,0)node[below]{$\phantom{i}$} to (0,0.5)node[above,yshift=-1pt]{$i$};
\end{tikzpicture}
,\quad
\text{red}:
\begin{tikzpicture}[anchorbase,smallnodes,rounded corners]
\draw[redstring] (0,0)node[below]{$i$} to (0,0.5)node[above,yshift=-1pt]{$\phantom{i}$};
\end{tikzpicture}
,\quad
\text{affine red}:
\begin{tikzpicture}[anchorbase,smallnodes,rounded corners]
\draw[affine] (0,0)node[below]{$i$} to (0,0.5)node[above,yshift=-1pt]{$\phantom{i}$};
\end{tikzpicture}
,
\end{gather}
where the label is under the solid and (affine) red strings and over the ghost strings.

An example of such a diagram is (the notation is explained below):
\begin{gather}\label{Eq:RecollectionRunningExample}
\begin{tikzpicture}[anchorbase]
\node[circle,inner sep=1.8pt,fill=DarkBlue] (0) at (0.87,-0.5){};
\node at (0.87,-0.75){$0$};
\node[circle,inner sep=1.8pt,fill=DarkBlue] (1) at (0,1){};
\node at (0,1.25){$1$};
\node[circle,inner sep=1.8pt,fill=DarkBlue] (2) at (-0.87,-0.5){};
\node at (-0.87,-0.75){$2$};
\draw[directed=0.5](1)--(0);
\draw[directed=0.5](1)--(2);
\draw[directed=0.5](0)--(2);
\node at (0,-0.75){$0.9$};
\node at (0.9,0.25){$0.75$};
\node at (-0.9,0.25){$1.25$};
\end{tikzpicture}
\rightsquigarrow
\begin{tikzpicture}[anchorbase,smallnodes,rounded corners]
\draw[ghost](1.9,1)node[above,yshift=-1pt]{$0$}--++(1,-1);
\draw[ghost](2.9,1)node[above,yshift=-1pt]{$0$}--++(-1,-1);
\draw[ghost](6.25,1)node[above,yshift=-1pt]{$1$}--++(-1,-1)node[below]{$\phantom{i}$};
\draw[ghost](6.75,1)node[above,yshift=-1pt]{$1$}--++(-1,-1)node[below]{$\phantom{i}$};
\draw[ghost](8.4,1)node[above,yshift=-1pt]{$0$}--++(-0.5,-0.5)--++(0.5,-0.5);
\draw[solid](1,1)node[above,yshift=-1pt]{$\phantom{i}$}--++(1,-1)node[below]{$0$};
\draw[solid](2,1)--++(-1,-1)node[below]{$0$};
\draw[solid](3.25,1)--++(0,-1)node[below]{$2$};
\draw[solid](3.5,1)node[above,yshift=-1pt]{$\phantom{i}$}--++(2,-1)node[below]{$2$};
\draw[solid](4.5,1)node[above,yshift=-1pt]{$\phantom{i}$}--++(-1,-1)node[below]{$2$};
\draw[solid](5.5,1)node[above,yshift=-1pt]{$\phantom{i}$}--++(-1,-1)node[below]{$1$};
\draw[solid](7,1)--++(0.5,-0.5)--++(-0.5,-0.5)node[below]{$2$};
\draw[solid](7.5,1)--++(-0.5,-0.5)--++(0.5,-0.5)node[below]{$0$};
\draw[redstring](2.25,1)--++(0,-1)node[below]{$0$};
\draw[redstring](4,1)--++(0,-1)node[below]{$2$};
\draw[affine](7.75,1)--++(0,-1)node[below]{$1$};
\draw[affine](8.25,1)--++(0,-1)node[below]{$1$};
\end{tikzpicture}
.
\end{gather}
We will explain the different features of this diagram below.

The weighted KLRW algebras depend on the following input:
\begin{enumerate}

\item An oriented quiver $\Gamma=(I,E)$ with countable vertex set $I$ and
countable edge set $E$ arising from a symmetrizable generalized Cartan matrix. (The choice
of orientation does not play an essential role, {\cf} \cite[Proposition 3A.1]{MaTu-klrw-algebras}.)
We let $e+1=\# I$ be the number of vertices, allowing $e=\infty$.
The most important edges for this paper are single
and double edges, written as $i\rightarrow j$ and $i\Rightarrow j$, respectively, which are relevant for the quivers in \autoref{Eq:BasesFunnyQuivers}.
If the multiplicity of an edge
does not play a role, then we write $i\rightsquigarrow j$ if there is an edge from $i$ to $j$.

The key examples the reader might want to keep in mind while reading below are the $\bad$ types listed in \autoref{Eq:BasesFunnyQuivers}.

\item Nonnegative integers $n,\ell\in\N$ and a tuple $\brho=(\rho_{1},\dots,\rho_{\ell})\in I^{\ell}$. These numbers give the number of solid strings,
the number of red strings
(called the \emph{level}) and the labels of the red strings, respectively.

For example, in \autoref{Eq:RecollectionRunningExample}
we have $n=8$, $\ell=2$ and $\brho=(0,2)$.

\item Let $u$ and $v$ be indeterminates over $R$. For $i,j\in I$ define \emph{$Q$-polynomials}
\begin{gather}\label{Eq:RecollectionQPoly}
Q_{ij}(u,v)=
\left\{
\begin{array}{CC@{\quad}|@{\quad}CC}
$u-v$ & if $i\rightarrow j$, & $v-u$ & if $i\leftarrow j$,
\\
$u-v^{2}$ & if $i\Rightarrow j$, & $v-u^{2}$ & if $i\Leftarrow j$,
\\
$0$ & if $i=j$, & $1$ & otherwise.
\\
\end{array}
\right.
\end{gather}
These are $Q$-polynomials as in
\cite[Section 3.2.3]{Ro-2-kac-moody}, \cite{Ro-quiver-hecke} or \cite[Section 2.1]{We-weighted-klr}. Recall also that $Q_{iji}(u,v,w)=\frac{Q_{ij}(u,v)-Q_{ij}(u,w)}{w-v}$.

\item Various types of data determining the positions of the strings.
That is, a \emph{solid positioning} $\bx=(x_{1},\dots,x_{n})\in\R^{n}$,
a \emph{ghost shift} $\bsig=(\sigma_{\epsilon}\in\R_{\neq 0})_{\epsilon\in E}$
which is a weighting of $\Gamma$ (illustrated as edge labels), a \emph{charge}
$\charge=(\kappa_{1},\dots,\kappa_{\ell})\in\R^{\ell}$ such that
$\kappa_{1}<\dots<\kappa_{\ell}$. The $\bx$, $\bsig$ and $\charge$ are
used to determine
the boundary points of solid, ghost and red strings in the diagrams, and are such that there are no overlapping strings.

The $n$-tuple $\bx$ gives the $x$-coordinates of the $n$ solid strings. For example, in \autoref{Eq:RecollectionRunningExample} above we have $\bx=(1,2,3.25,3.5,4.5,5.5,7,7.5)$). The ghost shift gives the distances between a solid string and its ghost, {\eg} in \autoref{Eq:RecollectionRunningExample}, the ghost shifts are $0.75$ and $1.25$ for the two edges starting a vertex $1$ and $0.9$ for other edge. Finally, $\charge$ gives the $x$-coordinates of the red strings and in \autoref{Eq:RecollectionRunningExample} we have $\charge=(0.75,4)$. In general, the $x$-coordinates of the strings at the top and bottom of the diagrams will be different.

\item A set $X$ of \emph{positions}, which is the set of possible $x$-coordinates for the strings.
We stress that the rank of the weighted KLRW algebra, which is the number of solid strings,
strongly depends on the choice of $X$.

\item Diagrams consisting of $n$ solid strings labeled by $i\in I$, their ghost strings (explained in the next bullet point), and $\ell$ red strings labeled by $\brho$. Solid and ghost strings can be decorated with dots.
These diagrams are such that their internal points have local
neighborhoods of the form
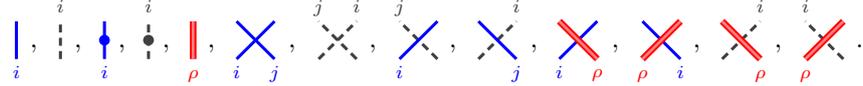
\begin{gather*}
\begin{tikzpicture}[anchorbase,smallnodes,rounded corners]
\draw[solid] (0,0)node[below]{$i$}node[below]{$\phantom{j}$} to (0,0.5)node[above,yshift=-1pt]{$\phantom{i}$};
\end{tikzpicture}
,\,
\begin{tikzpicture}[anchorbase,smallnodes,rounded corners]
\draw[ghost] (0,0)node[below]{$\phantom{i}$}node[below]{$\phantom{j}$} to (0,0.5)node[above,yshift=-1pt]{$i$};
\end{tikzpicture}
,\,
\begin{tikzpicture}[anchorbase,smallnodes,rounded corners]
\draw[solid,dot] (0,0)node[below]{$i$}node[below]{$\phantom{j}$} to (0,0.5)node[above,yshift=-1pt]{$\phantom{i}$};
\end{tikzpicture}
,\,
\begin{tikzpicture}[anchorbase,smallnodes,rounded corners]
\draw[ghost,dot] (0,0)node[below]{$\phantom{i}$}node[below]{$\phantom{j}$} to (0,0.5)node[above,yshift=-1pt]{$i$};
\end{tikzpicture}
,\,
\begin{tikzpicture}[anchorbase,smallnodes,rounded corners]
\draw[redstring] (0,0)node[below]{$\rho$}node[below]{$\phantom{j}$} to (0,0.5)node[above,yshift=-1pt]{$\phantom{i}$};
\end{tikzpicture}
,\,
\begin{tikzpicture}[anchorbase,smallnodes,rounded corners]
\draw[solid] (0,0)node[below]{$i$} to (0.5,0.5);
\draw[solid] (0.5,0)node[below]{$j$} to (0,0.5)node[above,yshift=-1pt]{$\phantom{i}$};
\end{tikzpicture}
,\,
\begin{tikzpicture}[anchorbase,smallnodes,rounded corners]
\draw[ghost] (0,0)node[below]{$\phantom{i}$} to (0.5,0.5)node[above,yshift=-1pt]{$i$};
\draw[ghost] (0.5,0)node[below]{$\phantom{i}$} to (0,0.5)node[above,yshift=-1pt]{$j$};
\end{tikzpicture}
,\,
\begin{tikzpicture}[anchorbase,smallnodes,rounded corners]
\draw[ghost] (0.5,0)node[below]{\phantom{i}} to (0,0.5)node[above,yshift=-1pt]{$j$};
\draw[solid] (0,0)node[below]{$i$} to (0.5,0.5);
\end{tikzpicture}
,\,
\begin{tikzpicture}[anchorbase,smallnodes,rounded corners]
\draw[ghost] (0,0)node[below]{\phantom{i}} to (0.5,0.5)node[above,yshift=-1pt]{$i$};
\draw[solid] (0.5,0)node[below]{$j$} to (0,0.5)node[above,yshift=-1pt]{$\phantom{i}$};
\end{tikzpicture}
,\,
\begin{tikzpicture}[anchorbase,smallnodes,rounded corners]
\draw[solid] (0,0)node[below]{$i$} to (0.5,0.5);
\draw[redstring] (0.5,0)node[below]{$\rho$} to (0,0.5)node[above,yshift=-1pt]{$\phantom{i}$};
\end{tikzpicture}
,\,
\begin{tikzpicture}[anchorbase,smallnodes,rounded corners]
\draw[solid] (0.5,0)node[below]{$i$} to (0,0.5)node[above,yshift=-1pt]{$\phantom{i}$};
\draw[redstring] (0,0)node[below]{$\rho$} to (0.5,0.5);
\end{tikzpicture}
,\,
\begin{tikzpicture}[anchorbase,smallnodes,rounded corners]
\draw[ghost] (0,0)node[below]{$\phantom{i}$} to (0.5,0.5)node[above,yshift=-1pt]{$i$};
\draw[redstring] (0.5,0)node[below]{$\rho$} to (0,0.5)node[above,yshift=-1pt]{$\phantom{i}$};
\end{tikzpicture}
,\,
\begin{tikzpicture}[anchorbase,smallnodes,rounded corners]
\draw[ghost] (0.5,0)node[below]{$\phantom{i}$} to (0,0.5)node[above,yshift=-1pt]{$i$};
\draw[redstring] (0,0)node[below]{$\rho$} to (0.5,0.5);
\end{tikzpicture}
.
\end{gather*}

\item Let $\epsilon\colon i\rightsquigarrow j$ be an edge. If $\sigma_{\epsilon}>0$ then each solid $i$-string has a ghost $i$-string that is shifted $\sigma_{\epsilon}$ units to the right. For each dot on the solid string there is a corresponding dot on the ghost string, which is shifted $\sigma_{\epsilon}$ units to the right. The ghost string is just a translate of the solid string, so it mimics its solid string. Similarly, if $\sigma_{\epsilon}<0$ then each solid $j$-string has a ghost $j$-string that is shifted $-\sigma_{\epsilon}$ units to the right, together with any corresponding dots.

In \autoref{Eq:RecollectionRunningExample}, each solid $0$-string has one ghost shifted 
$0.9$ units to the right.
The solid $1$-string has two ghosts, as there are two edges starting at $1$, that are shifted $0.75$
and $1.25$ unit to the right. The solid $2$-strings do not have any associated ghost strings, as no edges start from $2$.
Here are some more examples of these diagrams with a weighted quiver, where the edge weights gives the ghost shifts:
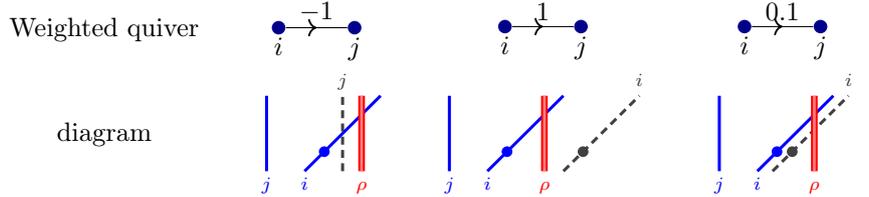
\begin{gather*}
\begin{array}{c@{\qquad}c@{\qquad}c@{\qquad}c}
\text{Weighted quiver}
&
\begin{tikzpicture}[anchorbase]
\node[circle,inner sep=1.8pt,fill=DarkBlue] (0) at (0,0){};
\node at (0,-0.25){$i$};
\node[circle,inner sep=1.8pt,fill=DarkBlue] (1) at (1,0){};
\node at (1,-0.28){$j$};
\draw[directed=0.5](0) to node[above,yshift=-1pt]{$-1$} node[below]{$\phantom{-1}$} (1);
\end{tikzpicture}
&
\begin{tikzpicture}[anchorbase]
\node[circle,inner sep=1.8pt,fill=DarkBlue] (0) at (0,0){};
\node at (0,-0.25){$i$};
\node[circle,inner sep=1.8pt,fill=DarkBlue] (1) at (1,0){};
\node at (1,-0.28){$j$};
\draw[directed=0.5](0) to node[above,yshift=-1pt]{$1$} node[below]{$\phantom{-1}$} (1);
\end{tikzpicture}
&
\begin{tikzpicture}[anchorbase]
\node[circle,inner sep=1.8pt,fill=DarkBlue] (0) at (0,0){};
\node at (0,-0.25){$i$};
\node[circle,inner sep=1.8pt,fill=DarkBlue] (1) at (1,0){};
\node at (1,-0.28){$j$};
\draw[directed=0.5](0) to node[above,yshift=-1pt]{$0.1$} node[below]{$\phantom{-1}$} (1);
\end{tikzpicture}
\\[0.25cm]
\text{diagram}
&
\begin{tikzpicture}[anchorbase,smallnodes,rounded corners]
\draw[ghost](0.5,0)node[below]{$\phantom{i}$}--++(0,1)node[above,yshift=-1pt]{$j$};
\draw[solid,dot=0.25](0,0)node[below]{$i$}--++(1,1);
\draw[solid](-0.5,0)node[below]{$j$}--++(0,1);
\draw[redstring](0.75,0)node[below]{$\rho$}--++(0,1);
\end{tikzpicture}
&
\begin{tikzpicture}[anchorbase,smallnodes,rounded corners]
\draw[ghost,dot=0.25](1,0)node[below]{$\phantom{i}$}--++ (1,1)node[above,yshift=-1pt]{$i$};
\draw[solid,dot=0.25](0,0)node[below]{$i$}--++(1,1);
\draw[solid](-0.5,0)node[below]{$j$}--++(0,1);
\draw[redstring](0.75,0)node[below]{$\rho$}--++(0,1);
\end{tikzpicture}
&
\begin{tikzpicture}[anchorbase,smallnodes,rounded corners]
\draw[ghost,dot=0.25](0.2,0)node[below]{$\phantom{i}$}--++ (1,1)node[above,yshift=-1pt]{$i$};
\draw[solid,dot=0.25](0,0)node[below]{$i$}--++(1,1);
\draw[solid](-0.5,0)node[below]{$j$}--++(0,1);
\draw[redstring](0.75,0)node[below]{$\rho$}--++(0,1);
\end{tikzpicture}
\end{array}
.
\end{gather*}

\item There is a degree function on these diagrams given by the local
rules
\begin{gather*}
\deg\begin{tikzpicture}[anchorbase,smallnodes,rounded corners]
\draw[solid,dot] (0,0)node[below]{$i$} to (0,0.5)node[above,yshift=-1pt]{$\phantom{i}$};
\end{tikzpicture}
=2d_{i}
,\quad
\deg\begin{tikzpicture}[anchorbase,smallnodes,rounded corners]
\draw[ghost,dot] (0,0)node[below]{$\phantom{i}$} to (0,0.5)node[above,yshift=-1pt]{$i$};
\end{tikzpicture}=0,
\quad
\deg\begin{tikzpicture}[anchorbase,smallnodes,rounded corners]
\draw[solid] (0,0)node[below]{$i$} to (0.5,0.5);
\draw[solid] (0.5,0)node[below]{$j$} to (0,0.5)node[above,yshift=-1pt]{$\phantom{i}$};
\end{tikzpicture}
=-\delta_{i,j}2d_{i}
,\quad
\deg\begin{tikzpicture}[anchorbase,smallnodes,rounded corners]
\draw[ghost] (0,0)node[below]{$\phantom{i}$} to (0.5,0.5)node[above,yshift=-1pt]{$i$};
\draw[solid] (0.5,0)node[below]{$j$} to (0,0.5)node[above,yshift=-1pt]{$\phantom{i}$};
\end{tikzpicture}
=
\deg\begin{tikzpicture}[anchorbase,smallnodes,rounded corners]
\draw[ghost] (0.5,0)node[below]{$\phantom{i}$} to (0,0.5)node[above,yshift=-1pt]{$i$};
\draw[solid] (0,0)node[below]{$j$} to (0.5,0.5)node[above,yshift=-1pt]{$\phantom{i}$};
\end{tikzpicture}
=
\begin{cases}
\<\alpha_{i},\alpha_{j}\>
&\text{if $i\rightsquigarrow j$},
\\
0&\text{else},
\end{cases}
\\
\deg\begin{tikzpicture}[anchorbase,smallnodes,rounded corners]
\draw[ghost] (0,0)node[below]{$\phantom{i}$} to (0.5,0.5)node[above,yshift=-1pt]{$i$};
\draw[ghost] (0.5,0)node[below]{$\phantom{i}$} to (0,0.5)node[above,yshift=-1pt]{$j$};
\end{tikzpicture}
=0
,\quad
\deg\begin{tikzpicture}[anchorbase,smallnodes,rounded corners]
\draw[solid] (0,0)node[below]{$i$} to (0.5,0.5)node[above,yshift=-1pt]{$\phantom{i}$};
\draw[redstring] (0.5,0)node[below]{$j$} to (0,0.5);
\end{tikzpicture}
=
\deg\begin{tikzpicture}[anchorbase,smallnodes,rounded corners]
\draw[solid] (0.5,0)node[below]{$i$} to (0,0.5)node[above,yshift=-1pt]{$\phantom{i}$};
\draw[redstring] (0,0)node[below]{$j$} to (0.5,0.5);
\end{tikzpicture}
=\tfrac{1}{2}\delta_{i,j}\<\alpha_{i},\alpha_{i}\>
,\quad
\deg\begin{tikzpicture}[anchorbase,smallnodes,rounded corners]
\draw[ghost] (0,0)node[below]{$\phantom{i}$} to (0.5,0.5)node[above,yshift=-1pt]{$i$};
\draw[redstring] (0.5,0)node[below]{$j$} to (0,0.5)node[above,yshift=-1pt]{$\phantom{i}$};
\end{tikzpicture}
=
\deg\begin{tikzpicture}[anchorbase,smallnodes,rounded corners]
\draw[ghost] (0.5,0)node[below]{$\phantom{i}$} to (0,0.5)node[above,yshift=-1pt]{$i$};
\draw[redstring] (0,0)node[below]{$j$} to (0.5,0.5)node[above,yshift=-1pt]{$\phantom{i}$};
\end{tikzpicture}
=0.
\end{gather*}
Here, $\<\placeholder,\placeholder\>$ is the Cartan pairing
associated to $\Gamma$, and $(d_{0},\dots,d_{e})$ is the symmetrizer.

The diagram in
\autoref{Eq:RecollectionRunningExample} is of degree $-1$ with nonzero contributions given by
$-2$ for the crossing with solid $0$-strings, $1$ twice for the solid $2$-strings crossing the red $2$-string, $-2$ for the crossing with solid $2$-strings, and $1$ for the crossing of a ghost $1$ and a solid $2$-string.

\item Let $\bu=(u_{1},\dots,u_{n})$.
If $f(\bu)=u_{1}^{a_{1}}\dots u_{n}^{a_{n}}\in\N[u_{1},\dots,u_{n}]$ is a monomial, then $f(\bu)$ times a diagram is the same diagram but with $a_{k}$ dots added at the top of the $k$th solid string, and all of its ghosts, for $1\leq k\leq n$. Similarly, if $f(\bu)\in\N[u_{1},\dots,u_{n}]$, then $f(\bu)$ times a diagram is the corresponding linear combination of diagrams.

\end{enumerate}

\begin{Definition}\label{D:RecollectionwKLRW}
The \emph{weighted KLRW algebra} $\WA[n](X)$
is the graded unital associative $R$-algebra generated by
such diagrams with multiplication given by stacking the diagrams
and subject to the multilocal relations listed below.
\emph{Multilocal} relations need to be 
applied in local neighborhoods around the solid strings
and, simultaneously, in the corresponding local neighborhoods around the ghost strings. See \autoref{E:ProofsMultilocal} below for an explicit example.

\begin{enumerate}

\item The \emph{dot sliding relations} hold. That is, solid and ghost dots can pass through any crossing except:
\begin{gather}\label{Eq:RecollectionDotCrossing}
\begin{tikzpicture}[anchorbase,smallnodes,rounded corners]
\draw[solid](0.5,0.5)node[above,yshift=-1pt]{$\phantom{i}$}--(0,0) node[below]{$i$};
\draw[solid,dot=0.25](0,0.5)--(0.5,0) node[below]{$i$};
\end{tikzpicture}
-
\begin{tikzpicture}[anchorbase,smallnodes,rounded corners]
\draw[solid](0.5,0.5)node[above,yshift=-1pt]{$\phantom{i}$}--(0,0) node[below]{$i$};
\draw[solid,dot=0.75](0,0.5)--(0.5,0) node[below]{$i$};
\end{tikzpicture}
=
\begin{tikzpicture}[anchorbase,smallnodes,rounded corners]
\draw[solid](0,0.5)node[above,yshift=-1pt]{$\phantom{i}$}--(0,0) node[below]{$i$};
\draw[solid](0.5,0.5)--(0.5,0) node[below]{$i$};
\end{tikzpicture}
=
\begin{tikzpicture}[anchorbase,smallnodes,rounded corners]
\draw[solid,dot=0.75](0.5,0.5)node[above,yshift=-1pt]{$\phantom{i}$}--(0,0) node[below]{$i$};
\draw[solid](0,0.5)--(0.5,0) node[below]{$i$};
\end{tikzpicture}
-
\begin{tikzpicture}[anchorbase,smallnodes,rounded corners]
\draw[solid,dot=0.25](0.5,0.5)node[above,yshift=-1pt]{$\phantom{i}$}--(0,0) node[below]{$i$};
\draw[solid](0,0.5)--(0.5,0) node[below]{$i$};
\end{tikzpicture}
.
\end{gather}
Note that this relation may involve ghost strings, depending on whether edges in the quiver start from $i$. Here, and below, we omit strings for clarity unless they make a nontrivial contribution to the diagram.

\item The \emph{Reidemeister II relations} hold except in the following cases:
\begin{gather}\label{Eq:RecollectionReidemeisterII}
\begin{gathered}
\begin{tikzpicture}[anchorbase,smallnodes,rounded corners]
\draw[solid](0,1)--++(0.5,-0.5)--++(-0.5,-0.5) node[below]{$i$};
\draw[solid](0.5,1)node[above,yshift=-1pt]{$\phantom{i}$}--++(-0.5,-0.5)--++(0.5,-0.5) node[below]{$i$};
\end{tikzpicture}
=0
,\quad
\begin{tikzpicture}[anchorbase,smallnodes,rounded corners]
\draw[ghost](0,1)node[above,yshift=-1pt]{$i$}--++(0.5,-0.5)--++(-0.5,-0.5) node[below]{$\phantom{i}$};
\draw[solid](0.5,1)--++(-0.5,-0.5)--++(0.5,-0.5) node[below]{$j$};
\end{tikzpicture}
=Q_{ij}(\bu)
\begin{tikzpicture}[anchorbase,smallnodes,rounded corners]
\draw[ghost](0,1)node[above,yshift=-1pt]{$i$}--++(0,-1)node[below]{$\phantom{i}$};
\draw[solid](0.5,1)--++(0,-1)node[below]{$j$};
\end{tikzpicture}
\;\text{or}\;
\begin{tikzpicture}[anchorbase,smallnodes,rounded corners]
\draw[ghost](0.5,1)node[above,yshift=-1pt]{$i$}--++(-0.5,-0.5)--++(0.5,-0.5) node[below]{$\phantom{i}$};
\draw[solid](0,1)--++(0.5,-0.5)--++(-0.5,-0.5) node[below]{$j$};
\end{tikzpicture}
=Q_{ji}(\bu)
\begin{tikzpicture}[anchorbase,smallnodes,rounded corners]
\draw[ghost](0.5,1)node[above,yshift=-1pt]{$i$}--++(0,-1)node[below]{$\phantom{i}$};
\draw[solid](0,1)--++(0,-1)node[below]{$j$};
\end{tikzpicture}
\quad
\text{if $i\rightsquigarrow j$}
,\\
\begin{tikzpicture}[anchorbase,smallnodes,rounded corners]
\draw[solid](0.5,1)node[above,yshift=-1pt]{$\phantom{i}$}--++(-0.5,-0.5)--++(0.5,-0.5) node[below]{$i$};
\draw[redstring](0,1)--++(0.5,-0.5)--++(-0.5,-0.5) node[below]{$i$};
\end{tikzpicture}
=
\begin{tikzpicture}[anchorbase,smallnodes,rounded corners]
\draw[solid,dot](0.5,0)node[below]{$i$}--++(0,1)node[above,yshift=-1pt]{$\phantom{i}$};
\draw[redstring](0,0)node[below]{$i$}--++(0,1);
\end{tikzpicture}
,\quad
\begin{tikzpicture}[anchorbase,smallnodes,rounded corners]
\draw[solid](0,1)node[above,yshift=-1pt]{$\phantom{i}$}--++(0.5,-0.5)--++(-0.5,-0.5) node[below]{$i$};
\draw[redstring](0.5,1)--++(-0.5,-0.5)--++(0.5,-0.5) node[below]{$i$};
\end{tikzpicture}
=
\begin{tikzpicture}[anchorbase,smallnodes,rounded corners]
\draw[solid,dot](0,0)node[below]{$i$}--++(0,1)node[above,yshift=-1pt]{$\phantom{i}$};
\draw[redstring](0.5,0)node[below]{$i$}--++(0,1);
\end{tikzpicture}.
\end{gathered}
\end{gather}

\item The \emph{Reidemeister III relations} hold except in the following cases:
\begin{gather}\label{Eq:RecollectionReidemeisterIII}
\begin{gathered}
\begin{tikzpicture}[anchorbase,smallnodes,rounded corners]
\draw[ghost](1,1)node[above,yshift=-1pt]{$i$}--++(1,-1)node[below]{$\phantom{i}$};
\draw[ghost](2,1)node[above,yshift=-1pt]{$i$}--++(-1,-1)node[below]{$\phantom{i}$};
\draw[solid,smallnodes,rounded corners](1.5,1)--++(-0.5,-0.5)--++(0.5,-0.5)node[below]{$j$};
\end{tikzpicture}
{=}
\begin{tikzpicture}[anchorbase,smallnodes,rounded corners]
\draw[ghost](3,1)node[above,yshift=-1pt]{$\phantom{i}$}--++(1,-1)node[below]{$\phantom{i}$};
\draw[ghost](4,1)node[above,yshift=-1pt]{$i$}--++(-1,-1)node[below]{$\phantom{i}$};
\draw[solid,smallnodes,rounded corners](3.5,1)--++(0.5,-0.5)--++(-0.5,-0.5)node[below]{$j$};
\end{tikzpicture}
-Q_{iji}(\bu)
\begin{tikzpicture}[anchorbase,smallnodes,rounded corners]
\draw[ghost](6.2,1)node[above,yshift=-1pt]{$i$}--++(0,-1)node[below]{$\phantom{i}$};
\draw[ghost](7.2,1)node[above,yshift=-1pt]{$i$}--++(0,-1)node[below]{$\phantom{i}$};
\draw[solid](6.7,1)--++(0,-1)node[below]{$j$};
\end{tikzpicture}
\;\text{or}\;
\begin{tikzpicture}[anchorbase,smallnodes,rounded corners]
\draw[solid](1,1)node[above,yshift=-1pt]{$\phantom{i}$}--++(1,-1)node[below]{$j$};
\draw[solid](2,1)--++(-1,-1)node[below]{$j$};
\draw[ghost,smallnodes,rounded corners](1.5,1)node[above,yshift=-1pt]{$i$}--++(-0.5,-0.5)--++(0.5,-0.5)node[below]{$\phantom{i}$};
\end{tikzpicture}
{=}
\begin{tikzpicture}[anchorbase,smallnodes,rounded corners]
\draw[solid](3,1)node[above,yshift=-1pt]{$\phantom{i}$}--++(1,-1)node[below]{$j$};
\draw[solid](4,1)--++(-1,-1)node[below]{$j$};
\draw[ghost,smallnodes,rounded corners](3.5,1)node[above,yshift=-1pt]{$i$}--++(0.5,-0.5)--++(-0.5,-0.5)node[below]{$\phantom{i}$};
\end{tikzpicture}
+Q_{iji}(\bu)
\begin{tikzpicture}[anchorbase,smallnodes,rounded corners]
\draw[solid](7.2,1)node[above,yshift=-1pt]{$\phantom{i}$}--++(0,-1)node[below]{$j$};
\draw[solid](8.2,1)--++(0,-1)node[below]{$j$};
\draw[ghost](7.7,1)node[above,yshift=-1pt]{$i$}--++(0,-1)node[below]{$\phantom{i}$};
\end{tikzpicture}
\quad
\text{if $i\rightsquigarrow j$}
,\\
\begin{tikzpicture}[anchorbase,smallnodes,rounded corners]
\draw[solid](1,1)node[above,yshift=-1pt]{$\phantom{i}$}--++(1,-1)node[below]{$i$};
\draw[solid](2,1)--++(-1,-1)node[below]{$i$};
\draw[redstring](1.5,1)--++(-0.5,-0.5)--++(0.5,-0.5)node[below]{$i$};
\end{tikzpicture}
=
\begin{tikzpicture}[anchorbase,smallnodes,rounded corners]
\draw[solid](3,1)node[above,yshift=-1pt]{$\phantom{i}$}--++(1,-1)node[below]{$i$};
\draw[solid](4,1)--++(-1,-1)node[below]{$i$};
\draw[redstring](3.5,1)--++(0.5,-0.5)--++(-0.5,-0.5)node[below]{$i$};
\end{tikzpicture}
-
\begin{tikzpicture}[anchorbase,smallnodes,rounded corners]
\draw[solid](5,1)node[above,yshift=-1pt]{$\phantom{i}$}--++(0,-1)node[below]{$i$};
\draw[solid](6,1)--++(0,-1)node[below]{$i$};
\draw[redstring](5.5,1)--++(0,-1)node[below]{$i$};
\end{tikzpicture}
.
\end{gathered}
\end{gather}
\end{enumerate}
\end{Definition}

The Reidemeister II or III relations without error terms are called \emph{honest} Reidemeister relations.

The multilocal relations are sometimes tricky to apply but they are crucial
for the combinatorics to work:

\begin{Example}\label{E:ProofsMultilocal}
Consider the following two diagrams
(in these pictures both solid $i$-strings have one ghost).
\begin{gather*}
\begin{tikzpicture}[anchorbase,smallnodes,rounded corners,xscale=1.5]
\draw[ghost](1,0)--++(0,1)node[above,yshift=-1pt]{$i$};
\draw[ghost](1.2,0)--++(0,1)node[above,yshift=-1pt]{$i$};
\draw[solid](0,0)node[below]{$i$}--++(0,1);
\draw[solid](0.2,0)node[below]{$i$}--++(0,1);
\end{tikzpicture}
\stackrel{\autoref{Eq:RecollectionDotCrossing}}{=}
\begin{tikzpicture}[anchorbase,smallnodes,rounded corners,xscale=1.5]
\draw[ghost](1,0)--++(0.2,1)node[above,yshift=-1pt]{$i$};
\draw[ghost,dot=0.8](1.2,0)--++(-0.2,1)node[above,yshift=-1pt]{$i$};
\draw[solid](0,0)node[below]{$i$}--++(0.2,1);
\draw[solid,dot=0.8](0.2,0)node[below]{$i$}--++(-0.2,1);
\end{tikzpicture}
-
\begin{tikzpicture}[anchorbase,smallnodes,rounded corners,xscale=1.5]
\draw[ghost](1,0)--++(0.2,1)node[above,yshift=-1pt]{$i$};
\draw[ghost,dot=0.2](1.2,0)--++(-0.2,1)node[above,yshift=-1pt]{$i$};
\draw[solid](0,0)node[below]{$i$}--++(0.2,1);
\draw[solid,dot=0.2](0.2,0)node[below]{$i$}--++(-0.2,1);
\end{tikzpicture}
,\qquad
\begin{tikzpicture}[anchorbase,smallnodes,rounded corners,xscale=1.5]
\draw[ghost](1,0)--++(0,1)node[above,yshift=-1pt]{$i$};
\draw[ghost](1.2,0)--++(0,1)node[above,yshift=-1pt]{$i$};
\draw[solid](0,0)node[below]{$i$}--++(0,1);
\draw[solid](0.2,0)node[below]{$i$}--++(0,1);
\draw[solid](1.1,0)node[below]{$k$}--++(0,1);
\end{tikzpicture}
.
\end{gather*}
One cannot decide locally around just the two solid
$i$-strings whether \autoref{Eq:RecollectionDotCrossing} can be applied. In the left-hand diagram, the relation can be applied, giving the linear combination of diagrams above.
In contrast, in the right-hand diagram relation \autoref{Eq:RecollectionDotCrossing} cannot be applied
because the two solid $i$-strings cannot be pulled arbitrarily close together using isotopies as the solid $k$-string prevents this from happening.
\end{Example}

\begin{Remark}\label{R:RecollectionDifferentStrings}
Assume that a vertex $i\in I$ has multiple outgoing edges 
with positive weighting, so that there
can be many ghost $i$-strings. For example:
\begin{gather*}
\begin{tikzpicture}[anchorbase]
\node[circle,inner sep=1.8pt,fill=DarkBlue] (0) at (0,0){};
\node at (0,-0.25){$0$};
\node[circle,inner sep=1.8pt,fill=DarkBlue] (1) at (1,0){};
\node at (1,-0.25){$1$};
\node[circle,inner sep=1.8pt,fill=DarkBlue] (2) at (2,0){};
\node at (2,-0.25){$2$};
\draw[directed=0.5](1)--(0)node[above,xshift=0.5cm]{$0.75$};
\draw[directed=0.5](1)--(2)node[above,xshift=-0.5cm]{$1.25$};
\end{tikzpicture}
\quad\rightsquigarrow\quad
\begin{tikzpicture}[anchorbase,smallnodes,rounded corners]
\draw[ghost](0.75,0)node[below]{$\phantom{i}$}--++(0,1)node[above,yshift=-0.05cm]{$1$};
\draw[ghost](1.25,0)node[below]{$\phantom{i}$}--++(0,1)node[above,yshift=-0.05cm]{$1$};
\draw[solid](0,0)node[below]{$1$}--++(0,1)node[above,yshift=-1pt]{$\phantom{i}$};
\end{tikzpicture}
.
\end{gather*}
Note that these two ghosts strings are different strings that behave differently because
the relations depend on the edges and on the residues.
In the example above, the two ghost $1$-strings
play different roles. For example, one of them has nontrivial
Reidemeister II relations with the $0$-strings and the other has nontrivial Reidemeister II relations with the $2$-strings.
\end{Remark}

\begin{Remark}\label{R:RecollectionNNotBeta}
Note that \cite{MaTu-klrw-algebras} mostly works with
$\WA(X)$, where $\beta$ determines the
labels of the solid strings. This difference is not very important since $\WA[n](X)=
\bigoplus_{\beta\in Q^{+}_{n}}\WA(X)$, where $Q^{+}_{n}$ is the set of positive roots of height $n$, which corresponds
to the set of possible labels on the $n$ solid strings.
\end{Remark}

A diagram is \emph{unsteady}
if it contains a solid string that can be pulled arbitrarily far to the right when the the red strings belong to region bounded by $X$, otherwise it is \emph{steady}. For example, we have
\begin{gather}\label{Eq:RecollectionUnsteady}
\text{Unsteady}\colon
\begin{tikzpicture}[anchorbase,smallnodes]
\draw[ghost](1,0)--++(0,1)node[above,yshift=-1pt]{$i$};
\draw[solid](0.5,0)node[below]{$i$}--++(0,1)node[above,yshift=-1pt]{$\phantom{i}$};
\draw[redstring](0.25,0)node[below]{$\rho$}--++(0,1)node[above,yshift=-1pt]{$\phantom{i}$};
\draw[->,decorate,decoration={snake,amplitude=0.4mm,segment length=2mm,post length=1mm},spinach] (0.5,0.5) to (1.5,0.5)node[right]{pulls freely};
\end{tikzpicture}
,\quad
\text{steady}\colon
\begin{tikzpicture}[anchorbase,smallnodes]
\draw[ghost](0.5,0)--++(0,1)node[above,yshift=-1pt]{$i$};
\draw[solid](0,0)node[below]{$i$}--++(0,1)node[above,yshift=-1pt]{$\phantom{i}$};
\draw[redstring](0.25,0)node[below]{$\rho$}--++(0,1)node[above,yshift=-1pt]{$\phantom{i}$};
\draw[->,decorate,decoration={snake,amplitude=0.4mm,segment length=2mm,post length=1mm},spinach] (0,0.5) to (1,0.5)node[right]{stopped by the red string};
\end{tikzpicture}
.
\end{gather}
Note that the ghost string mimics its parent solid string, so the solid $i$-string
in the left diagram does indeed pull freely to the right.

Finally, the \emph{cyclotomic weighted KLRW algebra} $\WAc[n](X)$
is the finite dimensional
quotient of $\WA[n](X)$ by the two-sided ideal generated by all
diagrams that factor through an unsteady diagram.


\subsection{Duality and partners}\label{SS:RecollectionDual}


We use the usual
\emph{diagrammatic antiinvolution} $(\placeholder)^{\star}$ given
by (the $R$-linear extension of) reflecting diagrams in their horizontal axis.
This antiinvolution is the one we use for the
homogeneous (affine) sandwich cellular basis.
An illustration of the diagrammatic antiinvolution is:
\begin{gather*}
\left(
\begin{tikzpicture}[anchorbase,smallnodes,rounded corners]
\draw[ghost](1.5,1)node[above,yshift=-1pt]{$j$}--++(1,-1);
\draw[ghost](2.5,1)node[above,yshift=-1pt]{$i$}--++(-1,-1);
\draw[ghost](4,1)node[above,yshift=-1pt]{$k$}--++(0.5,-1);
\draw[solid](1,1)node[above,yshift=-1pt]{$\phantom{i}$}--++(1,-1)node[below]{$j$};
\draw[solid](2,1)--++(-1,-1)node[below]{$i$};
\draw[solid](3,1)--++(0.5,-1)node[below]{$k$};
\draw[redstring](1.25,1)--++(0,-1)node[below]{$\rho$};
\end{tikzpicture}
\right)^{\star}
=
\begin{tikzpicture}[anchorbase,smallnodes,rounded corners]
\draw[ghost](1.5,1)node[above,yshift=-1pt]{$i$}--++(1,-1);
\draw[ghost](2.5,1)node[above,yshift=-1pt]{$j$}--++(-1,-1);
\draw[ghost](4,0)--++(0.5,1)node[above,yshift=-1pt]{$k$};
\draw[solid](1,1)node[above,yshift=-1pt]{$\phantom{i}$}--++(1,-1)node[below]{$i$};
\draw[solid](2,1)--++(-1,-1)node[below]{$j$};
\draw[solid](3,0)node[below]{$k$}--++(0.5,1);
\draw[redstring](1.25,1)--++(0,-1)node[below]{$\rho$};
\end{tikzpicture}
.
\end{gather*}

There is a different kind of duality on diagrams, which exists
because the relations
\autoref{Eq:RecollectionDotCrossing}--\autoref{Eq:RecollectionReidemeisterIII} come in mirrored pairs. Explicitly, up to scalars, the set of relations is invariant under reflecting diagrams in their vertical axis and, in some circumstances, swapping solid strings with their ghosts.
Most of the relations we will use have this type of duality, which we call \emph{partner relations}. Note that the partner relations can have different scalars, however, this will not play a role for us because we only use partner relations when similar arguments can be applied to both relations.

For example, if $i\rightsquigarrow j$, then
\begin{gather*}
\begin{tikzpicture}[anchorbase,smallnodes,rounded corners]
\draw[ghost](1,1)node[above,yshift=-1pt]{$i$}--++(1,-1)node[below]{$\phantom{i}$};
\draw[ghost](2,1)node[above,yshift=-1pt]{$i$}--++(-1,-1)node[below]{$\phantom{i}$};
\draw[solid,smallnodes,rounded corners](1.5,1)--++(-0.5,-0.5)--++(0.5,-0.5)node[below]{$j$};
\end{tikzpicture}
{=}
\begin{tikzpicture}[anchorbase,smallnodes,rounded corners]
\draw[ghost](3,1)node[above,yshift=-1pt]{$\phantom{i}$}--++(1,-1)node[below]{$\phantom{i}$};
\draw[ghost](4,1)node[above,yshift=-1pt]{$i$}--++(-1,-1)node[below]{$\phantom{i}$};
\draw[solid,smallnodes,rounded corners](3.5,1)--++(0.5,-0.5)--++(-0.5,-0.5)node[below]{$j$};
\end{tikzpicture}
-Q_{iji}(\bu)
\begin{tikzpicture}[anchorbase,smallnodes,rounded corners]
\draw[ghost](6.2,1)node[above,yshift=-1pt]{$i$}--++(0,-1)node[below]{$\phantom{i}$};
\draw[ghost](7.2,1)node[above,yshift=-1pt]{$i$}--++(0,-1)node[below]{$\phantom{i}$};
\draw[solid](6.7,1)--++(0,-1)node[below]{$j$};
\end{tikzpicture}
\xleftrightarrow[\text{relations}]{\text{partner}}
\begin{tikzpicture}[anchorbase,smallnodes,rounded corners]
\draw[solid](1,1)node[above,yshift=-1pt]{$\phantom{i}$}--++(1,-1)node[below]{$j$};
\draw[solid](2,1)--++(-1,-1)node[below]{$j$};
\draw[ghost,smallnodes,rounded corners](1.5,1)node[above,yshift=-1pt]{$i$}--++(-0.5,-0.5)--++(0.5,-0.5)node[below]{$\phantom{i}$};
\end{tikzpicture}
{=}
\begin{tikzpicture}[anchorbase,smallnodes,rounded corners]
\draw[solid](3,1)node[above,yshift=-1pt]{$\phantom{i}$}--++(1,-1)node[below]{$j$};
\draw[solid](4,1)--++(-1,-1)node[below]{$j$};
\draw[ghost,smallnodes,rounded corners](3.5,1)node[above,yshift=-1pt]{$i$}--++(0.5,-0.5)--++(-0.5,-0.5)node[below]{$\phantom{i}$};
\end{tikzpicture}
+Q_{iji}(\bu)
\begin{tikzpicture}[anchorbase,smallnodes,rounded corners]
\draw[solid](7.2,1)node[above,yshift=-1pt]{$\phantom{i}$}--++(0,-1)node[below]{$j$};
\draw[solid](8.2,1)--++(0,-1)node[below]{$j$};
\draw[ghost](7.7,1)node[above,yshift=-1pt]{$i$}--++(0,-1)node[below]{$\phantom{i}$};
\end{tikzpicture}
,
\end{gather*}
is an example of partner relations. Note the change of sign when going from left to right.


\subsection{Some diagrams that we need}\label{SS:RecollectionDiagrams}


We need certain special diagrams:
\begin{enumerate}

\item \emph{Idempotent diagrams} are diagrams
with no dots and no crossings, and fixed $x$-coordinates for all strings.

\item \emph{Straight line diagrams}  are diagrams
with no dots and no crossings.

\item \emph{Dotted idempotent diagrams}
and \emph{dotted straight line diagrams} are dotted diagrams
of the respective type.

\item \emph{Permutation diagrams} are diagrams with no dots such the solid strings correspond to a reduced expression of a permutation; see \cite[Definition 3B.1]{MaTu-klrw-algebras} for a more precise definition.

\end{enumerate}

\begin{Example}\label{E:RecollectionSomeDiagrams}
Examples of these types of diagrams are:
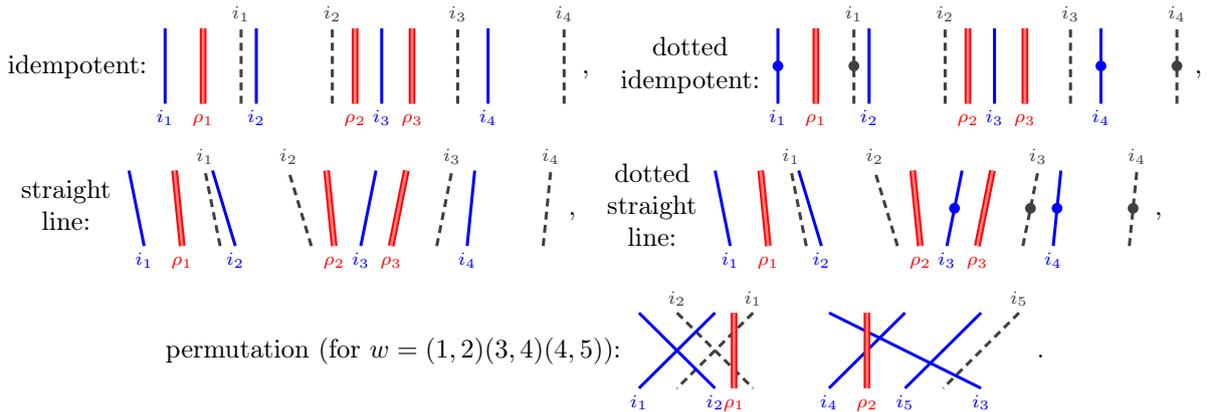
\begin{gather*}
\text{idempotent:}
\begin{tikzpicture}[anchorbase,smallnodes,rounded corners]
\draw[ghost](0,0)node[below]{$\phantom{i}$}--++(0,1)node[above,yshift=-1pt]{$i_{1}$};
\draw[ghost](1.2,0)node[below]{$\phantom{i}$}--++(0,1)node[above,yshift=-1pt]{$i_{2}$};
\draw[ghost](2.85,0)node[below]{$\phantom{i}$}--++(0,1)node[above,yshift=-1pt]{$i_{3}$};
\draw[ghost](4.25,0)node[below]{$\phantom{i}$}--++(0,1)node[above,yshift=-1pt]{$i_{4}$};
\draw[solid](-1,0)node[below]{$i_{1}$}--++(0,1) node[above,yshift=-1pt]{$\phantom{i_{1}}$};
\draw[solid](0.2,0)node[below]{$i_{2}$}--++(0,1);
\draw[solid](1.85,0)node[below]{$i_{3}$}--++(0,1);
\draw[solid](3.25,0)node[below]{$i_{4}$}--++(0,1);
\draw[redstring](-0.5,0)node[below]{$\rho_{1}$}--++(0,1);
\draw[redstring](1.5,0)node[below]{$\rho_{2}$}--++(0,1);
\draw[redstring](2.25,0)node[below]{$\rho_{3}$}--++(0,1);
\end{tikzpicture}
,\quad
\begin{gathered}
\text{dotted}
\\[-0.1cm]
\text{idempotent:}
\end{gathered}
\begin{tikzpicture}[anchorbase,smallnodes,rounded corners]
\draw[ghost,dot](0,0)node[below]{$\phantom{i}$}--++(0,1)node[above,yshift=-1pt]{$i_{1}$};
\draw[ghost](1.2,0)node[below]{$\phantom{i}$}--++(0,1)node[above,yshift=-1pt]{$i_{2}$};
\draw[ghost](2.85,0)node[below]{$\phantom{i}$}--++(0,1)node[above,yshift=-1pt]{$i_{3}$};
\draw[ghost,dot](4.25,0)node[below]{$\phantom{i}$}--++(0,1)node[above,yshift=-1pt]{$i_{4}$};
\draw[solid,dot](-1,0)node[below]{$i_{1}$}--++(0,1) node[above,yshift=-1pt]{$\phantom{i_{1}}$};
\draw[solid](0.2,0)node[below]{$i_{2}$}--++(0,1);
\draw[solid](1.85,0)node[below]{$i_{3}$}--++(0,1);
\draw[solid,dot](3.25,0)node[below]{$i_{4}$}--++(0,1);
\draw[redstring](-0.5,0)node[below]{$\rho_{1}$}--++(0,1);
\draw[redstring](1.5,0)node[below]{$\rho_{2}$}--++(0,1);
\draw[redstring](2.25,0)node[below]{$\rho_{3}$}--++(0,1);
\end{tikzpicture}
,
\\
\begin{gathered}
\text{straight}
\\[-0.1cm]
\text{line:}
\end{gathered}
\begin{tikzpicture}[anchorbase,smallnodes,rounded corners]
\draw[ghost](0,0)node[below]{$\phantom{i}$}--++(-0.2,1)node[above,yshift=-1pt]{$i_{1}$};
\draw[ghost](1.2,0)node[below]{$\phantom{i}$}--++(-0.3,1)node[above,yshift=-1pt]{$i_{2}$};
\draw[ghost](2.85,0)node[below]{$\phantom{i}$}--++(0.2,1)node[above,yshift=-1pt]{$i_{3}$};
\draw[ghost](4.25,0)node[below]{$\phantom{i}$}--++(0.1,1)node[above,yshift=-1pt]{$i_{4}$};
\draw[solid](-1,0)node[below]{$i_{1}$}--++(-0.2,1) node[above,yshift=-1pt]{$\phantom{i_{1}}$};
\draw[solid](0.2,0)node[below]{$i_{2}$}--++(-0.3,1);
\draw[solid](1.85,0)node[below]{$i_{3}$}--++(0.2,1);
\draw[solid](3.25,0)node[below]{$i_{4}$}--++(0.1,1);
\draw[redstring](-0.5,0)node[below]{$\rho_{1}$}--++(-0.1,1);
\draw[redstring](1.5,0)node[below]{$\rho_{2}$}--++(-0.1,1);
\draw[redstring](2.25,0)node[below]{$\rho_{3}$}--++(0.2,1);
\end{tikzpicture}
,\quad
\begin{gathered}
\text{dotted}
\\[-0.1cm]
\text{straight}
\\[-0.1cm]
\text{line:}
\end{gathered}
\begin{tikzpicture}[anchorbase,smallnodes,rounded corners]
\draw[ghost](0,0)node[below]{$\phantom{i}$}--++(-0.2,1)node[above,yshift=-1pt]{$i_{1}$};
\draw[ghost](1.2,0)node[below]{$\phantom{i}$}--++(-0.3,1)node[above,yshift=-1pt]{$i_{2}$};
\draw[ghost,dot](2.85,0)node[below]{$\phantom{i}$}--++(0.2,1)node[above,yshift=-1pt]{$i_{3}$};
\draw[ghost,dot](4.25,0)node[below]{$\phantom{i}$}--++(0.1,1)node[above,yshift=-1pt]{$i_{4}$};
\draw[solid](-1,0)node[below]{$i_{1}$}--++(-0.2,1) node[above,yshift=-1pt]{$\phantom{i_{1}}$};
\draw[solid](0.2,0)node[below]{$i_{2}$}--++(-0.3,1);
\draw[solid,dot](1.85,0)node[below]{$i_{3}$}--++(0.2,1);
\draw[solid,dot](3.25,0)node[below]{$i_{4}$}--++(0.1,1);
\draw[redstring](-0.5,0)node[below]{$\rho_{1}$}--++(-0.1,1);
\draw[redstring](1.5,0)node[below]{$\rho_{2}$}--++(-0.1,1);
\draw[redstring](2.25,0)node[below]{$\rho_{3}$}--++(0.2,1);
\end{tikzpicture}
,\quad
\\
\text{permutation (for $w=(1,2)(3,4)(4,5)$):}
\begin{tikzpicture}[anchorbase,smallnodes,rounded corners]
\draw[ghost](1.5,1)node[above,yshift=-1pt]{$i_{2}$}--++(1,-1);
\draw[ghost](2.5,1)node[above,yshift=-1pt]{$i_{1}$}--++(-1,-1);
\draw[ghost](6,1)node[above,yshift=-1pt]{$i_{5}$}--++(-1,-1)node[below]{$\phantom{i}$};
\draw[solid](1,1)node[above,yshift=-1pt]{$\phantom{i}$}--++(1,-1)node[below]{$i_{2}$};
\draw[solid](2,1)--++(-1,-1)node[below]{$i_{1}$};
\draw[solid](3.5,1)node[above,yshift=-1pt]{$\phantom{i}$}--++(2,-1)node[below]{$i_{3}$};
\draw[solid](4.5,1)node[above,yshift=-1pt]{$\phantom{i}$}--++(-1,-1)node[below]{$i_{4}$};
\draw[solid](5.5,1)node[above,yshift=-1pt]{$\phantom{i}$}--++(-1,-1)node[below]{$i_{5}$};
\draw[redstring](2.25,1)--++(0,-1)node[below]{$\rho_{1}$};
\draw[redstring](4,1)--++(0,-1)node[below]{$\rho_{2}$};
\end{tikzpicture}
.
\end{gather*}
These illustrate (dotted) idempotents,
(dotted) straight line and permutation diagrams.
\end{Example}

\begin{Notation}\label{N:RecollectionNotationNeeded2}
Let $\idem[{\bx,\bi}]$ be the idempotent diagram with bottom boundary given by $(\bx,\bi)$,
for $x\in X$ and $\bi\in I^{n}$. (That is, the positions and labels of the solid strings are given by $\bx$ and $\bi$, respectively, when read from left to right.) As in \autoref{SS:RecollectionDefinition}(i),
given a polynomial $p(\by)\in R[y_{1},\dots,y_{n}]$] write
$p(\by)\idem[{\bx,\bi}]$ for the corresponding linear combination of (dotted) diagrams $\idem[{\bx,\bi}]$.

Additionally, let $\Sym=\langle s_{1},\dots,s_{n-1}\rangle$ be the symmetric group on $\set{1,\dots,n}$ with $s_{i}=(i,i+1)$, which we think of as swapping the $i$th and $(i+1)$th solid strings.
For all $w\in\Sym$ fix a reduced expression and let $D(w)$ be the
associated permutation diagram. As often happens in the KLR world, the diagram $D(w)$ is only well-defined up to a choice of reduced expression and, as explained in \cite[Definition 3B.1]{MaTu-klrw-algebras}, some care needs to be taken. (These details are not important in this paper.)
\end{Notation}

\begin{Example}\label{E:RecollectionNotationNeeded2}
In pictures:
\begin{gather*}
\idem[{(0,0.5),(i,j)}]=
\begin{tikzpicture}[anchorbase,smallnodes,rounded corners]
\draw[ghost](1,0)--++(0,1)node[above,yshift=-1pt]{$i$};
\draw[ghost](1.5,0)--++(0,1)node[above,yshift=-1pt]{$j$};
\draw[solid](0,0)node[below]{$i$}--++(0,1);
\draw[solid](0.5,0)node[below]{$j$}--++(0,1);
\end{tikzpicture}
,\quad
y_{2}\idem[{(0,0.5),(i,j)}]=
\begin{tikzpicture}[anchorbase,smallnodes,rounded corners]
\draw[ghost](1,0)--++(0,1)node[above,yshift=-1pt]{$i$};
\draw[ghost,dot](1.5,0)--++(0,1)node[above,yshift=-1pt]{$j$};
\draw[solid](0,0)node[below]{$i$}--++(0,1);
\draw[solid,dot](0.5,0)node[below]{$j$}--++(0,1);
\end{tikzpicture}
,\quad
D(s_{1})\idem[{(0,0.5),(i,j)}]
=
\begin{tikzpicture}[anchorbase,smallnodes,rounded corners]
\draw[ghost](1,0)--++(0.51,1)node[above,yshift=-1pt]{$i$};
\draw[ghost](1.5,0)--++(-0.5,1)node[above,yshift=-1pt]{$j$};
\draw[solid](0,0)node[below]{$i$}--++(0.5,1);
\draw[solid](0.5,0)node[below]{$j$}--++(-0.51,1);
\end{tikzpicture}
.
\end{gather*}
Here we have not drawn red strings while the two solid
strings have $x$-coordinates $0$ and $0.5$.
\end{Example}

\begin{Definition}\label{D:RecollectionFactors}
A diagram $D$ \emph{factors through} $S$ if $D=D^{\prime}SD^{\prime\prime}$,
for some $D^{\prime},D^{\prime\prime}\in\WA[n](X)$.
\end{Definition}

Given a straight line diagram $S$ and some $0<\varepsilon\ll 1$, then there
exists an idempotent
diagram $L(S)$, called the \emph{left justification} of $S$,
such that $S$ factors through $L(S)$, and $L(S)$ has its strings
as far to the left as possible while its coordinates are
within the interval defined by $X$ and strings are at least $\varepsilon$ apart.
(The $\varepsilon>0$ is only needed to make $L(S)$ well-defined,
and we omit it below.)
See \cite[Section 6D]{MaTu-klrw-algebras} for a detailed account.

\begin{Example}\label{E:RecollectionLeftJustification}
Here is an example of a left justification (and also an example of factorization of diagrams):
\begin{gather*}
S=
\begin{tikzpicture}[anchorbase,smallnodes,rounded corners]
\draw[ghost](0,0)node[below]{$\phantom{i}$}--++(-0.2,1)node[above,yshift=-1pt]{$i_{1}$};
\draw[ghost](1.2,0)node[below]{$\phantom{i}$}--++(-0.3,1)node[above,yshift=-1pt]{$i_{2}$};
\draw[ghost](2.85,0)node[below]{$\phantom{i}$}--++(0.2,1)node[above,yshift=-1pt]{$i_{3}$};
\draw[ghost](4.25,0)node[below]{$\phantom{i}$}--++(0.1,1)node[above,yshift=-1pt]{$i_{4}$};
\draw[solid](-1,0)node[below]{$i_{1}$}--++(-0.2,1) node[above,yshift=-1pt]{$\phantom{i_{1}}$};
\draw[solid](0.2,0)node[below]{$i_{2}$}--++(-0.3,1);
\draw[solid](1.85,0)node[below]{$i_{3}$}--++(0.2,1);
\draw[solid](3.25,0)node[below]{$i_{4}$}--++(0.1,1);
\draw[redstring](-0.5,0)node[below]{$\rho_{1}$}--++(-0.1,1);
\draw[redstring](1.5,0)node[below]{$\rho_{2}$}--++(-0.1,1);
\draw[redstring](2.25,0)node[below]{$\rho_{3}$}--++(0.2,1);
\draw[very thick,mygray,densely dotted] (-1.25,0.5) to (4.5,0.5);
\end{tikzpicture}
,\quad
L(S)=
\begin{tikzpicture}[anchorbase,smallnodes,rounded corners]
\draw[ghost](0,0)node[below]{$\phantom{i}$}--++(0,1)node[above,yshift=-1pt]{$i_{1}$};
\draw[ghost](1.1,0)node[below]{$\phantom{i}$}--++(0,1)node[above,yshift=-1pt]{$i_{2}$};
\draw[ghost](2.3,0)node[below]{$\phantom{i}$}--++(0,1)node[above,yshift=-1pt]{$i_{3}$};
\draw[ghost](3.4,0)node[below]{$\phantom{i}$}--++(0,1)node[above,yshift=-1pt]{$i_{4}$};
\draw[solid](-1.0,0)node[below]{$i_{1}$}--++(0,1) node[above,yshift=-1pt]{$\phantom{i_{1}}$};
\draw[solid](0.1,0)node[below]{$i_{2}$}--++(0,1) node[above,yshift=-1pt]{$\phantom{i_{1}}$};
\draw[solid](1.3,0)node[below,xshift=0.075cm]{$i_{3}$}--++(0,1) node[above,yshift=-1pt]{$\phantom{i_{1}}$};
\draw[solid](2.4,0)node[below]{$i_{4}$}--++(0,1) node[above,yshift=-1pt]{$\phantom{i_{1}}$};
\draw[redstring](-0.9,0)--++(0,1)node[above]{$\rho_{1}$};
\draw[redstring](1.2,0)node[below,xshift=-0.075cm]{$\rho_{2}$}--++(0,1);
\draw[redstring](1.4,0)--++(0,1)node[above]{$\rho_{3}$};
\end{tikzpicture}
.
\end{gather*}
In general, the left justification is obtained by
considering a collar neighborhood of a horizontal cut of $S$.
To construct $L(S)$, we use isotopies to inductively pull 
strings to the left inside this neighborhood.
\end{Example}

In a dotted straight line diagram two strings are \emph{close} if you can pull them arbitrarily close together in a neighborhood that does not contain any other strings using only (multilocal) isotopies. We use \emph{close and to the left and right} in the obvious way.

\begin{Example}\label{E:ProofsClose}
Let us repeat \cite[Example 6D.10]{MaTu-klrw-algebras}:
\begin{gather*}
\begin{tikzpicture}[anchorbase,smallnodes,rounded corners]
\draw[ghost](1,0)--++(0,1)node[above,yshift=-1pt]{$i$};
\draw[ghost](1.2,0)--++(0,1)node[above,yshift=-1pt]{$j$};
\draw[solid](0,0)node[below]{$i$}--++(0,1);
\draw[solid](0.2,0)node[below]{$j$}--++(0,1);
\end{tikzpicture}
,\quad
\begin{tikzpicture}[anchorbase,smallnodes,rounded corners]
\draw[ghost](1,0)--++(0,1)node[above,yshift=-1pt]{$i$};
\draw[ghost](1.2,0)--++(0,1)node[above,yshift=-1pt]{$j$};
\draw[solid](0,0)node[below]{$i$}--++(0,1);
\draw[solid](0.2,0)node[below]{$j$}--++(0,1);
\draw[solid](1.1,0)node[below]{$k$}--++(0,1);
\end{tikzpicture}
.
\end{gather*}
The solid $i$-string
is close to the solid $j$-string in the left but not in the right-hand diagram.
\end{Example}

One of the other crucial points about weighted KLRW algebras
is that they generalize KLR algebras, see \cite[Section 3F]{MaTu-klrw-algebras}.
The construction given therein uses a particular choice $X$ of \emph{KLRW positioning} such that for a special $\bx\in X$ the KLR algebra is graded
isomorphic to $\idem[{\bx}]\WA(X)\idem[{\bx}]$, where $\idem[{\bx}]=\sum_{\bi\in I}\idem[{\bx,\bi}]$. We return to this point in \autoref{SS:BasesConstruction} below.


\subsection{A basis and a faithful module}\label{SS:RecollectionFaithful}


The following \emph{standard basis} proposition works only for the infinite dimensional
weighted KLRW algebras. One of the main features of our approach is that we
can use the faithful representation to prove cellularity of the finite dimensional quotients as well.

\begin{Proposition}\label{P:RecollectionWABasis}
The algebra $\WA[n](X)$ is free as an $R$-module with homogeneous basis
\begin{gather*}
\WABasis=
\set{D(w)y_{1}^{a_{1}}\dots y_{n}^{a_{n}}\idem[{\bx,\bi}]|a_{1},\dots,a_{n}\in\N,w\in\Sym,\bx\in X,\bi\in I^{n}}.
\end{gather*}
\end{Proposition}

\begin{proof}
See \cite[Proposition 3B.12]{MaTu-klrw-algebras} for details.
\end{proof}

Borrowing the pictures from \autoref{S:BasesSandwich}, the standard
basis is of the form
\begin{gather*}
D(w)\underbrace{y_{1}^{a_{1}}\dots y_{n}^{a_{n}}\idem[{\bx,\bi}]}_{=d}
\leftrightsquigarrow
\begin{tikzpicture}[anchorbase,scale=1]
\draw[line width=0.75,color=black,fill=cream] (0,1) to (0.25,0.5) to (0.75,0.5) to (1,1) to (0,1);
\node at (0.5,0.75){\scalebox{0.75}{$D(w)$}};
\draw[line width=0.75,color=black,fill=cream] (0.25,0) to (0.25,0.5) to (0.75,0.5) to (0.75,0) to (0.25,0);
\node at (0.5,0.25){$d$};
\end{tikzpicture}
,
\quad
\begin{aligned}
\begin{tikzpicture}[anchorbase,scale=1]
\draw[line width=0.75,color=black,fill=cream] (0,1) to (0.25,0.5) to (0.75,0.5) to (1,1) to (0,1);
\node at (0.5,0.75){\scalebox{0.75}{$D(w)$}};
\end{tikzpicture}
&\text{ a permutation diagram,}
\\
\begin{tikzpicture}[anchorbase,scale=1]
\draw[white,ultra thin] (0,0) to (1,0);
\draw[line width=0.75,color=black,fill=cream] (0.25,0) to (0.25,0.5) to (0.75,0.5) to (0.75,0) to (0.25,0);
\node at (0.5,0.25){$d$};
\end{tikzpicture}
&\text{ dots on an idempotent}.
\end{aligned}
\end{gather*}
´
\begin{Remark}\label{R:RecollectionWABasis}
The proof of \autoref{P:RecollectionWABasis} in \cite{MaTu-klrw-algebras} uses a faithful action on
the polynomial ring
\begin{gather*}
P_{n}(X)=
\bigoplus_{\bx\in X,\bi\in I^{n}}
R[y_{1},\dots,y_{n}]\idem[{\bx,\bi}].
\end{gather*}
We will not recall the action here as it is quite standard in the field,
see \cite[Section 3C]{MaTu-klrw-algebras} for details.
\end{Remark}


\section{The strategy}\label{S:Strategy}
\stepcounter{subsection}

Our strategy to construct homogeneous (affine) sandwich cellular
bases for the weighted KLRW algebras is the same as in
\cite[Remark 6.1]{MaTu-klrw-algebras}, which we now recall.
We will define 
the relevant notions later. The ideas are summarized in \autoref{R:StrategyMain} and
\autoref{E:StrategyMain} below.

\begin{Remark}\label{R:StrategyMain}
A crucial ingredient is the idea of using a certain form of \emph{maximality}, with respect to string placement.
\begin{enumerate}

\item Fix a residue sequence $\bi_{\blam}=(i_{1},\dots,i_{n})\in I^{n}$. We construct an idempotent diagram $\idem$
by placing strings inductively as far to the right as possible. Here we use the fact that the relations in \autoref{D:RecollectionwKLRW} only allow strings to be pulled to the right in certain situations, such as when they carry a dot. In this way, we think of previously placed strings as keeping the new string in check. One can think of this process as follows. We place the solid strings, and their ghosts, one at a time as follows. The positions of the red strings are determined by $\charge$. Next, place the solid string of residue $i_{1}$ to the left of all of the red strings and pull it to the right using isotopies and the honest Reidemeister II relations until it is blocked by an (affine) red string. Now place the string of residue $i_{2}$ to the left of all of the strings and pull it to the right until it is blocked by an existing string. Repeat this process to place the remaining strings, starting at the left of the existing strings and pulling it to the right until it becomes stuck. The idempotent diagram $\idem$ is then a horizontal cut in the middle of the resulting diagram.

This ensures that $\idem$ is maximal 
with respect to placing the strings to the right. In fact,
this strategy is a greedy algorithm,
as it is designed to be locally maximal but it produces
a globally maximal diagram.

\item By \autoref{D:RecollectionwKLRW} again, 
when we put dots on the strings in $\idem$ it might be possible to pull some strings further to the right, so that the dotted idempotent diagram is no longer maximal with respect to string placement, and even more problematic, might become unsteady. In the finite dimensional case we consider only those dotted idempotent diagrams $y^{\bfi}y_{\blam}\idem$ that are steady. In the infinite dimensional case we even can add extra dots and get $y^{\ba}y^{\bfi}y_{\blam}\idem$. This gives a collection of dotted idempotent diagrams, which
form the \emph{middle} of the homogeneous (affine) sandwich cellular
basis, in the sense of \autoref{R:SandwichCellularAlgebra}.

\item The (sandwich) cellular basis
is then obtained by a standard construction for diagram
algebras, which in our case means putting semistandard permutation diagrams above and below $y^{\ba}y^{\bfi}y_{\blam}\idem$.

\end{enumerate}
The basis itself is maximal (in a certain sense), by construction, and it is not
hard to prove that it is indeed a homogeneous (affine) sandwich cellular
basis. For example, putting additional dots on the basis
elements allows one pull
strings and jump dots to the right, making the resulting diagram bigger.
This gives an inductive way of proving results.

This strategy works perfectly in types $\typea$ and $\aone$, see 
\cite{MaTu-klrw-algebras}, and requires small adjustments in the $\bad$ types.
\end{Remark}

\begin{Example}\label{E:StrategyMain}
We now give an example of our main construction. This example is meant
for the reader to come back as they read through the text, so not all of the notation is defined at this point. This example is in level $1$ and there is one $0$-red string.

\begin{enumerate}

\item Consider type $\dtwo[2]$ and the following sequence of
Young diagrams, adding nodes in
row reading order, with the nodes filled with their residues.
\begin{align*}
\emptyset
\quad&\xrightarrow{\ 0\ }\quad
(1)=\begin{ytableau}
*(tomato!50) 0
\end{ytableau}
\\
&\xrightarrow{\ 1\ }\quad
(2)=\begin{ytableau}
*(tomato!50) 0 & 1
\end{ytableau}
\\
&\xrightarrow{\ 2\ }\quad
(3)=\begin{ytableau}
*(tomato!50) 0 & 1 & *(spinach!50) 2
\end{ytableau}
\\
&\xrightarrow{\ 2\ }\quad
(4)=\begin{ytableau}
*(tomato!50) 0 & 1 & *(spinach!50) 2 & *(spinach!50) 2 \\
\none
\end{ytableau}
\\
&\xrightarrow{\ 0\ }\quad
(4,1)=\begin{ytableau}
*(tomato!50) 0 & 1 & *(spinach!50) 2 & *(spinach!50) 2 \\
\none & *(tomato!50) 0
\end{ytableau}
\\
&\xrightarrow{\ 1\ }\quad
(4,2)=\begin{ytableau}
*(tomato!50) 0 & 1 & *(spinach!50) 2 & *(spinach!50) 2 \\
\none & *(tomato!50) 0 & 1
\end{ytableau}
\quad=\blam
.
\end{align*}
The final Young diagram $\blam$ has this associated Young diagram sequence.
We also associate to $\blam$ a sequence of idempotent
diagrams following the Young diagram sequence of $\blam$:
\begin{align*}
\idem[(0)]=\begin{tikzpicture}[anchorbase]
\draw[redstring](0,1)node[above,yshift=-1pt]{$\phantom{i}$}--++(0,-1)node[below]{$0$};
\end{tikzpicture}
\quad&\xrightarrow{\ 0\ }\quad
\idem[(1)]=\raisebox{-0.23cm}{$\DottedIdempotentD[]{2}{0}$}
\\
&\xrightarrow{\ 1\ }\quad
\idem[(2)]=\DottedIdempotentD[]{2}{0,1}
\\
&\xrightarrow{\ 2\ }\quad
\idem[(3)]=\DottedIdempotentD[]{2}{0,1,2}
\\
&\xrightarrow{\ 2\ }\quad
\idem[(4)]=\DottedIdempotentD[]{2}{0,1,2,2}
\\
&\xrightarrow{\ 0\ }\quad
\idem[(4,1)]=\DottedIdempotentD[]{2}{0,1,2,2,0}
\\
&\xrightarrow{\ 1\ }\quad
\idem[(4,2)]=\DottedIdempotentD[]{2}{0,1,2,2,0,1}
=\idem.
\end{align*}
Note that the solid $1$-string has two ghosts since we are in type $\dtwo[2]$.
In each step, the new string is as far to the right as
possible in the sense that we start at the very far left and
then pull it to the right using isotopies and honest
Reidemeister II relations whenever they apply. For example, the
first solid $0$-string gets stuck at the red $0$-string
because there is no honest Reidemeister II relation applies. The
resulting diagram is of degree $0$ and this diagram is the \emph{idempotent diagram} $\idem$
associated to $\blam$.

\item In the basis that we construct, we will allow two close solid $2$-strings to cross in $\idem$. In order for \autoref{Eq:RecollectionReidemeisterII}
not to annihilate the diagram we put a dot on the right-hand solid $2$-string,
so that
\begin{gather*}
\dotidem=\DottedIdempotentD[3]{2}{0,1,2,2,0,1}
\end{gather*}
is the \emph{dotted idempotent diagram}
$\dotidem$ associated to $\blam$. This diagram is of degree $2$.
Flanking with two crossings on the two close solid $2$-strings, gives us four basis elements, so far, in total. Note that we recover the diagram without the dot from the diagram with crossings both above and below the dotted idempotent diagram.

\item Note that the ghost of the left-hand solid $1$-string in $\idem$ is blocked by a solid $0$-string (and also solid $2$-strings). This gives a nontrivial sandwiched part for the finite dimensional case: we allow a dot on the left-hand solid $1$-string in $\idem$, so we now have eight diagrams of ``shape'' $\blam$ in total. The diagrams without crossing are
\begin{gather*}
\dotidem=
\begin{tikzpicture}[anchorbase,smallnodes,rounded corners]
\draw[dghost](-1,0)--++(0,1)node[above,yshift=-1pt]{$1$};
\draw[dghost](-0.2,0)--++(0,1)node[above,yshift=-1pt]{$1$};
\draw[solid](-4.1,0)node[below]{$1$}--++(0,1)node[above,yshift=-1pt]{$\phantom{i}$};
\draw[solid](-3.3,0)node[below]{$1$}--++(0,1)node[above,yshift=-1pt]{$\phantom{i}$};
\draw[solid](-0.8,0)node[below]{$0$}--++(0,1)node[above,yshift=-1pt]{$\phantom{i}$};
\draw[solid](-0.6,0)node[below]{$2$}--++(0,1)node[above,yshift=-1pt]{$\phantom{i}$};
\draw[solid,dot](-0.4,0)node[below]{$2$}--++(0,1)node[above,yshift=-1pt]{$\phantom{i}$};
\draw[solid](0,0)node[below]{$0$}--++(0,1)node[above,yshift=-1pt]{$\phantom{i}$};
\draw[redstring](0.25,0)node[below]{$0$}--++(0,1);
\end{tikzpicture}
,\quad
y_{1}\dotidem=
\begin{tikzpicture}[anchorbase,smallnodes,rounded corners]
\draw[dghost,dot](-1,0)--++(0,1)node[above,yshift=-1pt]{$1$};
\draw[dghost](-0.2,0)--++(0,1)node[above,yshift=-1pt]{$1$};
\draw[solid,dot](-4.1,0)node[below]{$1$}--++(0,1)node[above,yshift=-1pt]{$\phantom{i}$};
\draw[solid](-3.3,0)node[below]{$1$}--++(0,1)node[above,yshift=-1pt]{$\phantom{i}$};
\draw[solid](-0.8,0)node[below]{$0$}--++(0,1)node[above,yshift=-1pt]{$\phantom{i}$};
\draw[solid](-0.6,0)node[below]{$2$}--++(0,1)node[above,yshift=-1pt]{$\phantom{i}$};
\draw[solid,dot](-0.4,0)node[below]{$2$}--++(0,1)node[above,yshift=-1pt]{$\phantom{i}$};
\draw[solid](0,0)node[below]{$0$}--++(0,1)node[above,yshift=-1pt]{$\phantom{i}$};
\draw[redstring](0.25,0)node[below]{$0$}--++(0,1);
\end{tikzpicture}
.
\end{gather*}
Write $\psi_{4}$ for the crossing of the solid $2$-strings, which are the fourth and fifth strings read from the left. Then
the basis set is
\begin{gather*}
\set{\dotidem,y_{1}\dotidem,\psi_{4}\dotidem,\psi_{4}y_{1}\dotidem,\dotidem\psi_{4},y_{1}\dotidem\psi_{4},\psi_{4}\dotidem\psi_{4},\psi_{4}y_{1}\dotidem\psi_{4}}.
\end{gather*}
For example, the diagram for $\psi_{4}y_{1}\dotidem\psi_{4}$ is
\begin{gather}\label{Eq:StrategyMainBig}
\psi_{4}y_{1}\dotidem\psi_{4}
=
\begin{gathered}
\begin{tikzpicture}[anchorbase,smallnodes,rounded corners]
\draw[dghost](-1,0)--++(0,1)node[above,yshift=-1pt]{$1$};
\draw[dghost](-0.2,0)--++(0,1)node[above,yshift=-1pt]{$1$};
\draw[solid](-4.1,0)node[below]{$1$}--++(0,1)node[above,yshift=-1pt]{$\phantom{i}$};
\draw[solid](-3.3,0)node[below]{$1$}--++(0,1)node[above,yshift=-1pt]{$\phantom{i}$};
\draw[solid](-0.8,0)node[below]{$0$}--++(0,1)node[above,yshift=-1pt]{$\phantom{i}$};
\draw[solid](-0.6,0)node[below]{$2$}--++(0.2,1)node[above,yshift=-1pt]{$\phantom{i}$};
\draw[solid](-0.4,0)node[below]{$2$}--++(-0.2,1)node[above,yshift=-1pt]{$\phantom{i}$};
\draw[solid](0,0)node[below]{$0$}--++(0,1)node[above,yshift=-1pt]{$\phantom{i}$};
\draw[redstring](0.25,0)node[below]{$0$}--++(0,1);
\end{tikzpicture}
\\[-0.5cm]
\begin{tikzpicture}[anchorbase,smallnodes,rounded corners]
\draw[dghost,dot](-1,0)--++(0,1)node[above,yshift=-1pt]{$1$};
\draw[dghost](-0.2,0)--++(0,1)node[above,yshift=-1pt]{$1$};
\draw[solid,dot](-4.1,0)node[below]{$1$}--++(0,1)node[above,yshift=-1pt]{$\phantom{i}$};
\draw[solid](-3.3,0)node[below]{$1$}--++(0,1)node[above,yshift=-1pt]{$\phantom{i}$};
\draw[solid](-0.8,0)node[below]{$0$}--++(0,1)node[above,yshift=-1pt]{$\phantom{i}$};
\draw[solid](-0.6,0)node[below]{$2$}--++(0,1)node[above,yshift=-1pt]{$\phantom{i}$};
\draw[solid](-0.4,0)node[below]{$2$}--++(0,1)node[above,yshift=-1pt]{$\phantom{i}$};
\draw[solid](0,0)node[below]{$0$}--++(0,1)node[above,yshift=-1pt]{$\phantom{i}$};
\draw[redstring](0.25,0)node[below]{$0$}--++(0,1);
\end{tikzpicture}
\\[-0.5cm]
\begin{tikzpicture}[anchorbase,smallnodes,rounded corners]
\draw[dghost](-1,0)--++(0,1)node[above,yshift=-1pt]{$1$};
\draw[dghost](-0.2,0)--++(0,1)node[above,yshift=-1pt]{$1$};
\draw[solid](-4.1,0)node[below]{$1$}--++(0,1)node[above,yshift=-1pt]{$\phantom{i}$};
\draw[solid](-3.3,0)node[below]{$1$}--++(0,1)node[above,yshift=-1pt]{$\phantom{i}$};
\draw[solid](-0.8,0)node[below]{$0$}--++(0,1)node[above,yshift=-1pt]{$\phantom{i}$};
\draw[solid](-0.6,0)node[below]{$2$}--++(0,1)node[above,yshift=-1pt]{$\phantom{i}$};
\draw[solid,dot](-0.4,0)node[below]{$2$}--++(0,1)node[above,yshift=-1pt]{$\phantom{i}$};
\draw[solid](0,0)node[below]{$0$}--++(0,1)node[above,yshift=-1pt]{$\phantom{i}$};
\draw[redstring](0.25,0)node[below]{$0$}--++(0,1);
\end{tikzpicture}
\\[-0.5cm]
\begin{tikzpicture}[anchorbase,smallnodes,rounded corners]
\draw[dghost](-1,0)--++(0,1)node[above,yshift=-1pt]{$1$};
\draw[dghost](-0.2,0)--++(0,1)node[above,yshift=-1pt]{$1$};
\draw[solid](-4.1,0)node[below]{$1$}--++(0,1)node[above,yshift=-1pt]{$\phantom{i}$};
\draw[solid](-3.3,0)node[below]{$1$}--++(0,1)node[above,yshift=-1pt]{$\phantom{i}$};
\draw[solid](-0.8,0)node[below]{$0$}--++(0,1)node[above,yshift=-1pt]{$\phantom{i}$};
\draw[solid](-0.6,0)node[below]{$2$}--++(0,1)node[above,yshift=-1pt]{$\phantom{i}$};
\draw[solid](-0.4,0)node[below]{$2$}--++(0,1)node[above,yshift=-1pt]{$\phantom{i}$};
\draw[solid](0,0)node[below]{$0$}--++(0,1)node[above,yshift=-1pt]{$\phantom{i}$};
\draw[redstring](0.25,0)node[below]{$0$}--++(0,1);
\end{tikzpicture}
\\[-0.5cm]
\begin{tikzpicture}[anchorbase,smallnodes,rounded corners]
\draw[dghost](-1,0)--++(0,1)node[above,yshift=-1pt]{$1$};
\draw[dghost](-0.2,0)--++(0,1)node[above,yshift=-1pt]{$1$};
\draw[solid](-4.1,0)node[below]{$1$}--++(0,1)node[above,yshift=-1pt]{$\phantom{i}$};
\draw[solid](-3.3,0)node[below]{$1$}--++(0,1)node[above,yshift=-1pt]{$\phantom{i}$};
\draw[solid](-0.8,0)node[below]{$0$}--++(0,1)node[above,yshift=-1pt]{$\phantom{i}$};
\draw[solid](-0.6,0)node[below]{$2$}--++(0.2,1)node[above,yshift=-1pt]{$\phantom{i}$};
\draw[solid](-0.4,0)node[below]{$2$}--++(-0.2,1)node[above,yshift=-1pt]{$\phantom{i}$};
\draw[solid](0,0)node[below]{$0$}--++(0,1)node[above,yshift=-1pt]{$\phantom{i}$};
\draw[redstring](0.25,0)node[below]{$0$}--++(0,1);
\end{tikzpicture}
\end{gathered}
.
\end{gather}
The sandwiched algebra is the $R$-algebra spanned by $\dotidem$
and $y_{1}\dotidem$, and this algebra is isomorphic to $R[X]/(X^{2})$, with $X$ in degree four,
via $\dotidem\mapsto 1$ and $y_{1}\dotidem\mapsto X$.

\end{enumerate}
The graded dimension of
the set of all cyclotomic weighted KLRW
diagrams with $(0,1,2,2,0,1)$ as their residue sequence
at bottom and top is thus
$q^{2}(1+q^{-2})^{2}(1+q^{4})$,
where we use $q$ to keep track of the grading.
\end{Example}

\begin{Notation}\label{N:StrategyAffineRed}
Recall from \autoref{SS:RecollectionDefinition} that we will draw \emph{affine red strings} in diagrams, illustrated by:
\begin{gather*}
\text{genuine red string}:
\begin{tikzpicture}[anchorbase,smallnodes,rounded corners]
\draw[redstring] (0,0)node[below]{$i$} to (0,0.5)node[above,yshift=-1pt]{$\phantom{i}$};
\end{tikzpicture}
,\quad
\text{affine red string}:
\begin{tikzpicture}[anchorbase,smallnodes,rounded corners]
\draw[affine] (0,0)node[below]{$i$} to (0,0.5)node[above,yshift=-1pt]{$\phantom{i}$};
\end{tikzpicture}
.
\end{gather*}
Importantly, the affine red strings are not part of the
diagrams and they are drawn only as a visual aid.
\end{Notation}

The affine red strings are drawn to the right of all of the genuine red strings and we think of them as satisfying the same relations as the red strings. Our strategy is to pull strings as far to the right as possible but it is possible to pull some strings arbitrarily far to the right. The affine red strings are necessary to catch the solid strings that can be pulled past all of the genuine red strings in the diagram (we add enough affine red strings so as to be able to block all solid strings). By definition, a diagram is unsteady if and only if it contains a solid string that is blocked by an affine red string. As a vector space, the finite dimensional quotients of weighted KLRW algebras are spanned by those diagrams that do not have any strings blocked by affine red strings.

\begin{Example}\label{E:StrategyAffine}
To show why the affine red strings are necessary consider a diagram with one red $0$-string and one solid $1$-string. For simplicity we do not draw the ghost strings. Using the
honest Reidemeister II relation, the
solid $1$-string pulls arbitrarily far to the right:
\begin{gather*}
\begin{tikzpicture}[anchorbase,smallnodes,rounded corners]
\draw[solid](0,0)node[below]{$1$}--++(0,1)node[above,yshift=-1pt]{$\phantom{i}$};
\draw[redstring](0.25,0)node[below]{$0$}--++(0,1);
\end{tikzpicture}
=
\begin{tikzpicture}[anchorbase,smallnodes,rounded corners]
\draw[solid](0,0)node[below]{$1$}--++(0.5,0.5)--++(-0.5,0.5)node[above,yshift=-1pt]{$\phantom{i}$};
\draw[redstring](0.25,0)node[below]{$0$}--++(0,1);
\end{tikzpicture}
=
\begin{tikzpicture}[anchorbase,smallnodes,rounded corners]
\draw[solid](0,0)node[below]{$1$}--++(1.1,0.5)--++(-1.1,0.5)node[above,yshift=-1pt]{$\phantom{i}$};
\draw[redstring](0.25,0)node[below]{$0$}--++(0,1);
\end{tikzpicture}
=
\begin{tikzpicture}[anchorbase,smallnodes,rounded corners]
\draw[solid](0,0)node[below]{$1$}--++(2,0.5)--++(-2,0.5)node[above,yshift=-1pt]{$\phantom{i}$};
\draw[redstring](0.25,0)node[below]{$0$}--++(0,1);
\end{tikzpicture}
.
\end{gather*}
In order to place the $1$-string our strategy is to add an affine red $1$-string to the right of the red string:
\begin{gather*}
\begin{tikzpicture}[anchorbase,smallnodes,rounded corners]
\draw[solid](0,0)node[below]{$1$}--++(0,1)node[above,yshift=-1pt]{$\phantom{i}$};
\draw[redstring](0.25,0)node[below]{$0$}--++(0,1);
\draw[affine](1.25,0)node[below]{$1$}--++(0,1);
\end{tikzpicture}
=
\begin{tikzpicture}[anchorbase,smallnodes,rounded corners]
\draw[solid](0,0)node[below]{$1$}--++(0.5,0.5)--++(-0.5,0.5)node[above,yshift=-1pt]{$\phantom{i}$};
\draw[redstring](0.25,0)node[below]{$0$}--++(0,1);
\draw[affine](1.25,0)node[below]{$1$}--++(0,1);
\end{tikzpicture}
=
\begin{tikzpicture}[anchorbase,smallnodes,rounded corners]
\draw[solid](0,0)node[below]{$1$}--++(1.1,0.5)--++(-1.1,0.5)node[above,yshift=-1pt]{$\phantom{i}$};
\draw[redstring](0.25,0)node[below]{$0$}--++(0,1);
\draw[affine](1.25,0)node[below]{$1$}--++(0,1);
\end{tikzpicture}
.
\end{gather*}
The solid $1$-strings gets stuck at the affine red $1$-string because it does not satisfy an honest
Reidemeister~II relation. (We add enough affine red strings to ensure that all strings are blocked.)
\end{Example}

The following lemmas are the crucial diagrammatic
relations that we need to pull
strings and jump dots to the right. Here we are pulling the leftmost string to the right or jumping the leftmost  dot to the right. We highlight the strings where the action happens by coloring them.

\begin{Lemma}\label{L:StrategyMovingStringsDots}
For any quiver and any choice of $Q$-polynomials we have the following, plus partner relations:
\begin{gather}\label{Eq:StrategyTwoSolid}
\begin{tikzpicture}[anchorbase,smallnodes,rounded corners]
\draw[solid,spinach](0,1)node[above,yshift=-1pt]{$\phantom{i}$}--++(0,-1)node[below]{$i$};
\draw[solid](0.5,1)--++(0,-1)node[below]{$i$};
\end{tikzpicture}
=
\begin{tikzpicture}[anchorbase,smallnodes,rounded corners]
\draw[solid,dot](0,1)node[above,yshift=-1pt]{$\phantom{i}$}--++(0.8,-0.5)--++(-0.8,-0.5) node[below]{$i$};
\draw[solid,dot=0.1](0.5,1)--++(0,-1) node[below]{$i$};
\end{tikzpicture}
-
\begin{tikzpicture}[anchorbase,smallnodes,rounded corners]
\draw[solid,dot,dot=0.9](0,1)node[above,yshift=-1pt]{$\phantom{i}$}--++(0.8,-0.5)--++(-0.8,-0.5) node[below]{$i$};
\draw[solid](0.5,1)--++(0,-1) node[below]{$i$};
\end{tikzpicture}
,\quad
\begin{tikzpicture}[anchorbase,smallnodes,rounded corners]
\draw[solid,dot,spinach](0,1)node[above,yshift=-1pt]{$\phantom{i}$}--++(0,-1)node[below]{$i$};
\draw[solid](0.5,1)--++(0,-1)node[below]{$i$};
\end{tikzpicture}
=
\begin{tikzpicture}[anchorbase,smallnodes,rounded corners]
\draw[solid](0,1)node[above,yshift=-1pt]{$\phantom{i}$}--++(0,-1)node[below]{$i$};
\draw[solid,dot](0.5,1)--++(0,-1)node[below]{$i$};
\end{tikzpicture}
+
\begin{tikzpicture}[anchorbase,smallnodes,rounded corners]
\draw[solid,dot=0.42,dot=0.58](0,1)node[above,yshift=-1pt]{$\phantom{i}$}--++(0.8,-0.5)--++(-0.8,-0.5) node[below]{$i$};
\draw[solid,dot=0.1](0.5,1)--++(0,-1) node[below]{$i$};
\end{tikzpicture}
-
\begin{tikzpicture}[anchorbase,smallnodes,rounded corners]
\draw[solid,dot=0.42,dot=0.58,dot=0.9](0,1)node[above,yshift=-1pt]{$\phantom{i}$}--++(0.8,-0.5)--++(-0.8,-0.5)node[below]{$i$};
\draw[solid](0.5,1)--++(0,-1) node[below]{$i$};
\end{tikzpicture}
.
\end{gather}
For $i\to j$ edges and the choice of $Q$-polynomials in \autoref{Eq:RecollectionQPoly} we have the following, plus partner relations:
\begin{gather}\label{Eq:StrategyReidemeisterIIAndIIITypeA}
\begin{tikzpicture}[anchorbase,smallnodes,rounded corners]
\draw[ghost,dot,spinach](0,1)node[above,yshift=-1pt]{$\phantom{i}$}--++(0,-1)node[below]{$i$};
\draw[solid](0.5,1)--++(0,-1)node[below]{$j$};
\end{tikzpicture}
=
\begin{tikzpicture}[anchorbase,smallnodes,rounded corners]
\draw[ghost](0,1)node[above,yshift=-1pt]{$\phantom{i}$}--++(0.8,-0.5)--++(-0.8,-0.5) node[below]{$i$};
\draw[solid](0.5,1)--++(0,-1) node[below]{$j$};
\end{tikzpicture}
+
\begin{tikzpicture}[anchorbase,smallnodes,rounded corners]
\draw[ghost](0,1)node[above,yshift=-1pt]{$\phantom{i}$}--++(0,-1)node[below]{$i$};
\draw[solid,dot](0.5,1)--++(0,-1)node[below]{$j$};
\end{tikzpicture}
,\quad
\begin{tikzpicture}[anchorbase,smallnodes,rounded corners]
\draw[ghost](1,0)--++(0,1)node[above,yshift=-1pt]{$i$};
\draw[ghost](1.5,0)--++(0,1)node[above,yshift=-1pt]{$i$};
\draw[solid,spinach](0,0)node[below]{$i$}--++(0,1);
\draw[solid](0.5,0)node[below]{$i$}--++(0,1);
\draw[solid](1.25,0)node[below]{$j$}--++(0,1);
\end{tikzpicture}
=
-
\begin{tikzpicture}[anchorbase,smallnodes,rounded corners]
\draw[ghost](1,0)--++(0.5,0.5)--++(-0.5,0.5)node[above,yshift=-1pt]{$i$};
\draw[ghost,dot](1.5,0)--++(-0.5,0.5)--++(0.5,0.5)node[above,yshift=-1pt]{$i$};
\draw[solid](0,0)node[below]{$i$}--++(0.5,0.5)--++(-0.5,0.5);
\draw[solid,dot](0.5,0)node[below]{$i$}--++(-0.5,0.5)--++(0.5,0.5);
\draw[solid](1.25,0)node[below]{$j$}--++(-0.5,0.5)--++(0.5,0.5);
\end{tikzpicture}
-
\begin{tikzpicture}[anchorbase,smallnodes,rounded corners]
\draw[ghost,dot](1,0)--++(0.5,0.5)--++(-0.5,0.5)node[above,yshift=-1pt]{$i$};
\draw[ghost](1.5,0)--++(-0.5,0.5)--++(0.5,0.5)node[above,yshift=-1pt]{$i$};
\draw[solid,dot](0,0)node[below]{$i$}--++(0.5,0.5)--++(-0.5,0.5);
\draw[solid](0.5,0)node[below]{$i$}--++(-0.5,0.5)--++(0.5,0.5);
\draw[solid](1.25,0)node[below]{$j$}--++(0.5,0.5)--++(-0.5,0.5);
\end{tikzpicture}
.
\end{gather}
The right relation and its partner also hold for $i\Leftarrow j$.

For $i\Rightarrow j$ edges and the choice of $Q$-polynomials in \autoref{Eq:RecollectionQPoly} we have the following, plus partner relations:
\begin{gather}\label{Eq:StrategyReidemeisterIIAndIIITypesBAD}
\begin{tikzpicture}[anchorbase,smallnodes,rounded corners]
\draw[ghost,dot,spinach](0,1)node[above,yshift=-1pt]{$\phantom{i}$}--++(0,-1)node[below]{$i$};
\draw[solid](0.5,1)--++(0,-1)node[below]{$j$};
\end{tikzpicture}
=
\begin{tikzpicture}[anchorbase,smallnodes,rounded corners]
\draw[ghost](0,1)node[above,yshift=-1pt]{$\phantom{i}$}--++(0.8,-0.5)--++(-0.8,-0.5) node[below]{$i$};
\draw[solid](0.5,1)--++(0,-1) node[below]{$j$};
\end{tikzpicture}
+
\begin{tikzpicture}[anchorbase,smallnodes,rounded corners]
\draw[ghost](0,1)node[above,yshift=-1pt]{$\phantom{i}$}--++(0,-1)node[below]{$i$};
\draw[solid,dot=0.4,dot=0.6](0.5,1)--++(0,-1)node[below]{$j$};
\end{tikzpicture}
.
\end{gather}
\end{Lemma}

\begin{proof}
As in \cite[Lemmas 6D.1, 6D.4 and 7E.1]{MaTu-klrw-algebras}.
\end{proof}

The relations \autoref{Eq:StrategyReidemeisterIIAndIIITypeA}
and \autoref{Eq:StrategyReidemeisterIIAndIIITypesBAD}
motivate the dot placement in
\autoref{D:BasesSandwichpart} below: we construct our basis diagrams so that strings cannot be pulled, and dots cannot be jumped, further to the right, so we need to ensure that these relations cannot be applied.

\begin{Lemma}\label{L:StrategyQuadruple}
For $i\not\!\Leftarrow j$ and the choice of $Q$-polynomials in \autoref{Eq:RecollectionQPoly} we have the following. In any close situation of the form
\begin{gather*}
\begin{tikzpicture}[anchorbase,smallnodes,rounded corners]
\draw[ghost](0,0)--++(0,1)node[above,yshift=-1pt]{$j$};
\draw[solid,spinach](-0.4,0)node[below]{$i$}--++(0,1);
\draw[solid](-0.2,0)node[below]{$i$}--++(0,1);
\end{tikzpicture}
,
\end{gather*}
we can pull the marked string further to the right. Similarly for its partner relations.
\end{Lemma}

\begin{proof}
We first use \autoref{Eq:StrategyTwoSolid} to
pull the left solid $i$-string into the middle.
Note that this string also carries a dot after applying
\autoref{Eq:StrategyTwoSolid}.
If $i\not\!\Leftarrow j$ are not connected, then either an honest Reidemeister II relation applies, or either \autoref{Eq:StrategyReidemeisterIIAndIIITypeA}
or \autoref{Eq:StrategyReidemeisterIIAndIIITypesBAD} applies. In all cases, we can pull the leftmost solid $i$-string to the right.
\end{proof}

The next example should be compared with \autoref{R:BasesCombinatoricsDottedIdempotent} below.

\begin{Example}\label{E:StrategyStuck}
The following close configurations and their partners are stuck for $i\Rightarrow j$ and $i\Leftarrow j$, respectively, so that \autoref{L:StrategyMovingStringsDots} and \autoref{L:StrategyQuadruple}
do not apply:
\begin{gather*}
i\Rightarrow j\colon
\begin{tikzpicture}[anchorbase,smallnodes,rounded corners]
\draw[ghost](1,0)--++(0,1)node[above,yshift=-1pt]{$j$};
\draw[ghost](1.5,0)--++(0,1)node[above,yshift=-1pt]{$j$};
\draw[solid](1.25,0)node[below]{$i$}--++(0,1);
\end{tikzpicture}
,\quad
i\Leftarrow j\colon
\begin{tikzpicture}[anchorbase,smallnodes,rounded corners]
\draw[ghost](1,0)--++(0,1)node[above,yshift=-1pt]{$j$};
\draw[ghost](1.75,0)--++(0,1)node[above,yshift=-1pt]{$j$};
\draw[solid](1.25,0)node[below]{$i$}--++(0,1);
\draw[solid](1.5,0)node[below]{$i$}--++(0,1);
\end{tikzpicture}
.
\end{gather*}
These diagrams are therefore good in the sense that we can not pull strings or jump dots further to the right. These configurations will appear whenever $i$ corresponds to a leaf of $\Gamma$.
\end{Example}


\section{The bases}\label{S:Bases}


We now explain the main constructions
of this paper. In \autoref{R:ProofGeneral} we summarize the
parts of the arguments that are general and those that
depend on the underlying quiver.


\subsection{The basis in a nutshell}\label{SS:BasesGeneralStrategy}


The main aim of the paper is the construction of the elements in \autoref{Eq:BasesTheBasis} which form a homogeneous (affine)
sandwich cellular basis for the weighted KLRW algebras.

The construction in \autoref{R:SandwichCellularAlgebra} is quite general and works for many
algebras: we sandwich a class of
diagrams between another class of diagrams.
In more detail, the relevant picture for \autoref{Eq:BasesTheBasis} is:
\begin{gather*}
D_{\bS\bT}^{\ba,\bfi}
\leftrightsquigarrow
\begin{tikzpicture}[anchorbase,scale=1]
\draw[line width=0.75,color=black,fill=cream] (0,-0.5) to (0.25,0) to (0.75,0) to (1,-0.5) to (0,-0.5);
\node at (0.5,-0.25){$\bT$};
\draw[line width=0.75,color=black,fill=cream] (0,1) to (0.25,0.5) to (0.75,0.5) to (1,1) to (0,1);
\node at (0.5,0.75){$\bS$};
\draw[line width=0.75,color=black,fill=cream] (0.25,0) to (0.25,0.5) to (0.75,0.5) to (0.75,0) to (0.25,0);
\node at (0.5,0.25){\scalebox{0.7}{$\sandwich{\blam}{\ba}{\bfi}$}};
\end{tikzpicture}
=
\begin{tikzpicture}[anchorbase,scale=1]
\draw[line width=0.75,color=black,fill=cream] (0,-0.5) to (0.25,0) to (0.75,0) to (1,-0.5) to (0,-0.5);
\node at (0.5,-0.25){$\bT$};
\draw[line width=0.75,color=black,fill=cream] (0,2) to (0.25,1.5) to (0.75,1.5) to (1,2) to (0,2);
\node at (0.5,1.75){$\bS$};
\draw[line width=0.75,color=black,fill=cream] (0.25,0) to (0.25,0.5) to (0.75,0.5) to (0.75,0) to (0.25,0);
\node at (0.5,0.25){$\dotidem$};
\draw[line width=0.75,color=black,fill=cream] (0.25,0.5) to (0.25,1) to (0.75,1) to (0.75,0.5) to (0.25,0.5);
\node at (0.5,0.75){$\zeetwo$};
\draw[line width=0.75,color=black,fill=cream] (0.25,1) to (0.25,1.5) to (0.75,1.5) to (0.75,1) to (0.25,1);
\node at (0.5,1.25){$\daffine$};
\end{tikzpicture},
\quad
\begin{aligned}
\begin{tikzpicture}[anchorbase,scale=1]
\draw[line width=0.75,color=black,fill=cream] (0,1) to (0.25,0.5) to (0.75,0.5) to (1,1) to (0,1);
\node at (0.5,0.75){$\bS$};
\end{tikzpicture}
&\text{ a permutation diagram,}
\\
\begin{tikzpicture}[anchorbase,scale=1]
\draw[white,ultra thin] (0,0) to (1,0);
\draw[line width=0.75,color=black,fill=cream] (0.25,0) to (0.25,0.5) to (0.75,0.5) to (0.75,0) to (0.25,0);
\node at (0.5,0.25){$\daffine$};
\end{tikzpicture}
&\text{ an ``affine'' dot placement,}
\\
\begin{tikzpicture}[anchorbase,scale=1]
\draw[white,ultra thin] (0,0) to (1,0);
\draw[line width=0.75,color=black,fill=cream] (0.25,0) to (0.25,0.5) to (0.75,0.5) to (0.75,0) to (0.25,0);
\node at (0.5,0.25){$\zeetwo$};
\end{tikzpicture}
&\text{ a ``finite'' dot placement,}
\\
\begin{tikzpicture}[anchorbase,scale=1]
\draw[white,ultra thin] (0,0) to (1,0);
\draw[line width=0.75,color=black,fill=cream] (0.25,0) to (0.25,0.5) to (0.75,0.5) to (0.75,0) to (0.25,0);
\node at (0.5,0.25){$\dotidem$};
\end{tikzpicture}
&\text{ a dotted idempotent,}
\\
\begin{tikzpicture}[anchorbase,scale=1]
\draw[line width=0.75,color=black,fill=cream] (0,-0.5) to (0.25,0) to (0.75,0) to (1,-0.5) to (0,-0.5);
\node at (0.5,-0.25){$\bT$};
\end{tikzpicture}
&\text{ a permutation diagram.}
\end{aligned}
\end{gather*}
See also \autoref{Eq:StrategyMainBig} which gives an explicit factorization in this form. The affine dots do not appear in the basis for the finite dimensional cyclotomic quotients and, in this case, only steady diagrams appear.
The construction of the sandwich cellular basis for the $\bad$ types is almost exactly the same as the construction for the $\ac$ types in \cite{MaTu-klrw-algebras} except that the finite dot placement is trivial for the $\ac$ types.

Below we explain the ingredients needed to construct these bases: the partition combinatorics, which give the indexing poset, the dotted idempotents, which gives $\dotidem$,
the finite and affine dot placements, giving
$\zeetwo$ and $\daffine$, and permutation diagrams.


\subsection{Some notation}\label{SS:BasesSomeNotation}


We fix some conventions.

\begin{Notation}\label{N:BasesAffineHell}
Unless we are in specific example, for the duration of the paper we fix
arbitrary $n,\ell\in\N$, $e\in\Z_{\geq 2}$,
$\charge\in\Z^{\ell}$, with $\kappa_{1}<\dots<\kappa_{\ell}$, and $\brho\in I^{\ell}$. We let
\begin{gather*}
\hell=\ell+n(e+1)
\end{gather*}
be the \emph{affine level}, which counts the maximum number of genuine and affine red strings in the diagram. More generally, we will use the
underline notation to indicate definitions that only play a role in the
affine case.
\end{Notation}

\begin{Notation}\label{N:BasesFunnyQuivers1}
Our constructions given in this section
work for the quivers below.
\begin{enumerate}

\item We use Kac's notation \cite{Ka-infdim-lie} for Dynkin quivers
(but we reflected the quivers left-to-right).
The main quivers of study in this paper are:
\begin{gather}\label{Eq:BasesFunnyQuivers}
\begin{aligned}
\typeb:\quad&
\begin{tikzpicture}[anchorbase]
\draw[directed=0.5,double,double distance=0.5mm](1,0)--(0,0);
\draw(1,0)--(6,0);
\node at (6.5,0){$\cdots$};
\foreach \x in {0,...,6} {
\node[circle,inner sep=1.8pt,fill=DarkBlue] (\x) at (\x,0){};
\node at (\x,-0.25){$\x$};
}
\end{tikzpicture},
\\[1mm]
\atwo:\quad&
\begin{tikzpicture}[anchorbase]
\draw[directed=0.5,double,double distance=0.5mm](0,0)--(1,0);
\draw[directed=0.5,double,double distance=0.5mm](5,0)--(6,0);
\draw(1,0)--(3,0);
\draw(4,0)--(5,0);
\node at (3.5,0){$\cdots$};
\foreach \x in {0,...,3} {
\node[circle,inner sep=1.8pt,fill=DarkBlue] (\x) at (\x,0){};
\node at (\x,-0.25){$\x$};
}
\foreach \x [evaluate=\x as \c using {int(6-\x)}] in {4,5} {
\node[circle,inner sep=1.8pt,fill=DarkBlue] (\x) at (\x,0){};
\node at (\x,-0.25){$e-\c$};
}
\node[circle,inner sep=1.8pt,fill=DarkBlue] (6) at (6,0){};
\node at (6,-0.25){$e$};
\end{tikzpicture},
\\[1mm]
\dtwo:\quad&
\begin{tikzpicture}[anchorbase]
\draw[directed=0.5,double,double distance=0.5mm](1,0)--(0,0);
\draw[directed=0.5,double,double distance=0.5mm](5,0)--(6,0);
\draw(1,0)--(3,0);
\draw(4,0)--(5,0);
\node at (3.5,0){$\cdots$};
\foreach \x in {0,...,3} {
\node[circle,inner sep=1.8pt,fill=DarkBlue] (\x) at (\x,0){};
\node at (\x,-0.25){$\x$};
}
\foreach \x [evaluate=\x as \c using {int(6-\x)}] in {4,5} {
\node[circle,inner sep=1.8pt,fill=DarkBlue] (\x) at (\x,0){};
\node at (\x,-0.25){$e-\c$};
}
\node[circle,inner sep=1.8pt,fill=DarkBlue] (6) at (6,0){};
\node at (6,-0.25){$e$};
\end{tikzpicture}
,
\end{aligned}
\end{gather}
where $e\in\Z_{\geq 2}$.
Up to minor changes in the notation, we can fix any orientation of the simply laced edges. Here 
we orient the simply laced edges $i\to(i+1)$ as this fits better with classical partition combinatorics.
We will omit $\typeb$ from the discussion because
this type can be viewed as $\dtwo$ for $e\to\infty$, but doing so needs
some (harmless) adjustments of the exposition since {\eg} the affine level becomes arbitrarily large.

\item Let $0<\varepsilon<\frac{1}{4n\hell}$ be a small shift.
We let the ghost shifts for edges in $\bad$ types be $1$ with the exception
of the edge $0\Leftarrow 1$ in type $\dtwo$ where the ghost shift
is $1-\varepsilon^{2}$.

\end{enumerate}
\end{Notation}

\begin{Definition}\label{D:BasesSinks}
A \emph{sink} is a vertex of a Dynkin quiver $\Gamma$ that is a
(graph-theoretical) sink. That is, the vertex is not the tail of any edge. A \emph{multisink} is a sink such that all edges are multi-laced.
\end{Definition}

Note that, if $i$ is a sink, then the solid $i$-strings do not have ghosts.

\begin{Example}\label{E:BasesSinks}
The vertex $e$ is a multisink for $\bad$ types and additionally, the vertex
$0$ is a multisink for type $\dtwo$. In contrast,
in type $\cone$ no vertex is a multisink, but some vertices are sinks (depending on the choice of orientation).
\end{Example}

In type $\dtwo$ the solid $1$-string has two ghosts that are very close
to one another by our choice of ghost shifts: the ghost string for the edge $0\Leftarrow 1$ has ghost shift $1-\varepsilon^{2}$ and the ghost string for the edge $1\to 2$ has ghost shift $1$. As the $1$-ghost strings are very close, we display
them as doubled lines:
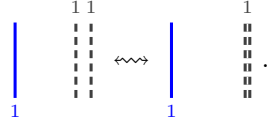
\begin{gather*}
\begin{tikzpicture}[anchorbase,smallnodes,rounded corners]
\draw[ghost](-0.2,0)--++(0,1)node[above,yshift=-1pt]{$1$};
\draw[ghost](0,0)--++(0,1)node[above,yshift=-1pt]{$1$};
\draw[solid](-1,0)node[below]{$1$}--++(0,1);
\end{tikzpicture}
\leftrightsquigarrow
\begin{tikzpicture}[anchorbase,smallnodes,rounded corners]
\draw[dghost](0,0)node[below]{$\phantom{i}$}--++(0,1)node[above,yshift=-1pt]{$1$};
\draw[solid](-1,0)node[below]{$1$}--++(0,1);
\end{tikzpicture}
.
\end{gather*}
We stress that these are two different ghost $1$-strings; {\cf} \autoref{R:RecollectionDifferentStrings}. As in \autoref{E:StrategyMain}, when these strings have dots, we put a single ghost dot on the doubled ghost string rather than a 
dot on each ghost string.

We will draw generic diagrams that are supposed to make sense in any
type, but the reader may need to remove or double some ghost strings to obtain the actual diagram for a particular type.

\begin{Definition}\label{D:BasesAffineNotation}
Define the \emph{affine charge}
$\affine{\charge}=
(\affine{\kappa}_{1},\dots,\affine{\kappa}_{\hell})\in\Z^{\hell}$
and the \emph{affine red labels}
$\affine{\brho}=(\affine{\rho}_{1},\dots,\affine{\rho}_{\hell})\in I^{\hell}$ by
\begin{gather*}
\affine{\kappa}_{m}=
\begin{cases*}
\kappa_{m} & if $1\leq m\leq\ell$,\\
\kappa_{\ell}+2n(m-\ell) & otherwise,
\end{cases*}
\quad\text{and}\quad
\affine{\rho}_{m}=
\begin{cases*}
\rho_{m} & if $1\leq m\leq\ell$,\\
\floor{m-\ell-1}{n}+(e+1)\Z & otherwise.
\end{cases*}
\end{gather*}
We call $\affine{\kappa}_{m}$ and $\affine{\rho}_{m}$ the \emph{position and residue} of a red string for $m\leq\ell$, and the \emph{position and residue} of an affine red string for $m>\ell$.
\end{Definition}

Note that the coordinates of the (affine) red strings
$\affine{\charge}$ are always integers.

\begin{Example}\label{E:BasesBasisExemplifiedPartI}
Take $n=3$, $e=2$ and $\ell=1$, so $\hell=1+3\cdot(2+1)=10$.
If $\charge=(2)$ and $\brho=(1)$, then $\affine{\charge}=(2,8,14,20,26,32,38,44,50,56)$
and $\affine{\brho}=(1,0,0,0,1,1,1,2,2,2)$. All entries, of $\affine{\charge}$ and of $\affine{\rho}$, except the first are affine.

Below we will use this data to describe a diagram with red strings and affine 
red strings for the reason explained in \autoref{E:StrategyAffine}. For example, the data above will give
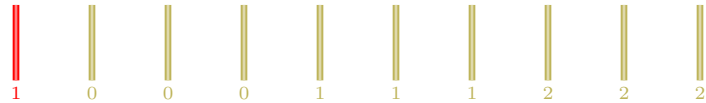
\begin{gather*}
\begin{tikzpicture}[anchorbase,smallnodes,rounded corners]
\draw[redstring](-1,0)node[below]{$1$}--++(0,1)node[above,yshift=-1pt]{$\phantom{i}$};
\draw[affine](0,0)node[below]{$0$}--++(0,1);
\draw[affine](1,0)node[below]{$0$}--++(0,1);
\draw[affine](2,0)node[below]{$0$}--++(0,1);
\draw[affine](3,0)node[below]{$1$}--++(0,1);
\draw[affine](4,0)node[below]{$1$}--++(0,1);
\draw[affine](5,0)node[below]{$1$}--++(0,1);
\draw[affine](6,0)node[below]{$2$}--++(0,1);
\draw[affine](7,0)node[below]{$2$}--++(0,1);
\draw[affine](8,0)node[below]{$2$}--++(0,1);
\end{tikzpicture}
\end{gather*}
In this example the positions of the strings are scaled to make it fit
onto the page, the real positions are given by $\affine{\charge}=(2,8,14,20,26,32,38,44,50,56)$. Note that $n=3$ and if we add $n$ solid strings to this diagram from the left, and their ghosts, then we will see that each of these strings gets stuck somewhere.
\end{Example}


\subsection{Partition combinatorics}\label{SS:BasesCombinatorics}


Before coming to our main definitions, we introduce
the tableau combinatorics
that arise in $\bad$ types. For the
standard tableau combinatorics that appears in the context
of KLR algebras we refer the reader to \cite[Section 3.3]{HuMa-klr-basis}.

\begin{Remark}\label{R:BasesTableaux}
The partition
combinatorics that we use is motivated by {\eg} \cite{ArPa-quiver-hecke-a2}
and \cite{ArPa-quiver-hecke-d} who use the
partitions of the crystal graphs for the affine fundamental weight. The associated weighted KLRW diagram combinatorics is a slight modification of
the combinatorics of type $\cone$ as in \cite[Section 7]{MaTu-klrw-algebras}.

The reader should be careful because, as we will see, the partition combinatorics depends on $\brho$. In particular, the combinatorics from \cite{ArPa-quiver-hecke-a2} and \cite{ArPa-quiver-hecke-d} only applies for cyclotomic KLR algebras for the fundamental weights associated to multisink vertices. This is related to the fact that the crystal graphs for fundamental weights indexed by
multisink vertices have different indexing sets to the crystal graphs
for the other fundamental weights.
\end{Remark}

Identify the vertices of $I$ with $\Z/(e+1)\Z$.
We use \emph{(usual) partitions} of $n$, following the
same conventions as in \cite[Section 6A]{MaTu-klrw-algebras}.
We also use \emph{shifted partitions} of $n$,
that is, partitions $\lambda$
with strictly decreasing parts $\lambda=(\lambda_{1}>\dots>\lambda_{k})$. (Shifted partitions are called \emph{strict} partitions by some authors.) Let $|\lambda|$ be the size of a partition or shifted partition.

\begin{Definition}\label{D:BasesCombinatoricsTableauxShifted}
The set of \emph{$\brho$-partitions} is
\begin{gather*}
\Parts=\Biggl\{(\lambda^{(1)}|\dots|\lambda^{(\ell)})\Biggm|
\begin{array}{l}|\lambda^{(1)}|+\dots+|\lambda^{(\ell)}|=n, \text{ where}
\\
\lambda^{(i)}\text{ is a usual partition if $\rho_{i}$ is not a multisink}
\\
\lambda^{(i)}\text{ is a shifted partition if $\rho_{i}$ is a multisink}
\end{array}\Biggr\}
.
\end{gather*}
We define $\hParts$ similarly using the indicated data.
\end{Definition}

We identify a $\brho$-partition $\blam$ with its
\emph{shifted $\brho$-Young diagram}, which is the set of nodes
$\set{(m,r,c)}$ where $1\leq m\leq\ell$, $\lambda^{(m)}_{r}>0$, and
\begin{gather*}
\begin{cases*}
1\leq c\leq\lambda^{(m)}_{r} & if $\rho_{m}$ is not a multisink,
\\
1+r\leq c\leq\lambda^{(m)}_{r}+r & if $\rho_{m}$ is a multisink.
\\
\end{cases*}
\end{gather*}

\begin{Remark}\label{R:BasesOrderOfMRC}
In this paper a node $(m,r,c)$ is in component $m$, row $r$ and column $c$. This is different to the convention of {\eg} \cite{DiJaMa-cyclotomic-q-schur}, where nodes are indexed in the order $(r,c,m)$. Our convention is better for this paper because in the representation theory the nodes are naturally ordered first by component and then by their row and column index, which is the lexicographic order on the nodes using our convention.
\end{Remark}

\begin{Notation}\label{N:BasesCombinatoricsTableaux}
We use the (shifted) English convention to illustrate
the associated \emph{(shifted) $\brho$-Young diagrams}. That is, we illustrate
these partitions by drawing them as boxes in the plane, with rows ordered from top to bottom, and columns left to right, and where the $r$th row is shifted $r$ positions to the right for shifted $\brho$-Young diagrams. (So the
$c$th column of the $r$th row is in position $c+r$.)
For example, a usual Young diagram and a shifted Young diagram are
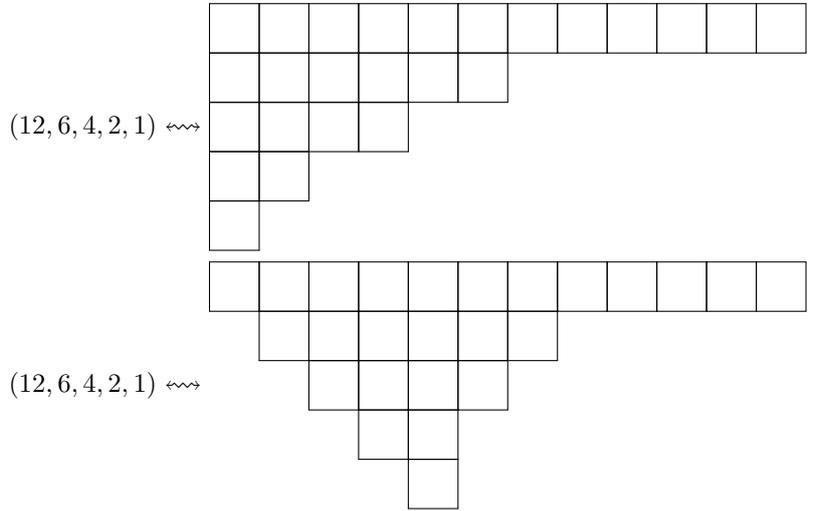
\begin{gather*}
(12,6,4,2,1)\leftrightsquigarrow
\ydiagram{12,6,4,2,1}
\;,
\\
(12,6,4,2,1)\leftrightsquigarrow
\ydiagram{12,1+6,2+4,3+2,4+1}
\;.
\end{gather*}
This is a different convention to that used
in \cite[Sections 6 and 7]{MaTu-klrw-algebras} where
the Russian convention is used. The Russian convention is
useful in types $\typea$ and $\aone$, but seems to be irrelevant in other types.
\end{Notation}

That is, if $\rho_{i}$
corresponds to a multisink vertex ({\ie} the vertices $e$ in
all $\bad$ types
and, additionally, the vertex $0$ in type $\dtwo$), then we consider shifted
partitions in the $i$th entry of $\blam$, and usual partitions otherwise.
We make similar definitions as above for $\affine{\brho}$-partitions.
The nodes for $m>\ell$ are called \emph{affine nodes}.
Identify $\Parts$ with the subset of $\hParts$ that consists of the $\affine{\brho}$ that do not contain any affine nodes.

\begin{Definition}\label{D:BasesResidues}
Let $a=1$ for type $\atwo$ and $a=2$ for type $\dtwo$.
Given an integer $k\in\Z$ we need the remainder of division
by $2e+a$ (viewed as an element of $\set{0,\dots,2e+a}\subset\N$),
which we denote by $\mathrm{Mod}(k,2e+a)$. Define the \emph{residue function} $\re\colon\Z\longrightarrow I$ by
\begin{gather*}
\re(k)=
\begin{cases*}
\mathrm{Mod}(k,2e+a)+(e+1)\Z & if $0\leq\mathrm{Mod}(k,2e+a)\leq e$,\\
2e+1-\mathrm{Mod}(k,2e+a)+(e+1)\Z & if $e<\mathrm{Mod}(k,2e+a)<2e+a$.
\end{cases*}
\end{gather*}
The \emph{($\affine{\brho}$-)residue} of the node $(m,r,c)$ is $\res(m,r,c)=\re(c-r)+\rho_{m}+(e+1)\Z$.
\end{Definition}

In illustrations we will often fill nodes with their residues.
As the next example shows, the formal definition is more complicated than the actual concept of a residue.

\begin{Example}\label{E:BasesCombinatoricsTableaux}
For $\blam=(10,9,8,7,6,5,4)$ and $e=3$, starting with $3$,
{\ie} $\brho=(3)$, we get:
\begin{gather*}
\atwo[3]:
\begin{ytableau}
*(spinach!50) 3 & 2 & 1 & 0 & 1 & 2 & *(spinach!50) 3 & *(spinach!50) 3 & 2 & 1 \\
\none & *(spinach!50) 3 & 2 & 1 & 0 & 1 & 2 & *(spinach!50) 3 & *(spinach!50) 3 & 2 \\
\none & \none & *(spinach!50) 3 & 2 & 1 & 0 & 1 & 2 & *(spinach!50) 3 & *(spinach!50) 3 \\
\none & \none & \none & *(spinach!50) 3 & 2 & 1 & 0 & 1 & 2 & *(spinach!50) 3 \\
\none & \none & \none & \none & *(spinach!50) 3 & 2 & 1 & 0 & 1 & 2 \\
\none & \none & \none & \none & \none & *(spinach!50) 3 & 2 & 1 & 0 & 1 \\
\none & \none & \none & \none & \none & \none & *(spinach!50) 3 & 2 & 1 & 0 \\
\end{ytableau}
\;,\quad
\dtwo[3]:
\begin{ytableau}
*(spinach!50) 3 & 2 & 1 & *(tomato!50) 0 & *(tomato!50) 0 & 1 & 2 & *(spinach!50) 3 & *(spinach!50) 3 & 2 \\
\none & *(spinach!50) 3 & 2 & 1 & *(tomato!50) 0 & *(tomato!50) 0 & 1 & 2 & *(spinach!50) 3 & *(spinach!50) 3 \\
\none & \none & *(spinach!50) 3 & 2 & 1 & *(tomato!50) 0 & *(tomato!50) 0 & 1 & 2 & *(spinach!50) 3 \\
\none & \none & \none & *(spinach!50) 3 & 2 & 1 & *(tomato!50) 0 & *(tomato!50) 0 & 1 & 2 \\
\none & \none & \none & \none & *(spinach!50) 3 & 2 & 1 & *(tomato!50) 0 & *(tomato!50) 0 & 1 \\
\none & \none & \none & \none & \none & *(spinach!50) 3 & 2 & 1 & *(tomato!50) 0 & *(tomato!50) 0 \\
\none & \none & \none & \none & \none & \none & *(spinach!50) 3 & 2 & 1 & *(tomato!50) 0 \\
\end{ytableau}
\;.
\end{gather*}
Note that colored/shaded residues, for multisinks, are doubled.
In words,
the residues increase along rows and columns until they hit $e$,
this value is doubled, and then
the residues bounce back until they hit $0$,
this value is doubled for $\dtwo$, and the process starts again.

In contrast, if $\brho=(1)$, then
\begin{gather*}
\atwo[3]:
\begin{ytableau}
1 & 2 & *(spinach!50) 3 & *(spinach!50) 3 & 2 & 1 & 0 & 1 & 2 & *(spinach!50) 3\\
0 & 1 & 2 & *(spinach!50) 3 & *(spinach!50) 3 & 2 & 1 & 0 & 1 \\
1 & 0 & 1 & 2 & *(spinach!50) 3 & *(spinach!50) 3 & 2 & 1 \\
2 & 1 & 0 & 1 & 2 & *(spinach!50) 3 & *(spinach!50) 3 \\
*(spinach!50) 3 & 2 & 1 & 0 & 1 & 2 \\
*(spinach!50) 3 & *(spinach!50) 3 & 2 & 1 & 0\\
2 & *(spinach!50) 3 & *(spinach!50) 3 & 2\\
\end{ytableau}
\;,\quad
\dtwo[3]:
\begin{ytableau}
1 & 2 & *(spinach!50) 3 & *(spinach!50) 3 & 2 & 1 & *(tomato!50) 0 & *(tomato!50) 0 & 1 & 2\\
*(tomato!50) 0 & 1 & 2 & *(spinach!50) 3 & *(spinach!50) 3 & 2 & 1 & *(tomato!50) 0 & *(tomato!50) 0 \\
*(tomato!50) 0 & *(tomato!50) 0 & 1 & 2 & *(spinach!50) 3 & *(spinach!50) 3 & 2 & 1 \\
1 & *(tomato!50) 0 & *(tomato!50) 0 & 1 & 2 & *(spinach!50) 3 & *(spinach!50) 3 \\
2 & 1 & *(tomato!50) 0 & *(tomato!50) 0 & 1 & 2 \\
*(spinach!50) 3 & 2 & 1 & *(tomato!50) 0 & *(tomato!50) 0\\
*(spinach!50) 3 & *(spinach!50) 3 & 2 & 1\\
\end{ytableau}
\;,
\end{gather*}
following \autoref{D:BasesCombinatoricsTableauxShifted}. For $\brho=(2)$ the situation is similar. For $\brho=(0)$ type $\atwo[3]$ has partition, and type $\dtwo[3]$ has shifted partition, combinatorics. We emphasize that, in general, there are more partitions than shifted partitions, so there are more of these diagrams when $\affine{\rho}_{i}$ is not a multisink.
\end{Example}

The coloring/shading from \autoref{E:BasesCombinatoricsTableaux} distinguishes between usual and shifted partition combinatorics: we use shifted partition combinatorics if and only if the first node is colored/shaded.

\begin{Example}\label{E:BasesBasisExemplifiedPartII}
We continue with \autoref{E:BasesBasisExemplifiedPartI}
and fix type $\atwo[2]$. There are
three usual and two shifted $1$-partitions of $3$
and eight in total, namely (filled with their residues):
\begin{gather*}
\lambda_{1}=
\begin{ytableau}
0 & 1 & *(spinach!50) 2
\end{ytableau}
\,,\quad
\lambda_{2}=
\begin{ytableau}
0 & 1 \\
1
\end{ytableau}
\,,\quad
\lambda_{3}=
\begin{ytableau}
0 \\ 1 \\ *(spinach!50) 2
\end{ytableau}
\,,
\\
\mu_{1}=
\begin{ytableau}
1 & *(spinach!50) 2 & *(spinach!50) 2
\end{ytableau}
\,,\quad
\mu_{2}=
\begin{ytableau}
1 & *(spinach!50) 2 \\
0
\end{ytableau}
\,,\quad
\mu_{3}=
\begin{ytableau}
1 \\ 0 \\ 1
\end{ytableau}
\,,\\
\nu_{1}=
\begin{ytableau}
*(spinach!50) 2 & 1 & 0
\end{ytableau}
\,\quad
\nu_{2}=
\begin{ytableau}
*(spinach!50) 2 & 1 \\
\none & *(spinach!50) 2
\end{ytableau}
\,.
\end{gather*}
We get $\Parts=\set{\lambda_{1},\lambda_{2},\lambda_{3}}$,
if $\brho=(0)$, $\Parts=\set{\mu_{1},\mu_{2},\mu_{3}}$,
if $\brho=(1)$, and $\Parts=\set{\nu_{1},\nu_{2}}$,
for $\brho=(2)$.

The set $\hPartsall$ is much bigger, and we will not list it here. Note, however, that the affine components of $\hPartsall$ have either usual or shifted partition combinatorics depending in the residue. More precisely, since $\affine{\brho}=(\rho,0,0,0,1,1,1,2,2,2)$ we get that the first six affine components of $\hPartsall$ use usual partitions and the last three components use shifted partition combinatorics.
\end{Example}

\begin{Remark}\label{R:BasesShifted}
In examples, and in the proofs below, we will mostly focus on the shifted partition combinatorics.
The arguments for the usual partition combinatorics requires only slight adjustments and is essentially identical to that used for type $\cone$ in \cite[Section 7]{MaTu-klrw-algebras}.
\end{Remark}


\subsection{The idempotents}\label{SS:BasesPositioning}


We now introduce some crucial definitions for defining the idempotent diagrams. For any
$r+(e+1)\Z\in I$, with $r\in\set{0,\dots,e}$, we sometimes abuse notation and identify $r+(e+1)\Z$
with the real number in $r\in\R$.
The meaning of the following definition is clarified in \autoref{R:BasesCombinatoricsPosition} below.

\begin{Definition}\label{D:BasesContent}
Define the \textbf{content function} (for the $m$th component) $\co^{m}\colon I\to\R$ by
\begin{gather*}
\co^{m}(r)=
\begin{cases}
2-\rho_{m}&\text{if $\rho_{m}\neq 0$, $r=0$ and we are in type $\dtwo$},
\\
r-2&\text{if $\rho_{m}=0$, $r\neq 0$ and we are in type $\dtwo$},
\\
r-\rho_{m}&\text{otherwise}.
\end{cases}
\end{gather*}
\end{Definition}

Recall the shift $\varepsilon$ from \autoref{N:BasesFunnyQuivers1}.
Below we define the \emph{positioning function},
which uses the \emph{row reading order} $o_{\blam}(m,r,c)=c+\sum_{i=1}^{r-1}\lambda^{(m)}_{i}$, which is the lexicographic order on the nodes.

\begin{Example}\label{E:BasesContent}
The function $o_{\blam}(\placeholder)$
returns the position of a node in a Young diagram
when reading along rows,
{\ie} reading first left to right, and then bottom to top.
So $o_{\blam}(m,r,c)=k$ if $(m,r,c)$ is the $k$th node in row reading order. This is be illustrated by
\begin{gather*}
\left(
\begin{ytableau}
1 & 2 & 3
\end{ytableau}
\,,
\begin{ytableau}
4 & 5 \\
6
\end{ytableau}
\,,
\begin{ytableau}
7 \\ 8 \\ 9
\end{ytableau}
\,\right),
\end{gather*}
where we filled the nodes with their value $o_{\blam}(m,r,c)$.
\end{Example}

\begin{Notation}\label{N:BasesContent}
From now on we always order nodes, and strings associated to nodes, by the row reading order.
\end{Notation}

\begin{Definition}\label{D:BasesCombinatoricsPosition}
Let $\blam\in\hPartsall$. The \emph{coordinate} of the node $(m,r,c)\in\blam$
is
\begin{gather*}
\hcoord(m,r,c)=
\affine{\kappa}_{m}-\tfrac{m}{\hell}+
\co^{m}\big(\re(k)\big)-k\varepsilon,\qquad
\text{where $k=o_{\blam}(m,r,c)$}.
\end{gather*}
The \emph{coordinates} $\hcoord(\blam)$ of $\blam$
is the ordered tuple
of the coordinates of its nodes listed in row reading order.
\end{Definition}

The real number $\hcoord(m,r,c)$ will be the $x$-coordinate of the corresponding string in the idempotent diagram~$\idem$. Moreover,
\autoref{D:BasesCombinatoricsPosition} looks more complicated than it actually is. The coordinate function simply places
strings in order, following the strategy outlined in
\autoref{R:StrategyMain}.

\begin{Remark}\label{R:BasesCombinatoricsPosition}
The contributions of the different ingredients of $\hcoord(m,r,c)$ are the following:
\begin{enumerate}

\item The $x$-coordinates of the nodes in the component $\lambda^{(m)}$ are centered around $\affine{\kappa}_{m}$, the position of the $m$th (affine) red string.

\item The $\tfrac{m}{\hell}$ differentiates between the components in the sense that the multiple of $\tfrac{1}{\hell}$ uniquely identifies the component of the node.

\item $\co^{m}\big(\re(k)\big)$ ensures that strings with the
same residue are placed within a certain region. The first two
cases in the definition of $\co^{m}$ takes care of the $0$ multisinks, which have special combinatorics, and the shift by $\rho_{m}$ ensures that coordinates of residue $\rho_{m}$ are the rightmost coordinates.

\item $-k\varepsilon$ ensures that nodes move a little bit
to the left as we read along rows.

\end{enumerate}
Note that in type $\dtwo$ the solid $0$-strings
and the solid $2$-strings have roughly the same coordinates. This is important because the solid $0$-strings do not have ghosts, but the solid $1$-strings have two ghosts.
\end{Remark}

\begin{Example}\label{E:BasesBasisExemplifiedPartIII}
We continue with \autoref{E:BasesBasisExemplifiedPartII}. Fix from now on
$\charge=(0)$. We
have the following coordinates of $\nu_{1}$ and $\nu_{2}$:
\begin{gather*}
\begin{array}{*4c}
\hcoord(\nu_{1})=(-\varepsilon,-1-2\varepsilon,-2-3\varepsilon)
,&
\begin{ytableau}
-\varepsilon & \scalebox{0.7}{$\text{-}1\text{-}2\varepsilon$} & \scalebox{0.7}{$\text{-}2\text{-}3\varepsilon$}
\end{ytableau}
&\xrightarrow{\varepsilon\to 0}&
\begin{ytableau}
0 & {-}1 & {-}2
\end{ytableau}
\,,
\\[3mm]
\hcoord(\nu_{2})=(-\varepsilon,-1-2\varepsilon,-3\varepsilon)
,&
\begin{ytableau}
-\varepsilon & \scalebox{0.7}{$\text{-}1\text{-}2\varepsilon$} \\
\none & -3\varepsilon
\end{ytableau}
&\xrightarrow{\varepsilon\to 0}&
\begin{ytableau}
0 & {-}1 \\
\none & 0
\end{ytableau}
\,.
\end{array}
\end{gather*}
We have illustrated the coordinates for the nodes
in the shifted Young diagrams, and also what happens in the
limit $\varepsilon\to 0$. Note hat we have omitted the additional
shift of $\tfrac{m}{\hell}=\tfrac{1}{10}$.
\end{Example}

\begin{Example}\label{E:BasesMain}
Let us continue with \autoref{E:StrategyMain}, where
$n=6$ and $e=2$.
As before, $\blam=(4,2)$. We label the nodes of $\blam$ in
row reading order and reproduce the associated
idempotent diagram:
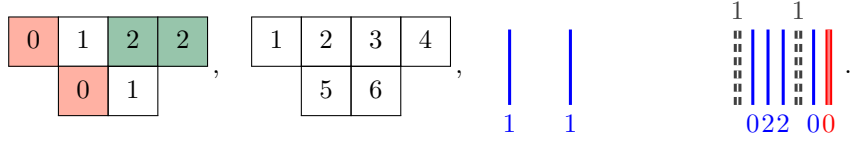
\begin{gather*}
\begin{ytableau}
*(tomato!50) 0 & 1 & *(spinach!50) 2 & *(spinach!50) 2 \\
\none & *(tomato!50) 0 & 1
\end{ytableau}
\,,\quad
\begin{ytableau}
1 & 2 & 3 & 4 \\
\none & 5 & 6
\end{ytableau}
\,,\quad
\DottedIdempotentD[]{2}{0,1,2,2,0,1}
.
\end{gather*}
Let us assume that the red $0$-string is placed at $0$
and this is the only component, {\ie} $\charge=(0)$, $\ell=1$ and $m=1$. We then have $\hell=1+6\cdot(2+1)=19$ and $\frac{m}{\hell}=\frac{1}{19}\approx 0.053$.
Let us fix $\varepsilon=0.001$. The solid strings displayed in
\autoref{E:StrategyMain}.(a) are at positions
\begin{gather*}
\hcoord(\blam)\approx
(\text{-}0.054,\text{-}1.055,\text{-}0.056,\text{-}0.057,\text{-}0.058,\text{-}1.059)
\xrightarrow[\text{left to right}]{\text{ordered}}
(\text{-}1.059,\text{-}1.055,\text{-}0.058,\text{-}0.057,\text{-}0.056,\text{-}0.054).
\end{gather*}
Here, $\hcoord(\blam)$ lists the positions in row reading order while,
on the right-hand side, we list the coordinates for the solid strings with increasing $x$-coordinates from left to right. Note that the illustrated idempotent diagram is scaled to avoid clutter: the distance between the two solid $1$-string is $\approx 0.01$ while the distance from the solid $1$-strings to the red $0$-string is $\approx 1$.

Finally, filling the nodes of $\blam$ with their positions in $\R$ gives:
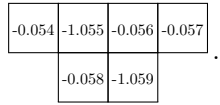
\begin{gather*}
\begin{ytableau}
\scalebox{0.6}{$\text{-}0.054$} & \scalebox{0.6}{$\text{-}1.055$} & \scalebox{0.6}{$\text{-}0.056$} & \scalebox{0.6}{$\text{-}0.057$} \\
\none & \scalebox{0.6}{$\text{-}0.058$} & \scalebox{0.6}{$\text{-}1.059$}
\end{ytableau}
\,.
\end{gather*}
In this sense one can think of $\hcoord(\blam)$ as a reading function on Young diagrams that associates a real number to each node in
a controlled way, which corresponds to the strategy
outlined in \autoref{R:StrategyMain}.
\end{Example}

We write $\hcoordf$, the \emph{maximal finite part},
for the maximum of $\hcoord$ on
$\Parts$, and call coordinates with $\hcoord(m,r,c)>\hcoordf$ \emph{affine}.

\begin{Lemma}\label{L:BasesAffine}
The nodes with affine coordinates are precisely the affine nodes.
\end{Lemma}

\begin{proof}
Easy and omitted.
\end{proof}

We now define the \emph{idempotent diagrams}
$\idem$ associated
to $\blam\in\hPartsall$. Dotted idempotents 
will appear in \autoref{SS:BasesPositioningDot}.

\begin{Definition}\label{D:BasesCombinatoricsIdempotent}
For $\blam\in\hPartsall$ let
$\idem$ be the idempotent diagram given by:
\begin{enumerate}

\item placing $\ell$ red strings with labels given by $\brho$ at coordinates given by $\charge$, and,

\item $n$ solid strings
with labels $\res(m,r,c)$ at coordinates $\hcoord(m,r,c)$, for
$(m,r,c)\in\blam$.

\end{enumerate}
\end{Definition}

The reader is now invited to go back to \autoref{E:StrategyMain} and compare the idempotent diagrams appearing there with \autoref{D:BasesCombinatoricsIdempotent}. See \autoref{E:BasesCombinatoricsDottedIdempotent} below for more examples.

As mentioned already, in diagrams we often draw affine red strings at positions $\kappa_{m}$, for $\ell<m\leq\hell$, even though these strings are not strictly part of the diagram.

\begin{Remark}\label{R:BasesCombinatoricsDottedIdempotent}
In $\ac$ types the crucial
illustrations that we used to show that our bases span the weighted KLRW algebras are \cite[(6A.10), (6A.11), (7A.8), (7A.9) and (7A.10)]{MaTu-klrw-algebras}. These diagrams identify local configurations of nodes in Young diagrams with local configurations of strings in diagrams. The analogs of these results the current setting, together with the Young diagram and the residues of the nodes, are as follows.

Assuming that the middle node is not $0$ or $e$, and we do not have $i-2=0$, $i=0$ or $i+2=0$ in type $\dtwo$, we have:
\begin{gather}\label{Eq:BasesImportantConfigurationOne}
\begin{ytableau}
i & i{+}1 & i{+}2
\end{ytableau}
\leftrightsquigarrow
\begin{tikzpicture}[anchorbase,smallnodes,rounded corners]
\draw[ghost](1,0)node[below]{$\phantom{i}$}--++(0,1)node[above,yshift=-0.05cm]{$i$};
\draw[ghost](1.9,0)node[below]{$\phantom{i}$}--++(0,1)node[above,yshift=-0.05cm]{$i{+}1$};
\draw[ghost](2.9,0)node[below]{$\phantom{i}$}--++(0,1)node[above,yshift=-0.05cm]{$i{+}2$};
\draw[solid](0,0)node[below]{$i$}--++(0,1)node[above,yshift=-1pt]{$\phantom{i}$};
\draw[solid](0.9,0)node[below]{$i{+}1$}--++(0,1)node[above,yshift=-1pt]{$\phantom{i}$};
\draw[solid](1.8,0)node[below]{$i{+}2$}--++(0,1)node[above,yshift=-1pt]{$\phantom{i}$};
\end{tikzpicture}
,\quad
\begin{ytableau}
i & i{-}1 & i{-}2
\end{ytableau}
\leftrightsquigarrow
\begin{tikzpicture}[anchorbase,smallnodes,rounded corners]
\draw[ghost](-1.2,0)node[below]{$\phantom{i}$}--++(0,1)node[above,yshift=-0.05cm]{$i{-}2$};
\draw[ghost](-0.1,0)node[below]{$\phantom{i}$}--++(0,1)node[above,yshift=-0.05cm]{$i{-}1$};
\draw[ghost](1,0)node[below]{$\phantom{i}$}--++(0,1)node[above,yshift=-0.05cm]{$i$};
\draw[solid](-2.2,0)node[below]{$i{-}2$}--++(0,1)node[above,yshift=-1pt]{$\phantom{i}$};
\draw[solid](-1.1,0)node[below]{$i{-}1$}--++(0,1)node[above,yshift=-1pt]{$\phantom{i}$};
\draw[solid](0,0)node[below]{$i$}--++(0,1)node[above,yshift=-1pt]{$\phantom{i}$};
\end{tikzpicture}
.
\end{gather}
The special cases in \autoref{D:BasesContent} for type $\dtwo$ correspond to the local string configurations:
\begin{gather}\label{Eq:BasesImportantConfigurationOneDTwo}
\begin{ytableau}
*(tomato!50) 0 & 1 & 2
\end{ytableau}
\leftrightsquigarrow
\begin{tikzpicture}[anchorbase,smallnodes,rounded corners]
\draw[dghost](-0.1,0)node[below]{$\phantom{i}$}--++(0,1)node[above,yshift=-0.05cm]{$1$};
\draw[ghost](0.8,0)node[below]{$\phantom{i}$}--++(0,1)node[above,yshift=-0.05cm]{$2$};
\draw[solid](-1.1,0)node[below]{$1$}--++(0,1)node[above,yshift=-1pt]{$\phantom{i}$};
\draw[solid](-0.2,0)node[below]{$2$}--++(0,1)node[above,yshift=-1pt]{$\phantom{i}$};
\draw[solid](0,0)node[below]{$0$}--++(0,1)node[above,yshift=-1pt]{$\phantom{i}$};
\end{tikzpicture}
,\quad
\begin{ytableau}
2 & 1 & *(tomato!50) 0
\end{ytableau}
\leftrightsquigarrow
\begin{tikzpicture}[anchorbase,smallnodes,rounded corners]
\draw[dghost](-0.1,0)node[below]{$\phantom{i}$}--++(0,1)node[above,yshift=-0.05cm]{$1$};
\draw[ghost](1,0)node[below]{$\phantom{i}$}--++(0,1)node[above,yshift=-0.05cm]{$2$};
\draw[solid](-1.1,0)node[below]{$1$}--++(0,1)node[above,yshift=-1pt]{$\phantom{i}$};
\draw[solid](-0.2,0)node[below]{$0$}--++(0,1)node[above,yshift=-1pt]{$\phantom{i}$};
\draw[solid](0,0)node[below]{$2$}--++(0,1)node[above,yshift=-1pt]{$\phantom{i}$};
\end{tikzpicture}
.
\end{gather}
That the two cases in \autoref{Eq:BasesImportantConfigurationOne} and \autoref{Eq:BasesImportantConfigurationOneDTwo} look different is an artifact of our conventions for ghost strings. However, this can not be avoided
(meaning that one always gets special behavior) in type
$\dtwo$ as there is no way to orient the quiver from left to right, or right to left.

When the middle residue is $0$ or $e$ we have:
\begin{gather}\label{Eq:BasesImportantConfigurationTwo}
\begin{gathered}
\atwo:
\begin{ytableau}
1 & 0 & 1
\end{ytableau}
\leftrightsquigarrow
\begin{tikzpicture}[anchorbase,smallnodes,rounded corners]
\draw[ghost](-0.1,0)node[below]{$\phantom{i}$}--++(0,1)node[above,yshift=-0.05cm]{$0$};
\draw[ghost](0.8,0)node[below]{$\phantom{i}$}--++(0,1)node[above,yshift=-0.05cm]{$1$};
\draw[ghost](1,0)node[below]{$\phantom{i}$}--++(0,1)node[above,yshift=-0.05cm]{$1$};
\draw[solid](-1.1,0)node[below]{$0$}--++(0,1)node[above,yshift=-1pt]{$\phantom{i}$};
\draw[solid](-0.2,0)node[below]{$1$}--++(0,1)node[above,yshift=-1pt]{$\phantom{i}$};
\draw[solid](0,0)node[below]{$1$}--++(0,1)node[above,yshift=-1pt]{$\phantom{i}$};
\end{tikzpicture}
,\quad
\dtwo:
\begin{ytableau}
1 & *(tomato!50) 0 & *(tomato!50) 0 & 1
\end{ytableau}
\leftrightsquigarrow
\begin{tikzpicture}[anchorbase,smallnodes,rounded corners]
\draw[dghost](0.6,0)node[below]{$\phantom{i}$}--++(0,1)node[above,yshift=-0.05cm]{$1$};
\draw[dghost](1,0)node[below]{$\phantom{i}$}--++(0,1)node[above,yshift=-0.05cm]{$1$};
\draw[solid](-0.4,0)node[below]{$1$}--++(0,1)node[above,yshift=-1pt]{$\phantom{i}$};
\draw[solid](0,0)node[below]{$1$}--++(0,1)node[above,yshift=-1pt]{$\phantom{i}$};
\draw[solid](0.7,0)node[below]{$0$}--++(0,1)node[above,yshift=-1pt]{$\phantom{i}$};
\draw[solid](0.9,0)node[below]{$0$}--++(0,1)node[above,yshift=-1pt]{$\phantom{i}$};
\end{tikzpicture}
,\\
\atwo,\dtwo:
\begin{ytableau}
\scalebox{0.9}{$e{-}1$} & *(spinach!50) e & *(spinach!50) e & \scalebox{0.9}{$e{-}1$}
\end{ytableau}
\leftrightsquigarrow
\begin{tikzpicture}[anchorbase,smallnodes,rounded corners]
\draw[ghost](0.6,0)node[below]{$\phantom{i}$}--++(0,1)node[above,yshift=-0.05cm,xshift=-0.1cm]{$e{-}1$};
\draw[ghost](1,0)node[below]{$\phantom{i}$}--++(0,1)node[above,yshift=-0.05cm,xshift=0.1cm]{$e{-}1$};
\draw[solid](-0.4,0)node[below,xshift=-0.1cm]{$e{-}1$}--++(0,1)node[above,yshift=-1pt]{$\phantom{i}$};
\draw[solid](0,0)node[below,xshift=0.1cm]{$e{-}1$}--++(0,1)node[above,yshift=-1pt]{$\phantom{i}$};
\draw[solid](0.7,0)node[below]{$e$}--++(0,1)node[above,yshift=-1pt]{$\phantom{i}$};
\draw[solid](0.9,0)node[below]{$e$}--++(0,1)node[above,yshift=-1pt]{$\phantom{i}$};
\end{tikzpicture}
.
\end{gathered}
\end{gather}
These should be compared to \autoref{E:StrategyStuck}. Note that in all cases we cannot pull the strings for the middle node further to the right using the Reidemeister II relations, which is the point of these diagrams.

Whenever we have a multisink \autoref{R:StrategyMain}.(b) fails because flanking two solid $i$-strings with two crossings
annihilates a diagram by \autoref{Eq:RecollectionReidemeisterII}. This is not desirable because we want these diagrams to form the middle of a sandwich cellular basis. In \autoref{D:BasesCombinatoricsDottedIdempotent} below, we will place a dot on such strands to avoid this problem. The picture to keep in mind is:
\begin{gather*}
\begin{tikzpicture}[anchorbase,smallnodes,rounded corners]
\draw[solid](0,1)--++(0.5,-0.5)--++(-0.5,-0.5) node[below]{$i$};
\draw[solid](0.5,1)node[above,yshift=-1pt]{$\phantom{i}$}--++(-0.5,-0.5)--++(0.5,-0.5) node[below]{$i$};
\draw (-0.1,0.5)node[left]{no dots};
\end{tikzpicture}
=0
\quad
\text{but}
\quad
\begin{tikzpicture}[anchorbase,smallnodes,rounded corners]
\draw[solid,dot](0,1)--++(0.5,-0.5)--++(-0.5,-0.5) node[below]{$i$};
\draw[solid](0.5,1)node[above,yshift=-1pt]{$\phantom{i}$}--++(-0.5,-0.5)--++(0.5,-0.5) node[below]{$i$};
\draw (-0.1,0.5)node[left]{with a dot};
\end{tikzpicture}
\neq 0
.
\end{gather*}

Finally, assuming that $|i-j|\leq 1$ in these pictures:
\begin{gather}\label{Eq:BasesMainTableaux}
\begin{gathered}
\begin{ytableau}
k{\mp}1 & \dots & j & i
\\
k
\end{ytableau}
\\
\text{or}
\\
\begin{ytableau}
k & \dots & j & i
\\
\none & k
\end{ytableau}
\end{gathered}
\leftrightsquigarrow
\begin{cases}
\begin{tikzpicture}[anchorbase,smallnodes,rounded corners]
\draw[ghost](0.6,0)node[below]{$\phantom{i}$}--++(0,1)node[above,yshift=-0.05cm]{$i$};
\draw[ghost](0.8,0)node[below]{$\phantom{i}$}--++(0,1)node[above,yshift=-0.05cm]{$i$};
\draw[ghost](1,0)node[below]{$\phantom{i}$}--++(0,1)node[above,yshift=-0.05cm]{$i$};
\draw[solid](-0.4,0)node[below]{$i$}--++(0,1)node[above,yshift=-1pt]{$\phantom{i}$};
\draw[solid](-0.2,0)node[below]{$i$}--++(0,1)node[above,yshift=-1pt]{$\phantom{i}$};
\draw[solid](0,0)node[below]{$i$}--++(0,1)node[above,yshift=-1pt]{$\phantom{i}$};
\end{tikzpicture}
&
\begin{array}{@{}l}
\text{if $i=k=j$}\\[0.05cm]
\text{is a multisink},
\end{array}
\\
\begin{tikzpicture}[anchorbase,smallnodes,rounded corners]
\draw[ghost](1.7,0)node[below]{$\phantom{i}$}--++(0,1)node[above,yshift=-0.05cm]{$i$};
\draw[ghost](1.9,0)node[below]{$\phantom{i}$}--++(0,1)node[above,yshift=-0.05cm]{$i$};
\draw[solid](0.7,0)node[below]{$i$}--++(0,1)node[above,yshift=-1pt]{$\phantom{i}$};
\draw[solid](0.9,0)node[below]{$i$}--++(0,1)node[above,yshift=-1pt]{$\phantom{i}$};
\end{tikzpicture}
&\text{if $i=k\neq j$},
\\
\begin{tikzpicture}[anchorbase,smallnodes,rounded corners]
\draw[ghost](1,0)node[below]{$\phantom{i}$}--++(0,1)node[above,yshift=-0.05cm]{$i$};
\draw[ghost](1.9,0)node[below]{$\phantom{i}$}--++(0,1)node[above,yshift=-0.05cm]{$k$};
\draw[solid](0,0)node[below]{$i$}--++(0,1)node[above,yshift=-1pt]{$\phantom{i}$};
\draw[solid](0.9,0)node[below]{$k$}--++(0,1)node[above,yshift=-1pt]{$\phantom{i}$};
\end{tikzpicture}
\quad\text{or}\quad
\begin{tikzpicture}[anchorbase,smallnodes,rounded corners]
\draw[ghost](-0.1,0)node[below]{$\phantom{i}$}--++(0,1)node[above,yshift=-0.05cm]{$k$};
\draw[ghost](1,0)node[below]{$\phantom{i}$}--++(0,1)node[above,yshift=-0.05cm]{$i$};
\draw[solid](-1.1,0)node[below]{$k$}--++(0,1)node[above,yshift=-1pt]{$\phantom{i}$};
\draw[solid](0,0)node[below]{$i$}--++(0,1)node[above,yshift=-1pt]{$\phantom{i}$};
\end{tikzpicture}
&\text{if $|k-i|=1$},
\\
\begin{gathered}
\text{the solid/ghost $k$-string is close}
\\
\text{to a ghost/solid $(k\pm 1)$-string}
\end{gathered}
&\text{otherwise}
,
\end{cases}
\end{gather}
where the $j$-string is only illustrated in the top diagram,
and in the last case the $i$ and $k$-strings
do not need to be close (neither the solids nor the ghosts).
Note that in the shifted Young diagram $k=\rho_{k}$.
\end{Remark}


\subsection{The dotted idempotents}\label{SS:BasesPositioningDot}


We now explain how we place dots on the idempotent diagram $\idem[\blam]$.

\begin{Definition}\label{D:BasesCombinatoricsDottedIdempotent}
Suppose $\blam\in\hPartsall$ and suppose that $(m,r,c),(m,r^{\prime},c^{\prime})\in\blam$ with $k=o_{\blam}(m,r,c)$ and $k+1=o_{\blam}(m,r^{\prime},c^{\prime})$.
If $\rho_{m}$ is not a multisink (so usual partitions), then
\begin{align}\label{Eq:BasesDottedIdempotent}
a_{k}&=
\begin{cases}
1 & \text{if $\res(m,r,c)=\res(m,r+1,1)$},
\\
0 & \text{otherwise.}
\end{cases}
\\
\intertext{If $\rho_{m}$ is a multisink (so shifted partitions), then}
\label{Eq:BasesDottedIdempotentTwo}
a_{k}&=
\begin{cases}
1 & \text{if $\res(m,r,c)=\res(m,r^{\prime},c^{\prime})$, or $r^{\prime}>1$ and $c^{\prime}=1$},
\\
0 & \text{otherwise.}
\end{cases}
\end{align}
The \emph{dotted idempotent} associated to $\blam$
is $\dotidem=y_{\blam}\idem$, where $y_{\blam}=y_{1}^{a_{1}}\dots
y_{n}^{a_{n}}\in R[y_{1},\dots,y_{n}]$.
\end{Definition}

In other words, \autoref{Eq:BasesDottedIdempotentTwo} places a dot
whenever the $k$th and the $(k+1)$th string are close and have the same residue, and takes care of new rows as in \autoref{Eq:BasesMainTableaux}, {\eg}:
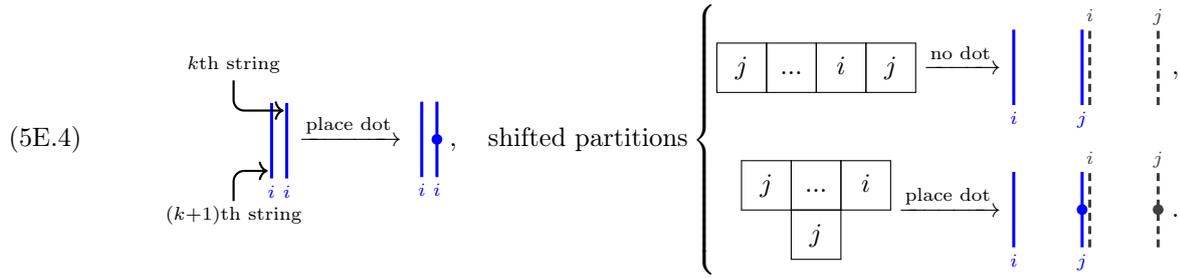
\begin{gather}\label{Eq:BasesDotPlacement}
\begin{tikzpicture}[anchorbase,smallnodes,rounded corners]
\draw[solid](0.7,0)node[below]{$i$}--++(0,1)node[above,yshift=-1pt]{$\phantom{i}$};
\draw[solid](0.9,0)node[below]{$i$}--++(0,1)node[above,yshift=-1pt]{$\phantom{i}$};
\draw[thick,->] (0.2,1.3)node[above]{$k$th string} to (0.2,0.9) to (0.86,0.9);
\draw[thick,->] (0.2,-0.3)node[below]{$(k{+}1)$th string} to (0.2,0.1) to (0.65,0.1);
\end{tikzpicture}
\hspace*{-0.25cm}\xrightarrow{\text{place dot}}
\begin{tikzpicture}[anchorbase,smallnodes,rounded corners]
\draw[solid](0.7,0)node[below]{$i$}--++(0,1)node[above,yshift=-1pt]{$\phantom{i}$};
\draw[solid,dot](0.9,0)node[below]{$i$}--++(0,1)node[above,yshift=-1pt]{$\phantom{i}$};
\end{tikzpicture}
,\quad
\text{shifted partitions}
\left\{
\begin{aligned}
\begin{ytableau}
j & \dots & i & j
\end{ytableau}
\xrightarrow{\text{no dot}}
&\begin{tikzpicture}[anchorbase,smallnodes,rounded corners]
\draw[ghost](1,0)node[below]{$\phantom{i}$}--++(0,1)node[above,yshift=-0.05cm]{$i$};
\draw[ghost](1.9,0)node[below]{$\phantom{i}$}--++(0,1)node[above,yshift=-0.05cm]{$j$};
\draw[solid](0,0)node[below]{$i$}--++(0,1)node[above,yshift=-1pt]{$\phantom{i}$};
\draw[solid](0.9,0)node[below]{$j$}--++(0,1)node[above,yshift=-1pt]{$\phantom{i}$};
\end{tikzpicture},
\\
\begin{ytableau}
j & \dots & i \\
\none & j
\end{ytableau}
\xrightarrow{\text{place dot}}
&\begin{tikzpicture}[anchorbase,smallnodes,rounded corners]
\draw[ghost](1,0)node[below]{$\phantom{i}$}--++(0,1)node[above,yshift=-0.05cm]{$i$};
\draw[ghost,dot](1.9,0)node[below]{$\phantom{i}$}--++(0,1)node[above,yshift=-0.05cm]{$j$};
\draw[solid](0,0)node[below]{$i$}--++(0,1)node[above,yshift=-1pt]{$\phantom{i}$};
\draw[solid,dot](0.9,0)node[below]{$j$}--++(0,1)node[above,yshift=-1pt]{$\phantom{i}$};
\end{tikzpicture}
.
\end{aligned}
\right.
\end{gather}
In the case of usual partitions the dot placement is the same as in type $\cone$ from \cite[Section 7]{MaTu-klrw-algebras}, see also
\autoref{Eq:BasesUsual} below.
Note that $\dotidem$ has zero or one dot on each strand.
For example, the top diagram in
\autoref{Eq:BasesMainTableaux} gets two dots, one on the middle and
one on the rightmost string. 

\begin{Remark}\label{R:BasesCombinatoricsDottedIdempotentDots}
The dot placement in \autoref{D:BasesCombinatoricsDottedIdempotent} is local to the  components of the (shifted) partitions. This is appears to be different to what is done in the KLR world where dots corresponding to later components are added to idempotents; see, for example, \cite[Definition~4.9]{HuMa-klr-basis}. In fact, the red strings implicitly add dots to the idempotent diagrams via \autoref{Eq:RecollectionReidemeisterII}, so these two pictures are compatible but different in flavor. Note since weighted KLRW algebras generalize KLR algebras the above dot placement gives a new way to think about cellular bases constructions in the KLR world.
\end{Remark}

\begin{Example}\label{E:BasesCombinatoricsDottedIdempotent}
We list a few examples of
the dotted idempotents $\dotidem$ for $\ell=1$ and $\charge=(0)$.
\begin{enumerate}

\item Consider $\blam=(9)$ for $e=3$ and $\brho=(3)$. We get:
\begin{gather*}
\atwo[3]:
\begin{ytableau}
*(spinach!50) 3 & 2 & 1 & 0 & 1 & 2 & *(spinach!50) 3 & *(spinach!50) 3 & 2 \\
\end{ytableau}
,\quad
\dtwo[3]:
\begin{ytableau}
*(spinach!50) 3 & 2 & 1 & *(tomato!50) 0 & *(tomato!50) 0 & 1 & 2 & *(spinach!50) 3 & *(spinach!50) 3 \\
\end{ytableau}
\,.
\end{gather*}
\begin{gather*}
\atwo[3]:
\dotidem=
\DottedIdempotentAA[7]{3}{3,2,1,0,1,2,3,3,2}
,
\\
\dtwo[3]:
\dotidem=
\DottedIdempotentD[4,8]{3}{3,2,1,0,0,1,2,3,3}
.
\end{gather*}
In contrast, if we instead have $\brho=(0)$, then we get
\begin{gather*}
\atwo[3]:
\begin{ytableau}
0 & 1 & 2 & *(spinach!50) 3 & *(spinach!50) 3 & 2 & 1 & 0 & 1 \\
\end{ytableau}
,\quad
\dtwo[3]:
\begin{ytableau}
*(tomato!50) 0 & 1 & 2 & *(spinach!50) 3 & *(spinach!50) 3 & 2 & 1 & *(tomato!50) 0 & *(tomato!50) 0 \\
\end{ytableau}
\,.
\end{gather*}
\begin{gather*}
\atwo[3]:
\dotidem=
\DottedIdempotentAA[2,4,7,9]{3}{0,1,2,3,3,2,1,0,1}
,
\\
\dtwo[3]:
\dotidem=
\DottedIdempotentD[4,8]{3}{0,1,2,3,3,2,1,0,0}
.
\end{gather*}
Note the different dot placement: for type $\atwo[3]$ we use \autoref{Eq:BasesDottedIdempotent}, while for $\dtwo[3]$ we use \autoref{Eq:BasesDottedIdempotentTwo}.

\item Next, we illustrate the dotted idempotents for
$\blam=(8,1)$, $e=2$ and $\brho=(2)$:
\begin{gather*}
\atwo[2]:
\begin{ytableau}
*(spinach!50) 2 & 1 & 0 & 1 & *(spinach!50) 2 & *(spinach!50) 2 & 1 & 0 \\
\none & *(spinach!50) 2
\end{ytableau}
,\quad
\dtwo[2]:
\begin{ytableau}
*(spinach!50) 2 & 1 & *(tomato!50) 0 & *(tomato!50) 0 & 1 & *(spinach!50) 2 & *(spinach!50) 2 & 1 \\
\none & *(spinach!50) 2
\end{ytableau}
,
\\
\atwo[2]:
\dotidem=
\DottedIdempotentAA[5]{2}{2,1,0,1,2,2,1,0,2}
,
\\
\dtwo[2]:
\dotidem=
\DottedIdempotentD[3,6]{2}{2,1,0,0,1,2,2,1,2}
.
\end{gather*}

\item Finally, consider $\blam=(7,2)$ for $e=2$ and $\brho=(0)$. Then:
\begin{gather*}
\atwo[2]:
\begin{ytableau}
0 & 1 & *(spinach!50) 2 & *(spinach!50) 2 & 1 & 0 & 1 \\
1 & 0
\end{ytableau}
,\quad
\dtwo[2]:
\begin{ytableau}
*(tomato!50) 0 & 1 & *(spinach!50) 2 & *(spinach!50) 2 & 1 & *(tomato!50) 0 & *(tomato!50) 0 \\
\none & *(tomato!50) 0 & 1
\end{ytableau}
,
\\
\atwo[2]:
\dotidem=
\DottedIdempotentAA[2,3,5,7]{2}{0,1,2,2,1,0,1,1,0}
,
\\
\dtwo[2]:
\dotidem=
\DottedIdempotentD[3,6,7]{2}{0,1,2,2,1,0,0,0,1}
.
\end{gather*}
Again, note the different dot placements for types $\atwo[2]$ and $\dtwo[2]$.

\end{enumerate}
All of these pictures were created using TikZ macros defined in the preamble.
\end{Example}

\begin{Example}\label{E:BasesBasisExemplifiedPartIV}
We continue with \autoref{E:BasesBasisExemplifiedPartIII}.
Recall that $n=3$, $e=2$, $\ell=1$ and $\charge=(0)$, and fix $\rho=(0)$, so $\hell=10$. Let $\blam_{1}=(3|\emptyset|\dots|\emptyset)$, $\blam_{2}=(2,1|\emptyset|\dots|\emptyset)$ and $\blam_{3}=(1^{3}|\emptyset|\dots|\emptyset)$, and also $\bmu_{1}=(2|1|\emptyset|\dots|\emptyset)$ and $\bmu_{2}=(1|1|1|\emptyset|\dots|\emptyset)$
These five diagrams have the following
dotted idempotent diagrams for type
$\atwo[2]$:
\begin{gather*}
\dotidem[\blam_{1}]
=\DottedIdempotentAA[2]{2}{0,1,2}
,\quad
\dotidem[\blam_{2}]
=\DottedIdempotentAA[2]{2}{0,1,1}
,\\
\dotidem[\blam_{3}]
=\DottedIdempotentAA[]{2}{0,1,2}
,\\
\dotidem[\bmu_{1}]
=\DottedIdempotentAA[2]{2}{0,1}\DottedIdempotentAA*(2.6)[]{2}{0}
,\\
\dotidem[\bmu_{2}]
=\DottedIdempotentAA[]{2}{0}\DottedIdempotentAA*(2.6)[]{2}{0}\DottedIdempotentAA*(2.6)[]{2}{0}
.
\end{gather*}
Here, we have drawn affine red strings only for the final two diagrams. The dotted idempotents
for type $\dtwo$ look similar, but since there is no ghost $0$-string
the $1$ and $2$ strings move to the left. For example
\begin{gather*}
\dotidem[\blam_{1}]
=\DottedIdempotentD{2}{0,1,2}
\end{gather*}
is the diagram for $\dotidem[\blam_{1}]$ for type $\dtwo[2]$.
\end{Example}

\begin{Lemma}\label{L:BasesCombinatoricsIdempotent}
We have $\dotidem[\blam]=\dotidem[\bmu]$ if and only if
$\blam=\bmu$.
\end{Lemma}

\begin{proof}
One direction is immediate, so let us assume that $\blam\neq\bmu$.
If $\res(\blam)\neq\res(\bmu)$, then $\dotidem[\blam]\neq\dotidem[\bmu]$
follows using the faithful polynomial module from \autoref{R:RecollectionWABasis}.
Moreover, a different dot placement on the same idempotent diagram $\idem$ can also be distinguished by \autoref{R:RecollectionWABasis}, so it remains
to show that $\dotidem[\blam]$ and $\dotidem[\bmu]$ have  different dot placements
if $\res(\blam)=\res(\bmu)$. To see this assume that the $k$th node, in row reading order, is
the first node that is different for $\blam$ and $\bmu$. By the minimality of $k$, we can assume that the $k$th node of $\bmu$ is
in a new row when compared to the $k$th node of $\blam$.
There are three cases to check now, depending on the residue $j$ of the
$k$th string and the residue $i$ of the $(k-1)$th string. When $i=j$ is a multisink the local diagrams for $\blam$ and $\bmu$ are:
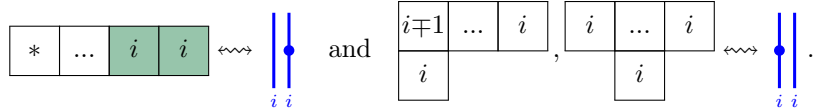
\begin{gather*}
\begin{ytableau}
\ast & \dots & *(spinach!50) i & *(spinach!50) i
\end{ytableau}
\leftrightsquigarrow
\begin{tikzpicture}[anchorbase,smallnodes,rounded corners]
\draw[solid](-0.2,0)node[below]{$i$}--++(0,1)node[above,yshift=-1pt]{$\phantom{i}$};
\draw[solid,dot](0,0)node[below]{$i$}--++(0,1)node[above,yshift=-1pt]{$\phantom{i}$};
\end{tikzpicture}
\quad\text{and}\quad
\begin{ytableau}
i{\mp}1 & \dots & i \\
i
\end{ytableau}
\,,
\begin{ytableau}
i & \dots & i \\
\none & i
\end{ytableau}
\leftrightsquigarrow
\begin{tikzpicture}[anchorbase,smallnodes,rounded corners]
\draw[solid,dot](-0.2,0)node[below]{$i$}--++(0,1)node[above,yshift=-1pt]{$\phantom{i}$};
\draw[solid](0,0)node[below]{$i$}--++(0,1)node[above,yshift=-1pt]{$\phantom{i}$};
\end{tikzpicture}
.
\end{gather*}
As required, these diagrams have different dot placements. The case $i\rightsquigarrow j$
for shifted partitions is illustrated in \autoref{Eq:BasesDotPlacement}, while
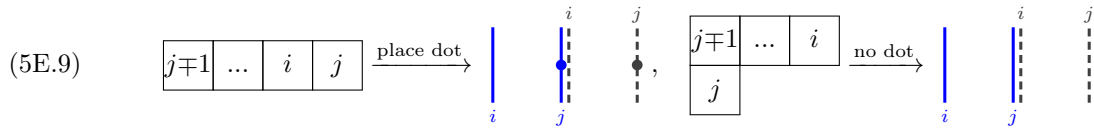
\begin{gather}\label{Eq:BasesUsual}
\begin{ytableau}
j{\mp}1 & \dots & i & j
\end{ytableau}
\xrightarrow{\text{place dot}}
\begin{tikzpicture}[anchorbase,smallnodes,rounded corners]
\draw[ghost](1,0)node[below]{$\phantom{i}$}--++(0,1)node[above,yshift=-0.05cm]{$i$};
\draw[ghost,dot](1.9,0)node[below]{$\phantom{i}$}--++(0,1)node[above,yshift=-0.05cm]{$j$};
\draw[solid](0,0)node[below]{$i$}--++(0,1)node[above,yshift=-1pt]{$\phantom{i}$};
\draw[solid,dot](0.9,0)node[below]{$j$}--++(0,1)node[above,yshift=-1pt]{$\phantom{i}$};
\end{tikzpicture}
,\quad
\begin{ytableau}
j{\mp}1 & \dots & i \\
j
\end{ytableau}
\xrightarrow{\text{no dot}}
\begin{tikzpicture}[anchorbase,smallnodes,rounded corners]
\draw[ghost](1,0)node[below]{$\phantom{i}$}--++(0,1)node[above,yshift=-0.05cm]{$i$};
\draw[ghost](1.9,0)node[below]{$\phantom{i}$}--++(0,1)node[above,yshift=-0.05cm]{$j$};
\draw[solid](0,0)node[below]{$i$}--++(0,1)node[above,yshift=-1pt]{$\phantom{i}$};
\draw[solid](0.9,0)node[below]{$j$}--++(0,1)node[above,yshift=-1pt]{$\phantom{i}$};
\end{tikzpicture}
\end{gather}
is the case for usual partitions.
The case $i\leftsquigarrow j$ is similar to \autoref{Eq:BasesDotPlacement}
and \autoref{Eq:BasesUsual}. Hence, in all cases the dot placements in $\dotidem[\blam]$ and $\dotidem[\bmu]$ are different, so the lemma is proved.
\end{proof}


\subsection{The sandwiched part}\label{SS:BasesSandwichPart}


The following configuration of close strings is needed in \autoref{D:BasesSandwichpart}:
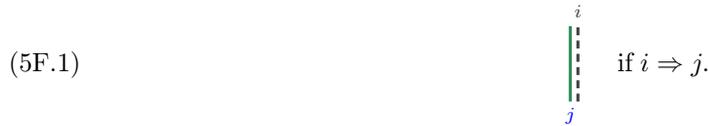
\begin{gather}\label{Eq:BasesSandwichDots}
\begin{tikzpicture}[anchorbase,smallnodes,rounded corners]
\draw[ghost](1,0)node[below]{$\phantom{i}$}--++(0,1)node[above,yshift=-0.05cm]{$i$};
\draw[solid,spinach](0.9,0)node[below]{$j$}--++(0,1)node[above,yshift=-1pt]{$\phantom{i}$};
\end{tikzpicture}
\quad\text{if}\;
i\Rightarrow j.
\end{gather}
(This diagram does not arise if $i$ has no ghost.)
In this diagram we have highlighted the $j$-string, which can get a dot in the sandwiched part.

\begin{Definition}\label{D:BasesSandwichpart}
Define the set of \emph{affine dots} as
\begin{gather*}
\Affch=\set[\big]{\ba=(a_{1},\dots,a_{n})\in\N^{n}|a_{k}=0\text{ whenever }
\hcoord(\blam)_{k}\leq\hcoordf}.
\end{gather*}
For $1\leq k\leq n$, define $c_{k}(\blam)=1$ if the $k$th string is close
as in \autoref{Eq:BasesSandwichDots}, or its partner relation, and this string does not
already have a dot in $\dotidem[\bmu]$ for some $\bmu\in\hPartsall$
with $\idem[\blam]=\idem[\bmu]$, and otherwise set $c_{k}(\blam)=0$.
Define the set of \emph{finite dots} to be
\begin{gather*}
\Finch=\set[\big]{\bfi=(f_{1},\dots,f_{n})\in\set{0,1}^{n}|0\leq f_{k}\leq c_{k}(\blam)}
.
\end{gather*}
\end{Definition}
The set of \emph{sandwiched dots} is
$\Sandch=\Affch\cup\Finch$. Note that $a_{k}\neq 0$ only
happens for affine coordinates $\hcoord(\blam)_{k}>\hcoordf$.

\begin{Example}\label{E:BasesBasisExemplifiedPartV}
We continue with \autoref{E:BasesBasisExemplifiedPartIV}. The condition $\hcoord(\blam)_{k}\leq\hcoordf$ implies that
$\Affch[\blam_{1}]=\Affch[\blam_{2}]=\Affch[\blam_{3}]=\set{(0,0,0)}$. Moreover,
all strings associated to the second and third component have affine coordinates
so $\Affch[\bmu_{1}]=\set{(0,0,n)|n\in\N}$ and
$\Affch[\bmu_{2}]=\set{(0,m,n)|m,n\in\N}$.

For the finite dots, in type $\atwo[2]$ we have $\Finch[\blam_{1}]=\Finch[\blam_{3}]=\set{(0,0,0),(0,0,1)}$
and $\Finch[\blam_{2}]=\Finch[\bmu_{1}]=\Finch[\bmu_{2}]=\set{(0,0,0)}$.
For type $\dtwo[3]$ the finite dots are same except that
$\Finch[\blam_{1}]=\set{(0,0,0),(0,0,1)}$.
\end{Example}


\subsection{Permutation diagrams}\label{SS:BasesPermuations}


We need the same tableaux as in $\ac$ types. That is,
let $\hcoorda(\placeholder)$ be the type~$A$ positioning
function, see \cite[Definition 6A.4]{MaTu-klrw-algebras}, which,
in the notation of this paper is
\begin{gather*}
\hcoorda(m,r,c)=
\affine{\kappa}_{m}
-\frac{m}{\hell}
+\re(k)
+k\varepsilon,
\end{gather*}
where $k=o_{\blam}(m,r,c)$ is the number given by the row reading order.

\begin{Definition}\label{D:BasesAllTheCombinatorics}
Let $\blam,\bmu\in\hPartsall$. A \emph{$\blam$-tableau of type $\bmu$} is a bijection $\bT\map\blam\hcoorda(\bmu)$.
Such a tableau is \emph{semistandard} if:
\begin{enumerate}

\item $\bT(m,1,1)\leq\kappa_{m}$ for $1\leq m\leq\ell$.

\item $\bT(m,r,c)+1>\bT(m,r-1,c)$ for all $(m,r,c),(m,r-1,c)\in\blam$.

\item $\bT(m,r,c)>\bT(m,r,c-1)+1$ for all?
$(m,r,c),(m,r,c-1)\in\blam$.

\end{enumerate}
Let $\hSStd(\blam,\bmu)$ be the set of semistandard $\blam$-tableaux of
type $\bmu$ and set $\hSStd(\blam)=\bigcup_{\bmu}\hSStd(\blam,\bmu)$.
\end{Definition}

\begin{Definition}\label{D:BasesAllTheCombinatoricsII}
For $\bT\in\hSStd(\blam,\bmu)$ define the permutation $w_{\bT}\in\Sym$
by requiring that
\begin{gather*}
x^{\bmu}_{w_{\bT}(k)}=\bT(m,r,c)\text{ whenever }
x^{\blam}_{k}=\hcoorda(m,r,c),
\quad\text{for $1\leq k\leq n$ and $(m,r,c)\in\blam$}.
\end{gather*}
This defines the \emph{permutation diagram} $D_{\bT}=D(w_{\bT})$
from $\hcoord(\bmu)$ to $\hcoord(\blam)$, as in \autoref{N:RecollectionNotationNeeded2}.
\end{Definition}

In other words, a semistandard $\blam$-tableau
$\bT$ of type $\bmu$ is a filling of the nodes of $\blam$ with the type $A$ coordinates of $\bmu$,
together with an anchor condition,
such that the fillings decrease along rows and columns
with an offset of $1$. The associated permutation has top points
defined by $\hcoord(\blam)$, bottom points by $\hcoord(\bmu)$ and
permutes them according to the entries of $\bT$.

\begin{Example}\label{E:BasesTableaux}
For $e=2$ let $n=6$, $\ell=1$, $\charge=(0)$ and $\varepsilon=0.05$.
Fix $\blam=(5,1)$ and $\bmu=(3,2,1)$ for $\brho=(2)$. Then
$\hcoorda(\blam)=(-0.15,-0.05,0.9,1.85,2.8,3.75)$ and
$\hcoorda(\bmu)=(-0.25,-0.15,-0.05,0.8,0.9,1.85)$. Two semistandard
$\blam$-tableaux, one of type $\blam$ and one of type $\bmu$, are:
\begin{gather*}
\bS=
\begin{ytableau}
\scriptstyle \text{-}0.05 & \scriptstyle 0.9 & \scriptstyle 1.85 & \scriptstyle 2.8 & \scriptstyle 3.75 \\
\none & \scriptstyle \text{-}0.15
\end{ytableau}
\;,\quad
D_{\bS}=\idem[\blam],
\quad
\bT=
\begin{ytableau}
\scriptstyle \text{-}0.15 & \scriptstyle 0.8 & \scriptstyle \text{-}0.05 & \scriptstyle 0.9 & \scriptstyle 1.85 \\
\none & \scriptstyle \text{-}0.25
\end{ytableau}
\leftrightsquigarrow
\begin{tikzpicture}[anchorbase,smallnodes,rounded corners]
\draw[solid](0,0) to (0,1);
\draw[solid](0.5,0) to (0.5,1);
\draw[solid](1,0) to (1.5,1);
\draw[solid](1.5,0) to (1,1);
\draw[solid](2,0) to (2,1);
\draw[solid](2.5,0) to (2.5,1);
\node at (3,-0.2) {$\hcoorda(\bmu)$};
\draw[very thick,->] (3,0) to (3,0.5)node[right]{read} to (3,1);
\node at (3,1.2) {$\hcoorda(\blam)$};
\end{tikzpicture}
.
\end{gather*}
The tableau $\bS$ is the \emph{canonical $\blam$-tableau} of type $\blam$, where all nodes are filled
with their coordinates. Its associated permutation diagram is the identity.
The permutation diagram $D_{\bT}$ is build from the
permutation illustrated above, which connects $\hcoorda(\bmu)$ to
$\hcoorda(\blam)$, using $\hcoord(\bmu)$ at the bottom
and $\hcoord(\blam)$ at the top.
\end{Example}

\begin{Example}\label{E:BasesTableauxHard}
Related to our main example from \autoref{E:StrategyMain}.
We use the numbers from \autoref{E:BasesMain}, {\eg} $\frac{m}{\hell}=\frac{1}{19}$ and $\varepsilon=0.001$.
In this case, the type $A$ positioning function gives
\begin{gather*}
\begin{ytableau}
\scalebox{0.6}{$\text{-}0.052$} & \scalebox{0.6}{$0.949$} & \scalebox{0.6}{$1.950$} & \scalebox{0.6}{$1.951$} \\
\none & \scalebox{0.6}{$\text{-}0.048$} & \scalebox{0.6}{$0.953$}
\end{ytableau}
\,.
\end{gather*}
In particular, going from this multitableau to itself, we can only permute the two solid $2$-strings without violating the condition of being standard. This gives precisely the four possible
permutation diagrams from \autoref{E:StrategyMain}.
\end{Example}


\subsection{Basis diagrams}\label{SS:BasesBasisDiagrams}


Before we can define our cellular bases we first need to fix a set $X$ of endpoints for our diagrams.

An \emph{addable node} of $\blam$ is a node $(m,r,c)\notin\blam$ such that $\blam\cup\set{(m,r,c)}$ is the diagram of a $\affine{\brho}$-partition.

\begin{Definition}\label{D:BasesEndPoints}
Let $\hParts[0]=\set{{(\emptyset|\dots|\emptyset)}}$.
For $n\geq 1$, let $\hParts$ be the subset of $\hPartsall$ defined by the condition that
$\blam\in\hParts$ only if $\blam=\bmu\cup\alpha$,
where $\bmu\in\hParts[n-1]$ and $\alpha$ is an addable $i$-node of
$\bmu$ such that:
\begin{gather*}
\text{whenever $\beta$ is an addable $i$-node of $\bmu$
with $\hcoord(\beta)<\hcoord(\alpha)$, then
$\hcoord(\beta)\leq\hcoordf$.}
\end{gather*}
Finally, let $X=\bigcup_{\blam}\hcoord(\blam)$ be the set of all
coordinates, for $\blam\in\hParts$.
\end{Definition}

\begin{Remark}\label{R:BasesEndPoints}
In general, only a handful of the elements of $\hPartsall$ belong to $\hParts$ as the rule in
\autoref{D:BasesEndPoints} only allows the addition of affine nodes that are not too far to the right.
\end{Remark}

We are ready for our main definition:

\begin{Definition}\label{R:BasesTheBasis}
For $\blam\in\hParts$ set $\sandwich{\blam}{\ba}{\bfi}=\daffine\zeetwo\dotidem$ and define
\begin{gather}\label{Eq:BasesTheBasisElements}
D_{\bS\bT}^{\ba,\bfi}=(D_{\bS})^{\star}\sandwich{\blam}{\ba}{\bfi}D_{\bT}
\end{gather}
for $\ba=(a_{1},\dots,a_{n})\in\Affch$, $\bfi=(f_{1},\dots,f_{n})\in\Finch$,
$\bS\in\hSStd(\blam,\bnu)$, $\bT\in\hSStd(\blam,\bmu)$.
\end{Definition}

\begin{Remark}\label{R:BasesSandwichCellularAlgebraNotation}
Note that we use $D^{\ba,\bfi}_{\bS\bT}$ in \autoref{Eq:BasesTheBasisElements}
to distinguish it from the abstract definition 
of \autoref{D:SandwichCellularAlgebra}. Of course, we will show that
these elements are examples of the $C^{b}_{ST}$ in this definition.
\end{Remark}

We call $\sandwich{\blam}{\ba}{\bfi}=\daffine\zeetwo\dotidem$
the \emph{sandwiched part},
$\daffine$ the \emph{affine part},
and $\zeetwo$ the \emph{finite part} of
$D_{\bS\bT}^{\ba,\bfi}$.
The following are the bases
that we consider.

\begin{Definition}\label{D:BasesTheBasis}
Let $\fsandwich{\blam}{\bfi}=\sandwich{\blam}{(0,\dots,0)}{\bfi}$ and  $D_{\bS\bT}^{\bfi}=D_{\bS\bT}^{(0,\dots,0),\bfi}$. We define
\begin{gather}\label{Eq:BasesTheBasis}
\BX=\set{D_{\bS\bT}^{\ba,\bfi}|\blam\in\hParts,\bS,\bT\in\SStd(\blam),\ba\in\Affch,\bfi\in\Finch}.
\end{gather}
\begin{gather}\label{Eq:BasesTheBasisFinite}
\BXc=\set{D_{\bS\bT}^{\bfi}|\blam\in\Parts,\bS,\bT\in\SStd(\blam),\bfi\in\Finch}.
\end{gather}
\end{Definition}


\subsection{Homogeneous (affine) sandwich cellular bases}\label{SS:BasesConstruction}


The sandwich cellular partial order, which we define below, is the same as in \cite[Definition 7C.1]{MaTu-klrw-algebras}. In essence, we measure how far solid strings are to the right, with diagrams becoming bigger in the partial order as we move solid strings to the right.

\begin{Definition}\label{D:BasesDominance}
Let $\blam,\bmu\in\hParts$. Then $\blam$ \emph{dominates} $\bmu$, written $\blam\gedom\bmu$, if there exists a bijection $d\colon\blam\to\bmu$ such that
$\hcoord(\alpha)\geq\hcoord(d(\alpha))$ and the solid string in $\idem$ at position $\hcoord(\alpha)$ has at least as many dots as the solid string in $\idem[\mu]$ at position $\hcoord(d(\alpha))$, for all $\alpha\in\blam$. Write $\blam\gdom\bmu$ if $\blam\gedom\bmu$ and $\blam\neq\bmu$.
\end{Definition}

We are now ready to define (involutive) bases for $\WA[n](X)$ and
$\WAc[n](X)$. Recall that we are working with
$\bad$ types. We also use the $Q$-polynomials
from \autoref{Eq:RecollectionQPoly}.

The cell datum $\mathscr{C}=(\hParts,T,\sand[],\sandbasis[],\BX[], \deg,(\placeholder)^{\ast})$ that we use is:
\begin{itemize}

\item The middle set is $\hParts=(\hParts,\gedom)$,

\item $T=\bigcup_{\blam\in\hParts}\SStd$ is the bottom/top set,

\item the bases of the sandwiched algebras are $\sandbasis[\blam]=\set[\big]{\sandwich{\blam}{\ba}{\bfi}|\ba\in\Affch,\bfi\in\Finch}$ for $\blam\in\hParts$, and the sandwiched algebras $\sand[\blam]$ are the subalgebras of the weighted KLRW algebra generated by these bases,

\item we take $\BX[]=\BX$ from \autoref{Eq:BasesTheBasis}, viewed as a map, as our basis,

\item the degree is $\bS\mapsto\deg D_{\bS}$,

\item the antiinvolution is the diagrammatic antiinvolution $(\placeholder)^{\star}$.

\end{itemize}

The proof of
the following theorem is postponed until \autoref{S:Proofs} below.

\begin{Theorem}\label{T:BasesMain}
The datum $\mathscr{C}$ is a graded affine sandwich cell datum for $\WA[n](X)$.
In particular, \autoref{Eq:BasesTheBasis} is a homogeneous affine sandwich
cellular basis for $\WA[n](X)$.
\end{Theorem}

If we replace $\sandbasis[\blam]=\set{\sandwich{\blam}{\ba}{\bfi}|\ba\in\Affch,\bfi\in\Finch}$ with
$\sandbasis[\blam]^{c}=\set{\fsandwich{\blam}{\bfi}|\bfi\in\Finch}$,
and $\BX$ with $\BXc$, then
\autoref{T:BasesMain} and the definitions directly imply:

\begin{Corollary}\label{C:BasesMain}
The datum $\mathscr{C}=(\hParts,T,\sand[]^{c},\sandbasis[]^{c},\BXc,\deg,(\placeholder)^{\ast})$ is a graded sandwich cell datum for $\WAc[n](X)$.
In particular, \autoref{Eq:BasesTheBasisFinite} is a homogeneous sandwich
cellular basis for $\WAc[n](X)$.
\end{Corollary}

Define a semistandard tableaux $S\in\SStd(\blam)$
to be \emph{standard} if it is of type $\bom=(n|0|\dots|0)$.
Let $\Std(\blam)$ be the set of standard $\blam$-tableaux.
Let $\TA[n]=\idem[\bom]\WA[n](X)\idem[\bom]$ and $\TAc=\idem[\bom]\WAc[n](X)\idem[\bom]$ be the associated KLR and cyclotomic KLR algebra, respectively. This terminology is justified at
the end of \autoref{SS:RecollectionDiagrams}.
We use~$E$ instead of $D$ to refer to the basis elements of the idempotent truncations. We will not
highlight the cell datum below.

\begin{Proposition}\label{P:BasesKLRCellular}
The set $\EX=\set[\big]{E^{\ba,\bfi}_{\bs\bt}|\blam\in\hParts,\bS,\bT\in\Std(\blam),\ba\in\Affch,\bfi\in\Finch}$
is a homogeneous affine sandwich cellular basis of
$\TA[n]$.
\end{Proposition}

\begin{proof}
Apply \cite[Proposition 3F.1 and Example 6A.11]{MaTu-klrw-algebras}.
\end{proof}

As before we obtain:

\begin{Corollary}\label{C:BasesKLRCellular}
The set
$\EXc=\set[\big]{E^{\bfi}_{\bs\bt}|\blam\in\hParts,\bS,\bT\in\Std(\blam),\bfi\in\Finch}$
is a homogeneous sandwich cellular basis of $\TAc[n]$.
\end{Corollary}

Let us now discuss the upshot of \autoref{T:BasesMain} and \autoref{C:BasesMain} for simple modules. To this end,
let $a(\blam)$ and $f(\blam)$ be the number of possible nonzero
positions of $\Affch$ and $\Finch$, respectively.

\begin{Lemma}\label{L:BasesAffineAndFiniteParts}
For all $\blam\in\hParts$ we have $\sand[\blam]\cong R[X_{1},\dots,X_{a(\blam)}]\otimes R[Y_{1},\dots,Y_{f(\blam)}]/(Y_{1}^{2},\dots,Y_{f(\blam)}^{2})$
and $\sand[\blam]^{c}\cong R[Y_{1},\dots,Y_{f(\blam)}]/(Y_{1}^{2},\dots,Y_{f(\blam)}^{2})$.
\end{Lemma}

\begin{proof}
The first claim, regarding the $X_{i}$, follows by using \autoref{P:RecollectionWABasis}.
For the second claim, regarding the $Y_{i}$, we pull strings and
jump dots to the right and use \autoref{C:BasesMain}. That is, if one of the strings, which corresponds to one of the possible nonzero
positions of $\Finch$, carries two dots we can use
\autoref{Eq:StrategyReidemeisterIIAndIIITypesBAD}
and the claim follows.
\end{proof}

\begin{Proposition}\label{P:BasesDimensions}
The algebra $\WAc[n](X)$ is free of rank $\sum_{\blam\in\Parts}2^{f(\blam)}
(\#\SStd)^{2}$, and $\TAc[n](X)$ is free of rank $\sum_{\blam\in\Parts}2^{f(\blam)}
(\#\Std)^{2}$.
\end{Proposition}

\begin{proof}
Directly from \autoref{L:BasesAffineAndFiniteParts} and the corollaries above.
\end{proof}

By the above, it is also easy to write down the graded dimensions
of these algebras.

\begin{Remark}\label{R:BasesDimensions}
\autoref{P:BasesDimensions} generalizes the dimension
formulas from \cite[Corollary 3.5]{ArPa-quiver-hecke-a2} and \cite[Corollary 3.3]{ArPa-quiver-hecke-d}.
\end{Remark}

\begin{Example}\label{E:BasesCompare}
It is illustrative to compare the combinatorics for
types $\cone[2]$, $\atwo[2]$
and $\dtwo[2]$. (See \cite[Section 7]{MaTu-klrw-algebras} for
the relevant constructions in type $\cone$.) We will ignore the affine part
in this example. Fix $\ell=1$, $\charge=(0)$
and $\brho=(0)$ and take $\beta=(0,1,2)\in I^{3}$. We are
first looking for all $1$-partitions of $3$ that have $\beta$ as their
\emph{residue sequence}, that is, $\blam\in\Parts[4]$ whose nodes have residues
$\beta$ in row reading order. We get the following $1$-partitions:
\begin{gather*}
\cone[2],\atwo[2],\dtwo[2]\colon
\blam_{1}=
\begin{ytableau}
0 & 1 & 2
\end{ytableau}
,\quad
\cone[2],\atwo[2]\colon
\blam_{3}=
\begin{ytableau}
0 \\ 1 \\ 2
\end{ytableau}
\,.
\end{gather*}
The dotted idempotents for types $\atwo[2]$
and $\dtwo[2]$ are displayed in
\autoref{E:BasesBasisExemplifiedPartIV},
and the ones for type $\cone[2]$ have the same form as the ones for
type $\atwo[2]$. The sets of finite dots for
types $\atwo[2]$
and $\dtwo[2]$ are listed in \autoref{E:BasesBasisExemplifiedPartV},
while the set of finite dots is always trivial in type $\cone[2]$. Thus, we get
by using \autoref{T:BasesMain} that
\begin{gather*}
\renewcommand*{\arraystretch}{1.25}
\begin{tabular}{c||c|c|c}
$\beta=(0,1,2)$ & $\cone[2]$ & $\atwo[2]$ & $\dtwo[2]$ \\
\hline
\# tableaux & $2$ & $2$ & $1$ \\
\hline
graded dim & $1+q^{2}$ & $(1+q^{2})(1+q^{4})$ & $1+q^{2}$ \\
\end{tabular}
,\quad
\begin{tabular}{c||c|c|c}
$\beta^{\prime}=(2,1,0)$ & $\cone[2]$ & $\atwo[2]$ & $\dtwo[2]$ \\
\hline
\# tableaux & $2$ & $1$ & $1$ \\
\hline
graded dim & $1+q^{2}$ & $1$ & $1+q^{2}$ \\
\end{tabular}
.
\end{gather*}
Here we have listed the graded dimension (using the usual $q$-notation
indicating the degree) of
the idempotent truncation of $\WAc(X)$ determined by $\beta$.
We also listed the relevant numbers for $\beta^{\prime}=(2,1,0)$ where $\brho=(2)$.
These numbers match \cite[Theorem 1.1]{HuSh-monomial-klr-basis},
which is expected as this case is the
cyclotomic KLR algebra of the respective types.
The reader is also invited to compare the graded dimension formula from \autoref{E:StrategyMain} to \cite[Theorem 1.1]{HuSh-monomial-klr-basis}; again we get the same result as expected.
\end{Example}

\begin{Proposition}\label{P:BasesSimplesC}
Suppose that $R$ is a field, and let $\hPartsgood$
or $\Partsgood$ be the
sets of apexes.
\begin{enumerate}

\item For a fixed apex $\blam\in\hPartsgood$ there exist a 1:1-correspondence
between simple $\WA[n](X)$-modules with apex $\blam$ and $R^{a(\blam)}$.
Moreover, up to isomorphism, there exists
exactly one graded simple $\WA[n](X)$-module of that apex.

\item For a fixed apex $\blam\in\Partsgood$ there exists
exactly one simple, and one graded simple, $\WAc[n](X)$-module of that apex up to isomorphism.

\end{enumerate}
\end{Proposition}

\begin{proof}
This is a combination of \autoref{T:SandwichMain} and
the results from this section. For example, the explicit
parametrization of the simple modules for fixed apexes follows from \autoref{L:BasesAffineAndFiniteParts}.
\end{proof}

Explicitly identifying the sets $\hPartsgood$
or $\Partsgood$ of apexes is nontrivial, and we do not consider this in this paper.


\section{Proof of cellularity}\label{S:Proofs}
\stepcounter{subsection}

We are now ready to prove \autoref{T:BasesMain}, which shows that the weighted KLRW algebras are sandwich cellular algebras.

\begin{Remark}\label{R:ProofsShifted}
As before, the combinatorics below is separated into
partition and shifted partition combinatorics. The former is very similar
to type $\cone$, which was covered in \cite[Section 7]{MaTu-klrw-algebras}, so we focus only on the shifted partition case.
\end{Remark}

As in \cite[Section 7E]{MaTu-klrw-algebras} the most important
notion that we need is that of Young equivalence.
To define it we need some preliminary definitions.

\begin{Definition}\label{D:ProofsCloseQuadruple}
For $i,j\in I$, a \emph{close $(i,i,i)$-triple}, respectively
a \emph{close $(i,j,k)$-triple} or
a \emph{close $(i,j,j,i)$-quadruple},
is a collection of close strings as in the following local configurations:
\begin{gather}\label{Eq:ProofTripleQuadruple}
\text{triple:}
\begin{tikzpicture}[anchorbase,smallnodes,rounded corners]
\draw[solid](0,0)node[below]{$i$}--++(0,1);
\draw[solid](0.1,0)node[below]{$i$}--++(0,1)node[above,yshift=-1pt]{$\phantom{i}$};
\draw[solid](0.2,0)node[below]{$i$}--++(0,1);
\end{tikzpicture}
,\quad
\text{triple:}
\left(
\begin{tikzpicture}[anchorbase,smallnodes,rounded corners]
\draw[ghost](0,0)--++(0,1)node[above,yshift=-1pt]{$i$};
\draw[solid](0.1,0)node[below]{$j$}--++(0,1)node[above,yshift=-1pt]{$\phantom{i}$};
\draw[ghost](0.2,0)--++(0,1)node[above,yshift=-1pt]{$k$};
\end{tikzpicture}
\quad\text{or}\quad
\begin{tikzpicture}[anchorbase,smallnodes,rounded corners]
\draw[solid](0,0)node[below]{$i$}--++(0,1);
\draw[ghost](0.1,0)--++(0,1)node[above,yshift=-1pt]{$j$};
\draw[solid](0.2,0)node[below]{$k$}--++(0,1);
\end{tikzpicture}
\right)
\quad\text{and}\quad
\left(
\text{quadruple:}
\begin{tikzpicture}[anchorbase,smallnodes,rounded corners]
\draw[ghost](0.9,0)--++(0,1)node[above,yshift=-1pt]{$i$};
\draw[ghost](1.2,0)--++(0,1)node[above,yshift=-1pt]{$i$};
\draw[solid](1,0)node[below]{$j$}--++(0,1);
\draw[solid](1.1,0)node[below]{$j$}--++(0,1);
\end{tikzpicture}
\quad\text{or}\quad
\begin{tikzpicture}[anchorbase,smallnodes,rounded corners]
\draw[ghost](1,0)--++(0,1)node[above,yshift=-1pt]{$j$};
\draw[ghost](1.1,0)--++(0,1)node[above,yshift=-1pt]{$j$};
\draw[solid](0.9,0)node[below]{$i$}--++(0,1);
\draw[solid](1.2,0)node[below]{$i$}--++(0,1);
\end{tikzpicture}
\right)
.
\end{gather}
\end{Definition}

We also need the following, which should be compared with
\autoref{R:BasesCombinatoricsDottedIdempotent}. Here we
consider the two ghost $1$-strings in type $\dtwo$ as one string.

\begin{Definition}\label{D:ProofsRowEquivalence}
Let $S$ be a dotted straight line diagram.
Solid $i$ and $j$-strings of $S$ are \emph{pseudo row equivalent} if
either:
\begin{enumerate}

\item $i\rightsquigarrow j$, there are no dots on the $i$ or $j$-strings, and the ghost $i$-string is close and to the left of the solid $j$-string;

\item $i\leftsquigarrow j$, there are no dots on the $i$ or $j$-strings, and the solid $i$-string is close and to the left of the ghost $j$-string;

\item $i=j$ is a multisink, the $i$-string carries a dot, and the solid $i$-string is close and to the right of the solid $j$-string;

\end{enumerate}
A \emph{row equivalence class} is a pseudo row equivalence class in which there are no close $(i,i,i)$-triples,
the only close $(i,j,k)$-triples are of the form
$(i,i+1,i+2)$ or $(i,i-1,i-2)$, or
$(1,0,1)$ in type $\atwo$, and if $(i,j,j,i)$ is a close quadruple,
then $j$ is a multisink
and either $i\rightsquigarrow j$ or $i\leftsquigarrow j$.
\end{Definition}

The illustrations for parts (a)--(c) of \autoref{D:ProofsRowEquivalence} are:
\begin{gather*}
\text{(a)}:
\begin{tikzpicture}[anchorbase,smallnodes,rounded corners]
\draw[ghost](1,0)node[below]{$\phantom{i}$}--++(0,1)node[above,yshift=-0.05cm]{$i$};
\draw[ghost](1.9,0)node[below]{$\phantom{i}$}--++(0,1)node[above,yshift=-0.05cm]{$j$};
\draw[solid](0,0)node[below]{$i$}--++(0,1)node[above,yshift=-1pt]{$\phantom{i}$};
\draw[solid](0.9,0)node[below]{$j$}--++(0,1)node[above,yshift=-1pt]{$\phantom{i}$};
\end{tikzpicture}
,\quad
\text{(b)}:
\begin{tikzpicture}[anchorbase,smallnodes,rounded corners]
\draw[ghost](-0.1,0)node[below]{$\phantom{i}$}--++(0,1)node[above,yshift=-0.05cm]{$j$};
\draw[ghost](1,0)node[below]{$\phantom{i}$}--++(0,1)node[above,yshift=-0.05cm]{$i$};
\draw[solid](-1.1,0)node[below]{$j$}--++(0,1)node[above,yshift=-1pt]{$\phantom{i}$};
\draw[solid](0,0)node[below]{$i$}--++(0,1)node[above,yshift=-1pt]{$\phantom{i}$};
\end{tikzpicture}
,
\quad
\text{(c)}:
\begin{tikzpicture}[anchorbase,smallnodes,rounded corners]
\draw[ghost](1,0)node[below]{$\phantom{i}$}--++(0,1)node[above,yshift=-0.05cm]{$i$};
\draw[ghost,dot](1.2,0)node[below]{$\phantom{i}$}--++(0,1)node[above,yshift=-0.05cm]{$i$};
\draw[solid](0,0)node[below]{$i$}--++(0,1)node[above,yshift=-1pt]{$\phantom{i}$};
\draw[solid,dot](0.2,0)node[below]{$i$}--++(0,1)node[above,yshift=-1pt]{$\phantom{i}$};
\end{tikzpicture}
.
\end{gather*}
These should be compared with \autoref{Eq:BasesImportantConfigurationOne}
and \autoref{Eq:BasesImportantConfigurationTwo}. Note that \autoref{Eq:BasesImportantConfigurationOneDTwo} is also included in the description since pseudo row equivalence classes only consider two solid strings at a time.

By definition, there is a unique \emph{first} and \emph{last} string in each row equivalence class because the strings in a row equivalence class are ordered by closeness, starting with the string that is not close and to the right of any other in the same class.

\begin{Definition}\label{D:ProofsYoung}
Assume that we in the case of shifted partitions.
A \emph{Young equivalence class} $Y$ is a disjoint union of row equivalence classes $R_{1}\cup\dots\cup R_{z}$ such that:
\begin{enumerate}

\item The first string in $R_{1}$ has no dot and is close to an (affine) red string of the same residue;

\item $|R_{1}|>|R_{2}|>\dots>|R_{z}|$;

\item The first string in $R_{a+1}$ is an $i$-string that is close to a dotted solid $i$-string of the same residue in $R_{a}$, or there is a $j$-string in $R_{a}$ that satisfies one of closeness conditions in (a) and (b) of \autoref{D:ProofsRowEquivalence}, with respect to this string.

\end{enumerate}
For partitions, we use the analog of \cite[Definition 7E.6]{MaTu-klrw-algebras}. That is, part (b) is replaced with $|R_{1}|\geq|R_{2}|\geq\dots\geq|R_{z}|$
and (c) mimics \autoref{Eq:BasesMainTableaux} with a dot on the $i$-string in the first two cases therein.
\end{Definition}

Recall that $L(S)$ is the left-justification of the dotted straight line diagram $S$ as {\eg} in
\autoref{E:RecollectionLeftJustification}.

\begin{Lemma}\label{L:ProofsVerticalHell}
Let $S$ be a dotted straight line diagram.
Then $L(S)=L(\dotidem)$, for some
$\blam\in\hParts$, if and only if the solid strings of $S$ are a disjoint union of Young equivalence classes.
\end{Lemma}

\begin{proof}
By construction, the solid strings in $\dotidem$
are a disjoint union of Young equivalence classes {\cf}
\autoref{R:BasesCombinatoricsDottedIdempotent}.

To prove the converse, given a dotted straight line diagram $S$,
we construct an $\hell$-partition $\blam$ by inductively associating the solid
strings in each Young equivalence class $Y$ to nodes in a component $\lambda^{(m)}$ of $\blam$.
We explain the situation of shifted partitions, which is
when $\rho_{m}$ is a multisink. The case of partitions can be proven similarly, following \cite[Lemma 7E.8]{MaTu-klrw-algebras}.

By \autoref{D:ProofsYoung}.(a), the first string of $Y$ is left adjacent to an (affine) red string of the same residue. If this is the $m$th red string, then identify the solid string with the node $(m,1,1)$. By induction we now assume that the $k$th solid $i$-string in $Y$ corresponds to the node $(m,r,c)\in\blam$ and consider the $(k+1)$st solid $j$-string. There are two cases to consider.

\textit{Case 1.} First, if $i$ is not the last
string in its row equivalence class, then
\autoref{Eq:BasesImportantConfigurationOne}, \autoref{Eq:BasesImportantConfigurationOneDTwo} and
\autoref{Eq:BasesImportantConfigurationTwo} correspond to
(a)--(c) of \autoref{D:ProofsRowEquivalence} and the condition on
close $(i,j,k)$-triples and close $(j,i,i,j)$-quadruples, with the
correct dot placement. Moreover, no other
configurations can appear, {\ie} there are no close $(i,i,i)$-triples
or any other triples or quadruples due to (a)--(c) of \autoref{D:ProofsRowEquivalence}.
Hence, the $(k+1)$st solid $j$-string corresponds to the node $(m,r,c+1)$.

\textit{Case 2.} If on the other hand $i$ is the last
string in its row equivalence class, then
we observe that \autoref{Eq:BasesMainTableaux}
corresponds to \autoref{D:ProofsYoung}.(c), and
the $(k+1)$st solid $j$-string corresponds to the node $(m,r+1,c)$.

Finally, note that the condition in \autoref{D:ProofsYoung}.(b)
ensures that the resulting diagram is a shifted $\hell$-partition.
\end{proof}

\begin{Proposition}\label{P:ProofsVerticalDominance}
Suppose that $D\in\WA[n](X)$ and that $D$ factors through the dotted idempotent diagram~$S$. Then there exists $\blam\in\hParts$ such that $D$ factors through $\dotidem$ and $\blam\gedom L(S)$.
\end{Proposition}

\begin{proof}
Without loss of generality we can assume that $D=S$. If $D=\dotidem$,
then there is nothing to prove by \autoref{L:ProofsVerticalHell}.
So assume that \autoref{L:ProofsVerticalHell} is not satisfied, so that
that $D$ is not a disjoint union of Young equivalence classes.

We can assume that all strings are
within $[\min X,\max X+1]\times[0,1]$, the \emph{region defined by $X$}.
Let $s$ be the rightmost solid string in $D$ that is not in any
Young equivalence class. We want to argue that we can pull $s$ to the right, jump dots on $s$ further to the right, or we can attach $s$ to a Young equivalence class, which implies
the claim by induction.
There are a few cases which we need to discuss.
We only consider shifted partition combinatorics since
the arguments for usual partition combinatorics are {\muta} as in \cite[Proposition 7E.9]{MaTu-klrw-algebras}.

\textit{Case 1a.} $s$ is the rightmost string in the sense that we can pull
it arbitrarily far to the right, and $s$ does not have a dot. In this case we can park it next to an affine red string of the same residue, and it is now part of a Young equivalence
class by \autoref{D:ProofsYoung}.(a).

\textit{Case 1b.} As in Case 1a, but now $s$ carries a dot. After pulling
the dot to the top of the diagram we are back in Case 1a.

We now assume that we are not in Cases 1a and 1b. Then, up to isotopy,
$s$ or its ghost
is close and to the right of a solid, ghost
or red string $t$. We focus on the situation when $s$ is close to $t$, where we again have several cases. The cases
where the ghost of $s$ is close to $t$ follow {\muta} and are omitted.
We also assume that $s$ does not carry a dot. If it does, then there is an additional extra argument one needs, as is explained in \cite[Proof of Proposition 7E.9, Case 5]{MaTu-klrw-algebras}
(this argument works {\muta} in the $\bad$ types), but in the end
the relation used below allow us to continue with the induction.

\textit{Case 2a.} Assume that $s$ is not in the Young equivalence
class of $t$ because $s$ does not satisfy the close $(i,j,k)$-triple condition because $s$ is the
leftmost string in \autoref{Eq:ProofTripleQuadruple}. In this case the
right-hand relation in \autoref{Eq:StrategyReidemeisterIIAndIIITypeA} applies, so we can pull $s$ further to the right.
(This works unless we in a close $(1,0,1)$-triple situation in type $\atwo$, in which case $s$ and $t$ are in the same Young equivalence class by \autoref{D:ProofsRowEquivalence}.)

\textit{Case 2b.} Similarly, assume that $s$ is not in the Young equivalence
class of $t$ because $s$ does not satisfy the close $(i,j,j,i)$-quadruple condition 
because $s$ is
the leftmost string in \autoref{Eq:ProofTripleQuadruple}.
Then we can use \autoref{L:StrategyQuadruple} to pull
$s$ further to the right.

\textit{Case 2c.} We now assume that $s$ is not
in the Young equivalence class of $t$ because the condition
$|R_{1}|>|R_{2}|>\dots>|R_{z}|$ is not satisfied.
For $i,j\in I$ and $i\rightsquigarrow j$ or $i\leftsquigarrow j$,
the crucial configurations are
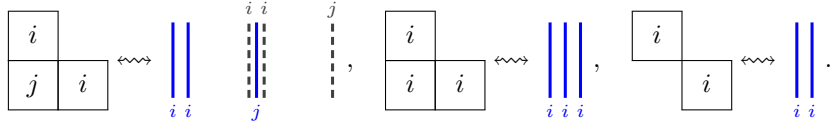
\begin{gather*}
\begin{ytableau}
i & \none \\ j & i
\end{ytableau}
\leftrightsquigarrow
\begin{tikzpicture}[anchorbase,smallnodes,rounded corners]
\draw[ghost](0.8,0)node[below]{$\phantom{i}$}--++(0,1)node[above,yshift=-0.05cm]{$i$};
\draw[ghost](1,0)node[below]{$\phantom{i}$}--++(0,1)node[above,yshift=-0.05cm]{$i$};
\draw[ghost](1.9,0)node[below]{$\phantom{i}$}--++(0,1)node[above,yshift=-0.05cm]{$j$};
\draw[solid](-0.2,0)node[below]{$i$}--++(0,1)node[above,yshift=-1pt]{$\phantom{i}$};
\draw[solid](0,0)node[below]{$i$}--++(0,1)node[above,yshift=-1pt]{$\phantom{i}$};
\draw[solid](0.9,0)node[below]{$j$}--++(0,1)node[above,yshift=-1pt]{$\phantom{i}$};
\end{tikzpicture}
,\quad
\begin{ytableau}
i & \none \\ i & i
\end{ytableau}
\leftrightsquigarrow
\begin{tikzpicture}[anchorbase,smallnodes,rounded corners]
\draw[solid](-0.2,0)node[below]{$i$}--++(0,1)node[above,yshift=-1pt]{$\phantom{i}$};
\draw[solid](0,0)node[below]{$i$}--++(0,1)node[above,yshift=-1pt]{$\phantom{i}$};
\draw[solid](0.2,0)node[below]{$i$}--++(0,1)node[above,yshift=-1pt]{$\phantom{i}$};
\end{tikzpicture}
,\quad
\begin{ytableau}
i & \none \\ \none & i
\end{ytableau}
\leftrightsquigarrow
\begin{tikzpicture}[anchorbase,smallnodes,rounded corners]
\draw[solid](-0.2,0)node[below]{$i$}--++(0,1)node[above,yshift=-1pt]{$\phantom{i}$};
\draw[solid](0,0)node[below]{$i$}--++(0,1)node[above,yshift=-1pt]{$\phantom{i}$};
\end{tikzpicture}
.
\end{gather*}
There are a few cases, but for all
these we can use \autoref{Eq:StrategyTwoSolid} or \autoref{Eq:StrategyReidemeisterIIAndIIITypeA}
to pull $s$ to the right.

Assume now that we are not in any of the cases above.

\textit{Case 3a.} If $t$ is an (affine) red string, then
an honest Reidemeister II relation pulls $s$ further to the
right. We can apply such a relation since the case where
$s$ has the same residue as $t$ is, which case $s$ is the first string in a row equivalence class.

\textit{Case 3b.} If $t$ is a solid string, then
an honest Reidemeister II relation applies unless $s$ has the same residue as $t$. In this latter case \autoref{Eq:StrategyTwoSolid} applies.

\textit{Case 3c.} Finally, if $t$ is a ghost string, then
we can use an honest Reidemeister II relation to pull $s$ to the right.
An honest Reidemeister II relation applies because the assumption that $s$ and $t$ are not in the same
Young equivalence class implies that $s$ and $t$ do not satisfy the conditions of (a) and (b) of \autoref{D:ProofsRowEquivalence}, or \autoref{D:ProofsYoung}.(c).

Hence, the result follows by induction.
\end{proof}

The rest of the proof of \autoref{T:BasesMain} is essentially
the same as in \cite[Section 7E]{MaTu-klrw-algebras}.
That is, applying dots or crossings to $\dotidem$
gives a linear combination of bigger elements. To this end, recall the definition of the finite dots and the integers $c_{m}(\blam)$ from \autoref{D:BasesSandwichpart}.

\begin{Lemma}\label{L:ProofsPullingDotsC}
Suppose that $\blam\in\Parts$ and $1\leq m\leq n$. Then $y_{m}y_{m}^{c_{m}(\blam)}\dotidem\in\WAlam*$.
\end{Lemma}

\begin{proof}
By \autoref{P:ProofsVerticalDominance}, the diagram $y_{m}y_{m}^{c_{m}(\blam)}\dotidem$ factors through $\dotidem[\bmu]$, for some $\bmu\gedom\blam$.
\end{proof}

By ignoring the components, consider $\blam$ as a composition, and let $\Sym[\blam]$ be the associated Young subgroup of $\Sym[n]$.

\begin{Lemma}\label{L:ProofsCosetRepsC}
Suppose that $\blam\in\hParts$ and $w\in\Sym[\blam]$. Then $D_{\blam}(w)\idem,\idem D_{\blam}(w)\in\WAlam*$.
\end{Lemma}

\begin{proof}
As in \cite[Lemma 6D.17]{MaTu-klrw-algebras}.
\end{proof}

\begin{proof}[Proof of \autoref{T:BasesMain}]
The arguments given in \cite[Sections 6D and 7E]{MaTu-klrw-algebras}, which are the
analogous statements for the $\ac$ types,
apply in $\bad$ types as well. In fact, these arguments are general and
use only \autoref{P:RecollectionWABasis} and
\autoref{R:RecollectionWABasis}, as well as the analogs
of the results proven above.
\end{proof}

\begin{Remark}\label{R:ProofGeneral}
Following the strategy in \autoref{R:StrategyMain},
the construction of the homogeneous (affine)
sandwich cellular basis and proof of \autoref{T:BasesMain} splits
into several parts:
\begin{enumerate}

\item Because the bottom and top of the cellular basis elements are be given
by permutation diagrams, the first step is
to find suitable tableau combinatorics for the quiver under study.
For a general quiver this is potentially hopeless, and one might want to use some other combinatorial data instead of tableaux, but for
a lot of quivers an answer is already in the literature.

\item The construction of the middle in \autoref{R:SandwichCellularAlgebra} is then crucial. This part is noncanonical,
although mostly dictated by \autoref{R:StrategyMain}. We hope to explain a more general approach in future work. Note that
additional dots might be
necessary to prevent basis elements being annihilated by \autoref{Eq:RecollectionReidemeisterII} and to ensure that we have the analog of \autoref{L:BasesCombinatoricsIdempotent}.

\item From here onwards the arguments are general
and do not depend on the quiver anymore:
\autoref{P:ProofsVerticalDominance} follows by analyzing
the combinatorics of the string placement of $\idem$, and this
proposition in turn directly implies \autoref{L:ProofsPullingDotsC}
and \autoref{L:ProofsCosetRepsC}. Once these two lemmas have been established
the proof of cellularity \autoref{T:BasesMain} is formal.
Linear independence follows using the faithful polynomial module in \autoref{R:RecollectionWABasis}, spanning using \autoref{L:ProofsPullingDotsC}
and \autoref{L:ProofsCosetRepsC} and the standard basis
in \autoref{P:RecollectionWABasis}. The latter arguments are independent of the underlying quiver.

\end{enumerate}
\end{Remark}


\section{Table of notation and central concepts}

In general we use an overline for notation that only plays
a role in the infinite dimensional setting.

\begin{xltabular}{\textwidth}{p{41mm}lX}\toprule
Name & Symbol & Description
\\
\midrule\endhead
\bottomrule\endfoot
Solid, ghost, red and affine strings
&
$\begin{tikzpicture}[anchorbase,smallnodes,rounded corners]
\draw[solid] (0,0)node[below]{$i$} to (0,0.5)node[above,yshift=-1pt]{$\phantom{i}$};
\end{tikzpicture}$
,
$\begin{tikzpicture}[anchorbase,smallnodes,rounded corners]
\draw[ghost] (0,0)node[below]{$\phantom{i}$} to (0,0.5)node[above,yshift=-1pt]{$i$};
\end{tikzpicture}$
,
$\begin{tikzpicture}[anchorbase,smallnodes,rounded corners]
\draw[redstring] (0,0)node[below]{$i$} to (0,0.5)node[above,yshift=-1pt]{$\phantom{i}$};
\end{tikzpicture}$
,
$\begin{tikzpicture}[anchorbase,smallnodes,rounded corners]
\draw[affine] (0,0)node[below]{$i$} to (0,0.5)node[above,yshift=-1pt]{$\phantom{i}$};
\end{tikzpicture}$
&
Strings in diagrams, see
\autoref{Eq:RecollectionStrings}
\\
\hline\mystrut
$\ell$=level, $\charge$=charge
& $n$, $\ell$, $e$, $\charge$, $\brho$, $X$
& The number of solid strings, the number of red strings, the number of vertices in the quiver, the positions of the red strings, the labels of the red strings, the positions and labels of the solid strings; all of these are fixed from the start, see \autoref{SS:RecollectionDefinition}
\\
Affine notation
&
$\hell$, $\affine{\charge}$, $\affine{\brho}$
&
We set $\hell=\ell+n(e+1)$. These ensure that we have enough affine red string to catch every solid string, see \autoref{SS:BasesSomeNotation}
\\
\hline\mystrut
Weighted KLRW algebra
& $\WA[n](X)$
& A diagram algebra with $n$ solid strings with top and bottom $x$-coordinates in $X$, see \autoref{SS:RecollectionDefinition}
\\
The cyclotomic version
& $\WAc[n](X)$
& A quotient of $\WA[n](X)$, {\cf}
\autoref{Eq:RecollectionUnsteady}
\\
\hline\mystrut
Idempotent diagram for $\blam$
&
$\idem$
&
The idempotent diagram created by inductively placing strings
while reading along the rows of
the $\affine{\brho}$-partition $\blam$, see \autoref{SS:BasesPositioning}
\\
The dotted version of $\idem$
&
$\dotidem$
&
We place dots $y_{\blam}$ on $\idem$, see \autoref{SS:BasesPositioningDot}
\\
Affine and finite part
&
$y^{\ba}$ and $y^{\bfi}$
&
These determine the sandwiched part of the main basis; the affine part only plays a role for $\WA[n](X)$, while the finite part appears for both
$\WA[n](X)$ and $\WAc[n](X)$, see \autoref{SS:BasesSandwichPart}
\\
Permutation diagram for $\bS$
&
$D_{\bS}$
&
A permutation diagram associated to
the $\affine{\brho}$-tableaux $\bS$, see
\autoref{SS:BasesPermuations}
\\
Sandwiched algebras
&
$\sand[\blam]$
&
The sandwiched algebras are tensor products of $R[X]$ and $R[X]/(X^{2})$, see
\autoref{SS:BasesConstruction}
\\
Cellular basis elements
&
$D_{\bS\bT}^{\ba,\bfi}$, $D_{\bS\bT}^{\bfi}$
&
The homogeneous cellular basis elements, see \autoref{SS:BasesConstruction}, for $\WA[n](X)$ and $\WAc[n](X)$
\\
Cellular basis sets
&
$\BX$, $\BXc$
&
The homogeneous cellular bases sets, see \autoref{SS:BasesConstruction}, for $\WA[n](X)$ and $\WAc[n](X)$
\\
\hline\mystrut
Cellular bases for the KLR(W) algebras
&
$\TA[n]$, $\TAc[n]$, $\EX$, $\EXc$
&
The KLR(W) algebras $\TA[n]$, $\TAc[n]$ arise via idempotent truncation
and $\EX$, $\EXc$
are the corresponding homogeneous cellular bases
\\
\hline\mystrut
Partitions
&
$\hParts$, $\Parts$
&
Indexing sets for the middle of the cellular bases of $\WA[n](X)$ and $\WAc[n](X)$, see \autoref{SS:BasesCombinatorics}
\\
Nodes
&
$(m,r,c)$
&
Nodes in (shifted) partitions: $m$ is the component index, $r$ the row index and $c$ the column index, see \autoref{R:BasesOrderOfMRC}
\\
Simple modules
&
$\hPartsgood$, $\Partsgood$
&
Indexing sets for the simple modules of $\WA[n](X)$ and $\WAc[n](X)$, see \autoref{SS:BasesConstruction}
\\
Tableaux
&
$\SStd$, $\Std$
&
Indexing sets for the bottom/top of the cellular bases of $\WA[n](X)$ and $\TA[n]$, {\cf} \autoref{SS:BasesBasisDiagrams}
\\
Positioning function
&
$\hcoord$
&
The function that assigns $x$-coordinates to nodes of Young diagrams, see \autoref{SS:BasesPositioning}
\\
Maximal finite part
&
$\max_{\R}(fin)$
&
Determines where the affine part of the diagrams begin, see \autoref{SS:BasesPositioning}
\end{xltabular}
\label{table:notation}

\end{document}